\documentclass[12pt]{amsart}
\usepackage{amssymb,amsfonts,amsthm,latexsym,bm}
\usepackage{hyperref}
\usepackage[centertags]{amsmath}
\usepackage{mathrsfs}
\usepackage{t1enc}
\usepackage{graphicx}
\usepackage{enumerate}
\usepackage{enumitem}
\setlist{leftmargin=0.5cm}

\usepackage{etoolbox}

\advance\headheight by 3pt

\patchcmd{\abstract}
{\leftmargin3pc \rightmargin\leftmargin}%
{\leftmargin1pc \rightmargin1pc}%
{}{}

\newtheorem{theorem}{Theorem}[section]
\newtheorem{corollary}[theorem]{Corollary}
\newtheorem{lemma}[theorem]{Lemma}

\newtheorem*{theorem*}{Theorem}
\newtheorem{proposition}[theorem]{Proposition}
\newtheorem{conjecture}[theorem]{Conjecture}
\newtheorem{prob}{Problem}

\newtheorem{innercustomthm}{Theorem}
\newenvironment{customthm}[1]
{\renewcommand\theinnercustomthm{#1}\innercustomthm}
{\endinnercustomthm}

\newtheorem{innercustomcor}{Corollary}
\newenvironment{customcor}[1]
{\renewcommand\theinnercustomcor{#1}\innercustomcor}
{\endinnercustomcor}

\numberwithin{equation}{section}

\def\pitem{\advance\leftskip3mm\advance\linewidth-3mm}
\def\mitem{\advance\leftskip-3mm\advance\linewidth3mm}

\begingroup\makeatletter\ifx\SetFigFont\undefined
\gdef\SetFigFont#1#2#3#4#5{
  \reset@font\fontsize{#1}{#2pt}
  \fontfamily{#3}\fontseries{#4}\fontshape{#5}
  \selectfont}
\fi\endgroup

\theoremstyle{definition}

\newtheorem*{rem}{Remark}

\renewcommand{\subjclass}[1]{\thanks{\emph{2020 Mathematics Subject Classification:}~#1}}
\renewcommand{\keywords}[1]{\thanks{\emph{Keywords and Phrases:}~#1}}
\renewcommand{\date}{\thanks{\today}}

\newcommand{\tv}{{\bf t}}

\newcommand{\CC}{\mathcal{C}}
\newcommand{\DD}{\mathcal{D}}

\newcommand{\II}{\mathcal{I}}

\newcommand{\MM}{\mathcal{M}}

\newcommand{\OO}{\mathcal{O}}
\newcommand{\PP}{\mathcal{P}}

\newcommand{\RR}{\mathcal{R}}
\renewcommand{\SS}{\mathcal{S}}
\newcommand{\TT}{\mathcal{T}}

\newcommand{\VV}{\mathcal{V}}

\newcommand{\fp}{\mathfrak{p}}

\newcommand{\fa}{\mathfrak{a}}

\newcommand{\fc}{\mathfrak{c}}
\newcommand{\fd}{\mathfrak{d}}

\newcommand{\fD}{\mathfrak{D}}

\newcommand{\fI}{\mathfrak{I}}

\newcommand{\Cc}{\mathbb{C}}
\newcommand{\Rr}{\mathbb{R}}
\newcommand{\Qq}{\mathbb{Q}}

\newcommand{\Zz}{\mathbb{Z}}

\newcommand{\Pp}{\mathbb{P}}
\newcommand{\Gg}{\mathbb{G}}

\newcommand{\vk}{\Bbbk}

\renewcommand{\ge}{\geq}
\renewcommand{\le}{\leq}

\newcommand{\ve}{\varepsilon}

\def\house#1{\setbox1=\hbox{$\,#1\,$}%
\dimen1=\ht1 \advance\dimen1 by 2pt \dimen2=\dp1 \advance\dimen2
by 2pt
\setbox1=\hbox{\vrule height\dimen1 depth\dimen2\box1\vrule}%
\setbox1=\vbox{\hrule\box1}%
\advance\dimen1 by .4pt \ht1=\dimen1 \advance\dimen2 by .4pt
\dp1=\dimen2 \box1\relax}

\newcommand{\kdots}{,\ldots ,}

\newcommand{\GL}{{\rm GL}}

\newcommand{\half}{\mbox{$\textstyle{\frac{1}{2}}$}}

\newcommand{\medfrac}[2]{\mbox{\large{$\textstyle{\frac{#1}{#2}}$}}}

\renewcommand{\gcd}{{\rm gcd}}

\newcommand{\rank}{{\rm rank}\,}

\newcommand{\cross}{{\rm cr}}
\newcommand{\Ht}{{\rm Ht}}

\DeclareRobustCommand{\SkipTocEntry}[5]{}

\title[Effective reduction theory of integral polynomials]{Effective reduction theory of integral polynomials of
given discriminant, and related topics 
\\[1.5mm]{\MakeLowercase{\textnormal{(survey with a brief historical overview)}}}
}
\dedicatory{\vskip-0.4cm{Dedicated to Charles Hermite (1822-1901) and Alan Baker (1939-2018), for their fundamental contributions to the field}}


\author[J.-H. Evertse]{Jan-Hendrik Evertse}
\address{J.-H. Evertse \newline
         \indent Universiteit Leiden, Mathematisch Instituut, \newline
         \indent Postbus 9512, 2300 RA Leiden, The Netherlands \newline
         \indent \textit{URL:} {\tt https://pub.math.leidenuniv.nl/$\sim$evertsejh}}
\email{evertse\char'100math.leidenuniv.nl}

\author[K. Gy\H ory]{K\'alm\'an Gy\H ory}
\address{K. Gy\H ory\newline
	\indent University of Debrecen, Institute of Mathematics, \newline
	\indent P.O. Box 400, 4002 Debrecen, Hungary \newline
	\indent \textit{URL:} {\tt https://math.unideb.hu/en/dr-kalman-gyory}}
\email{gyory\char'100science.unideb.hu}

\begin{document}

\begin{abstract}
For polynomials in $\Zz [X]$, the classical $\Zz$-equivalence (monic case) and $GL_2(\Zz)$-equivalence preserve the discriminant as an invariant.
The effective  reduction theory for polynomials of given degree and discriminant 
consists of results that give,
for a given polynomial $f\in \Zz [X]$, a $\Zz$-equivalent (in the monic case) or $GL_2(\Zz)$-equivalent polynomial $g$ whose coefficients are effectively bounded above in terms of only the degree and discriminant of $f$.
We discuss the classical results of this type of Lagrange (1773) and Hermite (1851) on quadratic and cubic polynomials, the general ineffective theorem of Birch and Merriman (1972),
 the general effective theorem of Gy\H{o}ry (1973) for monic polynomials, obtained independently, and that of Evertse and Gy\H{o}ry (1991) for arbitrary polynomials.
The proofs of these two effective theorems use Gy\H{o}ry's effective results on unit equations, which were proved by means of Baker's effective theory of logarithmic forms.
Later Evertse, Gy\H ory and others obtained several applications and generalizations; see the book Evertse and Gy\H ory (2017).
In his long-forgotten paper Hermite (1857), Hermite attempted  to extend the above results of Lagrange and Hermite  to polynomials of arbitrary degree. However, as was pointed out in our joint work BEGyRS (2023) with Bhargava, Remete and Swaminathan, Hermite (1857) proved an important result but with a weaker equivalence only. Thus, it was only by the above mentioned theorems of Győry (1973) and Evertse and Győry (1991) that Hermite's problem from 1857 was settled in full generality.
This and many other recent results inspired us to write this survey paper on the subject. We present here several older and recent generalizations and applications of the effective reduction theory, e.g., to monogenic number fields and monogenic and rationally monogenic orders.  
We also give an overview of bounds on the number of times a given order is monogenic or rationally monogenic. 
In the Appendix we discuss further related topics not strictly belonging to the reduction theory of integral polynomials.
\end{abstract}
    
\maketitle



\newpage
\setcounter{tocdepth}{1}
\tableofcontents

\newpage
\section{Introduction}\label{sec:0}

We give an overview of older and recent results on the reduction theory of integral polynomials of given discriminant, and its many consequences and applications. We first recall some definitions and notation.

\subsection{Preliminaries}\label{sec:0.1}
\vphantom{a}

\noindent
Two polynomials $f,g\in \Zz [X]$ of degree $n\ge 2$ are called \emph{$\Zz$-equivalent} if
\[
\text{$g(X)=f(X+a)$ or $g(X)=(-1)^nf(-X+a)$ for some $a\in\Zz$,}
\] 
and \emph{$GL_2(\Zz )$-equivalent}
if 
\[
\text{$g(X)=\pm (cX+d)^ nf\big(\medfrac{aX+b}{cX+d}\big)$ for some matrix $\big(\begin{smallmatrix}a&b\\c&d\end{smallmatrix}\big)\in GL_2(\Zz )$,}
\]
 i.e., $a,b,c,d\in\Zz$ and $ad-bc=\pm 1$. Clearly, $\Zz$-equivalence implies $GL_2(\Zz )$-equivalence. 
Polynomials that are $\Zz$-equiv\-alent to a monic polynomial 
are also monic.

The \emph{discriminant} of a polynomial 
$$
f=a_0X^n+\cdots +a_n=a_0\prod_{i=1}^n(X-\alpha_i),\ \ \text{with } a_0\not= 0
$$ 
is defined by 
$$
D(f):=a_0^{2n-2}\prod_{1\leq i<j\leq n} (\alpha_i-\alpha_j)^2.
$$
This is a homogeneous polynomial of degree $2n-2$ in $\Zz [a_0\kdots a_n]$;
thus, if $f\in \Zz [X]$ then $D(f)\in\Zz$. As one may easily verify,
polynomials that are $\Zz$-equivalent or $GL_2(\Zz )$-equivalent 
have the same discriminant.

We define the \emph{height} $H(f)$ of a polynomial
$f=a_0X^n+\cdots +a_n\in \Zz [X]$ by 
$$
H(f):=\max (|a_0|\kdots |a_n|).
$$

An \emph{invariant} is a function $\Zz [X]\to\Rr$ that assumes the same value at $GL_2(\Zz )$-equivalent polynomials.
In general, \textit{reduction theory of polynomials} is about results 
of the following type: given a set of invariants, $I_1\kdots I_t$, say, there exists for any $f\in \Zz [X]$
a polynomial $g\in \Zz [X]$ that is 
$GL_2(\Zz )$-equivalent (or $\Zz$-equivalent in the monic case) to $f$ and whose coefficients are bounded in terms of $I_1(f)\kdots I_t(f)$.
In this paper, we focus on results in which the height $H(g)$ of $g$ is bounded above in terms of $\deg f$ and $|D(f)|$, i.e., on \textit{reduction theory for polynomials of given degree and given discriminant}. Such results imply that up to $GL_2(\Zz )$-equivalence (resp. $\Zz$-equivalence if we restrict ourselves to monic polynomials) there are only finitely many polynomials $f\in \Zz [X]$ of degree $n$ and given discriminant $D\not= 0$.

In fact, most of the literature deals with reduction theory of \emph{binary forms} of given discriminant.
Recall that any binary form $F(X,Y)\in\Zz [X,Y]$ can be factored as $\prod_{i=1}^n (\alpha_iX-\beta_iY)$ with algebraic $\alpha_i ,\beta_i$, 
and that its discriminant is $D(F):=\prod_{1\leq i<j\leq n}(\alpha_i\beta_j-\alpha_j\beta_i)^2$. Two binary forms $F,G\in\Zz [X,Y]$ are called $GL_2(\Zz )$-equivalent if 
$G(X,Y)=\pm F(aX+bY,cX+dY)$ for some $\big(\begin{smallmatrix}a&b\\c&d\end{smallmatrix}\big)\in GL_2(\Zz )$, and clearly,
$GL_2(\Zz )$-equivalent binary forms have the same discriminant.
The results on reduction theory for binary forms $F$
can be translated immediately into similar results for univariate polynomials $f$ and vice-versa, using the correspondence $f(X)=F(X,1)$, $F(X,Y)=Y^{\deg f}f(X/Y)$. As in our joint paper BEGyRS (2023) with Bhargava, Remete and Swaminathan, we formulate our results in terms of univariate polynomials for convenience of presentation.

For definitions of \textit{effectively given} concepts, structures and \textit{effective determination, computation}, one can consult e.g. the corresponding sections of our books Evertse and Gy\H ory (2015, 2017, 2022).

\subsection{Summary}
\vphantom{a}

\noindent
Lagrange (1773) developed a
reduction theory of integral binary quadratic forms of given discriminant, which can be translated immediately into a reduction theory for quadratic polynomials of given non-zero discriminant.
His results imply that up to the classical $GL_2(\Zz )$-equivalence, resp. $\Zz$-equivalence (monic case) there are only finitely many quadratic polynomials in $\Zz [X]$ of given discriminant. Lagrange's result is \emph{effective} in the sense that one can
effectively determine the reduced polynomials. This was later made more precise by Gauss (1801).

Hermite (1848, 1851) introduced a reduction theory for binary forms, or equivalently univariate polynomials of arbitrary degree but using another invariant instead of the discriminant. In the case of cubic polynomials,
Hermite's invariant is up to a constant a power of the absolute value of the discriminant.
Thus, Hermite's reduction theory implies that up to $GL_2(\Zz )$-equivalence there are only finitely many \emph{cubic} polynomials  in $\Zz [X]$ of given discriminant.
Hermite was apparently interested to extend this to polynomials of arbitrary degree $n\geq 4$.
In Hermite (1857) he introduced a new equivalence relation  (called by us `\emph{Hermite equivalence}', see Section \ref{sec:2})  and proved in an 
ineffective way a finiteness result on the corresponding 
equivalence classes of integral polynomials of degree $n$ and discriminant $D$.
But he did not compare his equivalence relation to the classical
equivalence relations, i.e., to
$GL_2(\Zz )$-equivalence and $\Zz$-equivalence.
The result of Hermite (1857) does not appear to have been studied in the literature until the excellent book of Narkiewicz (2018), where Hermite equivalence was confused with the classical equivalence relations.

Hermite's apparent goal, i.e., the finiteness result with $GL_2(\Zz )$-equiv\-alence
instead of Hermite equivalence, was finally achieved
more than a century later by Birch and Merriman (1972) for arbitrary polynomials in an ineffective form 
and independently, for monic polynomials and in a more precise and effective form by Gy\H{o}ry (1973).
The general result of Birch and Merriman was subsequently made 
effective by Evertse and Gy\H{o}ry (1991a).
More precisely, 
Gy\H{o}ry (1973) and Evertse and Gy\H{o}ry (1991a) proved that there exists an effectively computable number $c(n,D)$ depending only on $n$ and $D$ such that every $f\in \Zz [X]$ of degree $n$ and discriminant $D\not= 0$ is $GL_2(\Zz )$-equivalent (and even $\Zz$-equivalent in the monic case) to a polynomial $g$ with height 
\begin{equation}\label{0}
H(g)\leq c(n,D).
\end{equation}
These results heavily depend on effective finiteness results for unit equations $ax+by=1$ with solutions $x,y$ from the unit group of the ring of integers of a number field, which were derived in turn using Baker's theory of logarithmic forms.
This solved the old problem of Hermite (1857)
mentioned above in an effective way,
and further resulted in many significant consequences
and applications.
For example, in the 1970's, Gy\H{o}ry  deduced 
from his paper from 1973 the first general effective algorithm
that decides monogenicity and existence of power integral bases of number fields,
and in fact finds all power integral bases.
For later applications and generalizations we refer to the monograph Evertse and Gy\H{o}ry (2017) and Sections \ref{sec:3}--\ref{sec:9} of the present paper.

In our recent paper BEGyRS (2023) with Bhargava, Remete and
Swaminathan we provided a thorough treatment of the notion of Hermite equivalence, and proved that $\Zz$-equivalence and $GL_2(\Zz )$-equival\-ence 
are much more precise
than Hermite equivalence. This confirmed that Hermite's result from 1857
was weaker than those of Birch and Merriman, Gy\H{o}ry,
and that of Evertse and Gy\H{o}ry mentioned above. It should of course be mentioned that unlike the last authors, Hermite didn't have the powerful Baker's theory of logarithmic forms and its application to unit equations  at his disposal.

In Section \ref{sec:1} we briefly recall the reduction theory of quadratic and cubic polynomials
of given non-zero discriminant. In Section \ref{sec:2}, following BEGyRS (2023), we deal with Hermite equivalence
and compare it with $\Zz$-equivalence and $GL_2(\Zz )$-equivalence.
In Section \ref{sec:3} we discuss in more detail the general results of Birch and Merriman (1972),
Gy\H{o}ry (1973), Evertse and Gy\H{o}ry (1991a), and those from the paper BEGyRS (2023).
We present the best known effective height estimates for the solutions of unit equations and $S$-unit equations.
We sketch how to deduce results of the type \eqref{0}.  
An important part of Section \ref{sec:3} is Subsection \ref{sec:4.4new}, which gives much stronger conjectural upper bounds for the height of $g$. These bounds follow from the abc-conjecture and related conjectures.
This is partly joint work with Rafael von K\"anel.
In Section \ref{sec:4} we present some consequences in algebraic number theory. In particular, 
we give an overview of effective finiteness results concerning algebraic numbers of given discriminant, resp. given index, and index form equations. 
Further, we deduce applications to monogenic number fields and orders,
and also generalizations to so-called rationally monogenic orders.
In Section \ref{sec:5} we discuss practical algorithms for solving index form equations, i.e., determining all power integral bases in number fields of degree $\leq 6$. In Section \ref{sec:6} 
we give applications to canonical number systems in number fields and orders,
and in Section \ref{sec:8} to some classical Diophantine equations. 
Section \ref{sec:9} gives a brief overview of generalizations, 
among others to the number field and $\fp$-adic case, and to results where the ground ring is of characteristic $0$ and finitely generated as a $\Zz$-algebra.
In Section \ref{sec:10} we give an overview of results concerning multiply monogenic and rationally monogenic orders, where we present uniform upper bounds for the multiplicity of (rational) monogenicity of orders, depending only on the degree of the underlying number field. In the Appendix we briefly discuss related topics not strictly belonging to reduction theory of integral polynomials, in particular statistical results for monogenic and rationally monogenic number fields, and Hasse's problem to give an arithmetic characterization of the monogenic number fields. 
\\[0.2cm]
\textbf{Remark.} Since the 1970's, the reduction theory of integral polynomials  of given discriminant has been constantly developing, with a growing number of results and applications. Except for Section \ref{sec:1}, the other sections contain results from this period. We propose some problems, whose solutions would yield considerable progress in the reduction theory. 
\\[0.2cm]
\textbf{Acknowledgments.} We are very much indebted to Professor Rafael von K\"{a}nel for his important contributions to Subsection \ref{sec:4.4new}, and Professor Attila Peth\H o for his useful comments on Section \ref{sec:6}. We are very grateful to Dr. Csan\'ad Bert\'ok for typing a substantial part of the manuscript.
The second named author was supported in part from the
Austrian-Hungarian joint project ANN130909 (FWF-NKFIH) and from NKFIH 150284.

\section{Reduction theory of integral 
quadratic and cubic polynomials of given non-zero discriminant}\label{sec:1}

As we mentioned, Lagrange (1773) was the first to develop a reduction theory for binary quadratic forms with integral coefficients. His theory was made more precise by Gauss (1801). For integral polynomials, their theories imply the following.
Recall that the \emph{height} $H(g)$ of a polynomial with integral coefficients
is the maximum of the absolute values of its coefficients.

\begin{theorem}[{Lagrange}, 1773; {{Gauss}}, 1801]\label{thm:1.1}
	For any quadratic polynomial $f\in\Zz[X]$ of discriminant $D\ne 0$, there exists $g\in\Zz[X]$, $GL_2(\Zz)$-equivalent to $f$, such that $H(g)\le c(D)$ with some effectively computable constant $c(D)$ depending only on $D$.
\end{theorem}

For monic polynomials, the following more precise variant is known.

\begin{theorem}\label{thm:1.2}
	For any monic quadratic polynomial $f\in\Zz[X]$ of discriminant $D\ne 0$, there exists $g\in\Zz[X]$, $\Zz$-equivalent to $f$, such that $H(g)\le c'(D)$ with some effectively computable constant $c'(D)$ depending only on $D$.
\end{theorem}

The above results have the following effective equivalent variants.

\begin{theorem}\label{thm:1.3}
	There are only finitely many $GL_2(\Zz)$-equivalence (resp. $\Zz$-equivalence) classes of quadratic (resp. monic quadratic) polynomials in $\Zz[X]$ of given discriminant $D\ne 0$. Further, each equivalence class has a representative of height at most $c(D)$ (resp. $c'(D)$).
\end{theorem}

Later, mostly these equivalent versions were investigated, used and generalized.

Hermite (1848, 1851) studied integral binary forms of degree larger than $2$. He developed an effective reduction theory for such forms which implies, among other things, the following:

\begin{theorem}[Hermite, 1848, 1851]\label{thm:1.4}
	There are only finitely many $GL_2(\Zz)$-equivalence classes of cubic polynomials in $\Zz[X]$ of given non-zero discriminant, and a full set of representatives of these classes can be effectively determined (in the sense that the proof provides an algorithm to determine, at least in principle, a full system of representatives).
\end{theorem}

In fact, Hermite (1848, 1851) introduced another invariant for polynomials $f\in\Zz [X]$ of arbitrary degree, which is in fact the discriminant $\Delta_f$ of a positive definite binary quadratic form $\Phi_f(X,Y)=AX^ 2+BXY+CY^2\in\Rr [X]$ associated with $f$. He called $f$ reduced if $\Phi_f$
is reduced in Gauss' sense, i.e., if $|B|\le A\le C$. He showed that $f$ is $GL_2(\Zz )$-equivalent to a reduced polynomial $g$, and that the coefficients of $g$ are bounded effectively in terms of $\Delta_f$. Hermite showed further that for cubic $f$, 
$\Delta_f =|27D(f)|^{1/4}$, implying Theorem \ref{thm:1.4}.
Hermite's theory was made more precise by {Julia} (1917).

For more details about reduction theories of integral binary forms and polynomials of low degree we refer to {{Dickson}}, Vol. 3 (1919, reprint\-ed 1971), {{Cremona}} (1999), {{Evertse}} and {{Gy\H ory}} (2017), Bhargava and Yang (2022), and for more general results and applications, also to Section \ref{sec:3} of the present paper and the references given there.

For the number of $\Zz$-equivalence classes of \textit{{cubic monic}} integral polynomials with given non-zero discriminant, no finiteness results were known before 1930. Then {{Delone}} and {{Nagell}} proved independently the following.

\begin{theorem}[{{Delone}}, 1930; {{Nagell}}, 1930]\label{thm:1.5} 
Up to $\Zz$-equivalence, there are only finitely many irreducible cubic monic polynomials in 
$\Zz[X]$ of given non-zero discriminant.
\end{theorem}

The proofs of Delone and Nagell of Theorem \ref{thm:1.5} were both \textit{{ineffective}}, 
in that they did not provide a method to determine 
the polynomials.
In fact, these proofs were based on a classical ineffective finiteness theorem of {{Thue}} (1909) on \textit{Thue equations}, i.e. on equations of the form $F(x,y)=m$, $x,y\in\mathbb{Z}$, where $F\in\mathbb{Z}[X,Y]$ is an irreducible binary form of degree $\ge 3$ and $m$ is an integer. In some concrete cases {{Delone}} and {{Faddeev}} (1940) made effective Theorem \ref{thm:1.5}, and posed the problem to make it effective for any irreducible cubic monic polynomial. An effective version of Theorem \ref{thm:1.5} follows from the famous effective result of {{Baker}} (1968b) on Thue equations.

%
%
%

\section{Hermite's attempt (1857) to extend the reduction results of polynomials of degree $\leq 3$ to polynomials of arbitrary degree}\label{sec:2}

\subsection{$\mathbf{GL_n(\Zz )}$-equivalence of decomposable forms}\label{sec:2.1}
\vphantom{a}

\noindent
Hermite tried to extend his theorem (1851) on cubic integral binary forms resp. polynomials to the case of any degree $n\ge 4$, but without success. Instead, he proved a finiteness theorem with a \textit{weaker equivalence}, see {Theorem \ref{thm:2.2}} below. Hermite's notion of equivalence (called by us `Hermite equivalence') is based on an equivalence relation for certain
\textit{{decomposable forms}}.

Consider decomposable forms of degree $n\ge 2$ in the same number $n$ of variables
$$F(\mathbf{X})=a_0\prod_{i=1}^n(\alpha_{i,1}X_1+\cdots+\alpha_{i,n}X_n)\in\Zz[X_1,\ldots,X_n],$$
where $a_0$ is a non-zero rational number and $\alpha_{i,j}$ are algebraic numbers, not all zero, for $i,j=1,\ldots,n$. The \textit{{discriminant}} of $F$ is defined as
$$D(F):=a_0^2(\det(\alpha_{i,j}))^2.$$
It is important to note that $D(F)$ is a rational integer.

Let $GL_n(\Zz )$ denote the multiplicative group of $n\times n$ integer matrices 
 of determinant $\pm 1$. Two decomposable forms  $F,G$ as above are called \emph{$GL_n(\Zz)$}-equivalent if 
$$G(\mathbf{X})=\pm F(U\mathbf{X})\ \text{ for some } U\in GL_n(\Zz),$$
where $\mathbf{X}$ denotes the column vector of variables
$(X_1,\ldots,X_n)^T$.

It is easy to see that two $GL_n(\Zz)$-equivalent decomposable forms in $n$ variables have the same discriminant.

Hermite proved the following.
\begin{theorem}[Hermite, 1851]\label{thm:2.1}
	Let $n$ and $D$ be integers with $n\ge 2$, $D\ne 0$. Then the decomposable forms in $\Zz[X_1,\ldots,X_n]$ of degree $n$ and discriminant $D$ lie in finitely many $GL_n(\Zz)$-equivalence classes.
\end{theorem}

\subsection{Hermite equivalence of polynomials and Hermite's finiteness theorem}\label{sec:2.2}
\vphantom{a}

\noindent
Let
$$f(X)=a_0(X-\alpha_1)\cdots(X-\alpha_n)\in\Zz[X]$$
be an integral polynomial with $a_0\in\Zz\setminus\{0\}$, and $\alpha_1,\ldots,\alpha_n\in\overline{\Qq}$. Then the \textit{{discriminant}} of $f$ is
$$D(f)=a_0^{2n-2}\prod_{1\le i<j\le n}(\alpha_i-\alpha_j)^2\in\Zz.$$
To $f$ we associate the \textit{decomposable form}
$$[f](\mathbf{X}):=a_0^{n-1}\prod_{i=1}^n(X_1+\alpha_iX_2+\cdots+\alpha_i^{n-1}X_n)\in\Zz[X_1,\ldots,X_n].$$
Using the properties of Vandermonde determinants, one can prove that 
\begin{equation}\label{vandermonde}
D([f])=D(f).
\end{equation}

The following equivalence relation was introduced by Hermite (1857):

\begin{itemize}
\item \textit{Two polynomials $f,g\in\Zz[X]$ of degree $n$ are said to be} Hermite equivalent \textit{if the associated decomposable forms $[f]$ and $[g]$ are $GL_n(\Zz)$-equivalent, i.e.,}
$$[g](\mathbf{X})=\pm[f](U\mathbf{X})\text{ \textit{for some} } U\in GL_n(\Zz).$$
\end{itemize}

From \eqref{vandermonde} it follows directly that Hermite equivalent polynomials in $\Zz[X]$ have the same discriminant.

Hermite's Theorem \ref{thm:2.1} on decomposable forms and identity \eqref{vandermonde} imply the following {finiteness theorem on polynomials.}

\begin{theorem}[Hermite, 1854, 1857]\label{thm:2.2}
	Let $n\ge 2$ and $D\ne 0$ be integers. Then the polynomials $f\in\Zz[X]$ of degree $n$ and of discriminant $D$ lie in finitely many Hermite equivalence classes.
\end{theorem}

Hermite's proof is \textit{ineffective}.

\subsection{Comparison between Hermite equivalence and $GL_2(\Zz)$-equival\-ence and $\Zz$-equivalence}\label{sec:2.3}
\vphantom{a}

\noindent
In our five authors paper with {{Bhargava}}, {{Remete}} and {{Swaminathan}} (BEGyRS, 2023) we have integrated Hermite's long-forgotten notion of equivalence and his finiteness theorem in the reduction theory, have corrected a faulty reference to Hermite's result in Narkiewicz' excellent book (2018) and compared Hermite's theorem with the most significant results of this area; see the next Section \ref{sec:3}.

In BEGyRS (2023) we proved that $GL_2(\Zz)$-equivalence and, in the monic case, $\Zz$-equivalence imply Hermite equivalence.

\begin{theorem}[BEGyRS, 2023]\label{prop:2.3}
	Let $f,g\in\Zz[X]$ be two $\Zz$-equivalent, resp. $GL_2(\Zz)$-equivalent integral polynomials. Then they are Hermite equivalent.
\end{theorem}

Since $\Zz$-equivalence implies $GL_2(\Zz)$-equivalence, it suffices to prove Theorem \ref{prop:2.3} for $GL_2(\Zz)$-equivalence. We recall the proof from BEGyRS (2023). 

\begin{proof}
	Let $f,g$ in $\Zz[X]$ be any two $GL_2(\Zz)$-equivalent polynomials. Then they can be written in the form
	$f(X)=a_0\prod_{i=1}^n(X-\alpha_i)$ and $g(X)=\pm(cX+d)^nf\left(\frac{aX+b}{cX+d}\right)$, where $A:=\left(\begin{smallmatrix}
		a&b\\c&d
	\end{smallmatrix}\right)\in GL_2(\Zz)$. Thus, we have
	$$g (X)=\pm a_0\prod_{i=1}^n(\beta_iX-\gamma_i)\text{, where $\beta_i=c-a\alpha_i,\, \gamma_i=-d+b\alpha_i$}$$
	for $i=1,\ldots,n$. Define the inner product of two column vectors
	$$\mathbf{x}=(x_1,\ldots,x_n)^T,\,\,\mathbf{y}=(y_1,\ldots,y_n)^T\text{ by } \langle \mathbf{x},\mathbf{y}\rangle:=x_1y_1+\cdots +x_ny_n.$$
	Let as before $\mathbf{X}=(X_1,\ldots,X_n)^T$. Thus,
	\begin{align*}
	&[f](\mathbf{X})=a_0^{n-1}\prod_{i=1}^n\langle\mathbf{a}_i,\mathbf{X}\rangle\text{, where } 
	\mathbf{a}_i=(1,\alpha_i,\ldots,\alpha_i^{n-1})^T,
 \\
	[g](\mathbf{X})&=\pm a_0^{n-1}\prod_{i=1}^n\langle\mathbf{b}_i,\mathbf{X}\rangle\text{, where } 
	\mathbf{b}_i=(\beta_i^{n-1},\beta_i^{n-2}\gamma_i,\ldots,\gamma_i^{n-1})^T.
	\end{align*}
	Then $\mathbf{b}_i=t(A)\mathbf{a}_i$ with some $t(A)\in GL_n(\Zz)$ for $i=1,\ldots,n$. So
	\begin{align*}
		[g](\mathbf{X})=&\pm c^{n-1}\prod_{i=1}^n\langle t(A)\mathbf{a}_i,\mathbf{X}\rangle=\\
		&=\pm c^{n-1}\prod_{i=1}^n\langle \mathbf{a}_i,t(A)^T\mathbf{X}\rangle =\pm [f](t(A)^T\mathbf{X}),
	\end{align*}
	i.e. $f$ and $g$ are indeed Hermite equivalent.
\end{proof}

For integral polynomials of degree $2$ and $3$, Hermite equivalence and $GL_2(\Zz)$-equivalence coincide. For quadratic polynomials this is trivial, while for cubic polynomials this follows from a result of {{Delone}} and {{Faddeev}} (1940).

In BEGyRS (2023) we gave, for every $n\ge 4$ and both for the non-monic and for the monic case, infinite collections of polynomials in $\Zz[X]$ with degree $n$ that are Hermite equivalent but not $GL_2(\Zz)$-equivalent. More precisely we proved the following.


\begin{theorem}[BEGyRS, 2023]\label{thm:2.4}
Let $n$ be an integer $\ge 4$. 
\begin{enumerate}
	\item[(i)] There exist infinitely many Hermite equivalence classes of properly non-monic\footnote{That is, not $GL_2(\Zz)$-equivalent to any monic polynomial} primitive
\footnote{An integral polynomial is called \emph{primitive} if its coefficients have greatest common divisor $1$},		
irreducible polynomials of degree $n$ that split into more than one $GL_2(\Zz)$-equivalence class.
	\item[(ii)] There exist infinitely many Hermite equivalence classes of monic irreducible polynomials of degree $n$ that split into more than one $GL_2(\Zz)$-equivalence class.
\end{enumerate}
\end{theorem}

In the monic case every $GL_2(\Zz)$-class contains a $\Zz$-equivalence class, hence in (ii) $GL_2(\Zz)$-equivalence can be replaced by $\Zz$-equivalence.\\

We proved Theorem \ref{thm:2.4} simultaneously for the cases (i) and (ii). We constructed, for every integer $n\ge 4$, an infinite parametric family of pairs $(f_{t,c}^{(n)}, g_{t,c}^{(n)})$ of primitive,
irreducible polynomials $f_{t,c}^{(n)}$, $g_{t,c}^{(n)}$ of degree $n$, where $c$ runs through $1$ and an infinite set of primes, and $t$ runs through an infinite set of primes with $t\ne c$ with the following properties:
\begin{align}
	&\text{for each $n$, $f_{t,c}^{(n)},\, g_{t,c}^{(n)}$ have leading coefficient $c$ and are }\label{eq:2.1}
 \\[-0.1cm]
  &\text{properly non-monic if $c>1$;}\nonumber
 \\
	&\text{for each $n$, $f_{t,c}^{(n)},\, g_{t,c}^{(n)}$ are Hermite equivalent;}\label{eq:2.2}
 \\
 	&\text{for each $n$, $f_{t,c}^{(n)},\, g_{t,c}^{(n)}$ are not $GL_2(\Zz)$-equivalent;}\label{eq:2.3}
 \\
  &\text{the pairs $(f_{t,c}^{(n)}, g_{t,c}^{(n)})$ ($n=1,2,\ldots$) lie in different Hermite }
\label{eq:2.4}
\\[-0.1cm]
&\text{equivalence classes.}\nonumber
\end{align}
The main steps of the proof are as follows. From the construction of $f_{t,c}^{(n)}$ and  $g_{t,c}^{(n)}$ it is easy to show that \eqref{eq:2.1} and \eqref{eq:2.2} hold. The proof of \eqref{eq:2.3} is more complicated. It requires the use of an irreducibility theorem of Dumas (1906), Chebotarev's density theorem, and Dirichlet's theorem on primes in arithmetic progressions. Finally, $f_{t,c}^{(n)}$ is so chosen that if we fix $n,c$ and let $t\to\infty$ then the absolute value of the discriminant of $f_{t,c}^{(n)}$  tends to $\infty$. By making a selection, we may assume that the discriminants of the polynomials $f_{t,c}^{(n)}$ are pairwise different. Since Hermite equivalent polynomials have the same discriminant, we obtain \eqref{eq:2.4}.

\begin{rem}
	We note that in our paper BEGyRS (2023) it turned out that the Hermite equivalence class of a polynomial has a very natural interpretation in terms of the so-called invariant order and invariant ideal associated with the polynomial, see Theorem \ref{thm:B5} in Subsection \ref{sec:5.5} for more details. This fact turned out to be important in the above proofs.
\end{rem}

{Theorem} \ref{prop:2.3} and {Theorem} \ref{thm:2.4} imply that $GL_2(\Zz)$-equivalence, resp. $\Zz$-equivalence are \textit{stronger} than Hermite equivalence, and hence that Hermite's {Theorem} \ref{thm:2.2} is weaker than the most significant results of this area presented in Section \ref{sec:3} below.

\section{Reduction theory of integral polynomials of given non-zero discriminant and of arbitrary degree}\label{sec:3}

As was mentioned in the Introduction, the breakthroughs in the reduction theory due to {Birch} and {Merriman} (1972), {Gy\H ory} (1973), and  {Evertse} and {Gy\H ory} (1991a)
settled the old problem of Hermite (1857),
to prove that for every given $n\ge 2$ and $D\not= 0$ there are up to $GL_2(\Zz )$-equivalence only finitely many polynomials
$f\in \Zz [X]$ of degree $n$ and discriminant $D$, and to determine these effectively. 
We state the results in more detail.

\subsection{The theorems of Birch and Merriman (1972), Gy\H ory (1973) and Evertse and Gy\H ory (1991a)}\label{sec:3.1}
\vphantom{a}

\noindent


\begin{theorem}[Birch and Merriman, 1972]\label{thm:3.1}
	Let $n\ge 2$ and $D\ne 0$. There are only finitely many $GL_2(\Zz)$-equivalence classes of polynomials in $\Zz[X]$ of degree $n$ and discriminant $D$.
\end{theorem}

Birch and Merriman established this theorem in an equivalent form, in terms of integral binary forms. Their proof uses the finiteness of the number of solutions of  \textit{unit equations} $ax+by=1$ in units $x,y$ of the ring of integers of a number field,
for which at the time effective proofs were available, but it combines this with some ineffective arguments. Consequently, Birch's and Merriman's proof of Theorem \ref{thm:3.1} is ineffective.

For monic polynomials, the corresponding result with $\Zz$-equivalence was proved \textit{independently} by {Gy\H ory} (1973) but in an \textit{effective} form. This turned out to be of crucial importance in many applications; see e.g. Sections \ref{sec:4} to \ref{sec:9} below and {Evertse} and {Gy\H ory} (2017).

\begin{theorem}[Gy\H ory, 1973]\label{thm:3.2}
	Let $f\in\Zz[X]$ be a monic polynomial of degree $n\ge 2$ with discriminant $D\ne 0$. Then 
 \begin{enumerate}
	    \item[(i)] $n\le c_1(|D|)$, and
            \item[(ii)] there is a monic $g\in\Zz[X]$, $\Zz$-equivalent to $f$, such that
		$$H(g)\le c_2(n,|D|),$$
 \end{enumerate}
    	where $c_1$ and $c_2$ are effectively computable positive numbers depending on $|D|$, resp. on $n$ and $|D|$.
\end{theorem}

This theorem was first proved and published in Gy\H ory's PhD dissertation Gy\H ory (1972a) and was utilized in Gy\H ory (1972b) as well.

\begin{corollary}[Gy\H ory, 1973]\label{cor:3.3}
	There are only finitely many $\Zz$-equivalence classes of monic polynomials in $\Zz[X]$ of given non-zero discriminant, and a full set of representatives of these classes can be at least in principle determined.
\end{corollary}

In {Gy\H ory} (1974), an explicit version was given; see below.

In his proof of Theorem \ref{thm:3.2}, {Gy\H ory} combined his own effective result on unit equations obtained by {Baker's method}, with his so-called `graph method'. We sketch below the proof of Theorem \ref{thm:3.2}.

Theorem \ref{thm:3.1}, resp. Theorem \ref{thm:3.2} and its Corollary \ref{cor:3.3} are generalizations of the corresponding results presented in Section \ref{sec:1} for polynomials of degree $n\le 3$; Theorem \ref{thm:3.1} gives an ineffective generalization of Theorem \ref{thm:1.4} for degree $n\ge 4$ and Theorem \ref{thm:3.2} is an effective generalization of Theorem \ref{thm:1.3} in the monic case for degree $n\ge 3$, and of Theorem \ref{thm:1.5} for any monic polynomial of degree $n\ge 3$.

In 1991, Evertse and Gy\H ory gave a new, effective proof for Birch's and Merriman's theorem, proving the following.

\begin{theorem}[Evertse and Gy\H ory, 1991a]\label{thm:3.4}
	Let $f\in\Zz[X]$ be a polynomial of degree $n\ge 2$ and discriminant $D\ne 0$. There is $g\in\Zz[X]$, $GL_2(\Zz)$-equivalent to $f$, such that
	$$H(g)\le c_3(n,|D|),$$
	where $c_3(n,|D|)$ is an effectively computable number, given explicitly in terms of $n$ and $|D|$.
\end{theorem}

This theorem was stated and proved in Evertse and Gy\H ory (1991a) in an equivalent form, in terms of integral binary forms.

As was mentioned above, Theorems \ref{thm:3.2} and \ref{thm:3.4} led to a general effective reduction theory of integral polynomials of given non-zero discriminant.

The main tool in our proof of Theorem \ref{thm:3.4} is an effective result of {Gy\H ory} (1974) \textit{on homogeneous unit equations} in \textit{three unknowns}, whose proof is also based on {Baker's theory of logarithmic forms}.

 We note that Theorems \ref{thm:3.1} and \ref{thm:3.4} were established directly  in a more general form, in the number field and $\mathfrak{p}$-adic case. For such  and other generalizations of Theorem \ref{thm:3.2}, (ii), see Gy\H ory (1978b, 1984) and Section \ref{sec:9} below.

Theorems \ref{thm:3.2} and \ref{thm:3.4}, their \textit{explicit} versions below and their various generalizations have a great number of consequences and applications; see our book {Evertse} and {Gy\H ory} (2017) and Sections \ref{sec:4} to \ref{sec:9} below.

\subsection{Explicit versions of theorems of Gy\H ory (1973) and Evertse and Gy\H ory (1991a)}
\vphantom{a}

\noindent
First we present explicit versions of Theorem \ref{thm:1.1}, Theorem \ref{thm:1.2} and Theorem \ref{thm:1.4} in the quadratic and cubic cases. An explicit version of Theorem \ref{thm:1.1} is the following.

\begin{customthm}{2.1*}\label{thm:3.5}
	Let $f\in\Zz[X]$ be a quadratic polynomial of discriminant $D\ne 0$. Then $f$ is $GL_2(\Zz)$-equivalent to a quadratic polynomial $g\in\Zz[X]$ such that
	\begin{enumerate}
		\item[(i)] $H(g)\le |D|/3$ if $D<0$;
		\item[(ii)] $H(g)\le |D|/4$ if $D>0$ and $f$ is irreducible;
		\item[(iii)] $H(g)\le D^{1/2}$ if $D>0$ and $f$ is reducible.
	\end{enumerate}
\end{customthm}

In the cubic case, we have the following.

\begin{customthm}{2.4*}\label{thm:3.6}
	Let $f\in\Zz[X]$ be a cubic polynomial of discriminant $D\ne 0$. Then $f$ is $GL_2(\Zz)$-equivalent to a cubic polynomial $g\in\Zz[X]$ such that
	\begin{enumerate}
		\item[(i)] $H(g)\le\frac{64}{27}|D|^{1/2}$ if $f$ is irreducible;
		\item[(ii)] $H(g)\le\frac{64}{3\sqrt{3}}|D|$ if $f$ is reducible.
	\end{enumerate}
\end{customthm}

We note that the arguments in the proofs of Theorems \ref{thm:3.5} and \ref{thm:3.6} are a variation on the arguments in {Julia} (1917). For the details we refer to Subsection 13.1 of the book of {Evertse} and {Gy\H ory} (2017).

In the monic case, it is relatively simple to prove the following explicit version of Theorem \ref{thm:1.2}.

\begin{customthm}{2.2*}\label{thm:3.7}
	For any monic quadratic polynomial $f\in\Zz[X]$ with discriminant $D\ne 0$, there exist $g\in\Zz[X]$, $\Zz$-equivalent to $f$, such that
	$$H(g)\le |D|/4 + 1.$$
\end{customthm}

As was mentioned above, the first explicit version of Theorem \ref{thm:3.2} was given in {Gy\H ory} (1974). The height estimate was improved in 2017 by the authors.

We use the notation $\log^\ast x:=\max(1,\log x)$ for $x>0$.

\begin{customthm}{4.2*}[Evertse and Gy\H ory, 2017]\label{thm:3.8}
	Let $f\in\Zz[X]$ be a monic polynomial of degree $n\ge 2$ and discriminant $D\ne 0$. Then $f$ is $\Zz$-equivalent to a polynomial $g\in\Zz[X]$ for which
	\begin{align}\label{eq:3.1}
		H(g)\le\exp\{n^{20}8^{n^2+19}(|D|(\log^\ast |D|)^n)^{n-1}\}.
	\end{align}
\end{customthm}

This is in fact Theorem 6.6.2 from Evertse and Gy\H ory (2017) with a slightly larger, simplified constant in terms of $n$.

A completely explicit, improved version of Theorem \ref{thm:3.4} was also established by the authors.

\begin{customthm}{4.4*}[Evertse and Gy\H ory, 2017, Theorem 14.1.1]\label{thm:3.9}
	Let $f\in\Zz[X]$ be a polynomial of degree $n\ge 2$ and discriminant $D\ne 0$. Then $f$ is $GL_2(\Zz)$-equivalent to a polynomial $g\in\Zz[X]$ for which
	\begin{align}\label{eq:3.2}
	H(g)\le\exp\{(4^2n^3)^{25n^2}\cdot|D|^{5n-3}\}.
	\end{align}	
\end{customthm}

In both Theorems \ref{thm:3.8} and \ref{thm:3.9}, the degree $n$ of $f$ can also be explicitly estimated from above in terms of $|D|$.
\begin{theorem}[Gy\H ory, 1974]\label{thm:3.10}
	Every polynomial $f\in\Zz[X]$ with discriminant $D\ne 0$ has degree at most
	$$3+2\log|D|/\log 3.$$
\end{theorem}

For monic polynomials $f\in\Zz[X]$, the upper bound can be improved slightly to $2+2\log|D|/\log 3$.

Theorem \ref{thm:3.4} together with Theorem \ref{thm:3.10} implies the following analogue of Corollary \ref{cor:3.3}.
\begin{corollary}[Evertse and Gy\H ory, 1991a]\label{cor:3.10}
	There are only finitely many $GL_2(\Zz)$-equivalence classes of polynomials in $\Zz[X]$ of given non-zero discriminant, and a full set of representatives of these classes can be at least in principle effectively determined.
\end{corollary}

\subsection{Consequences of Theorems \ref{thm:3.9}, \ref{thm:3.8} and Theorem \ref{prop:2.3} for Hermite equivalence classes}\label{sec:4.3new}
\vphantom{a}

\noindent
As was pointed out in BEGyRS (2023),an important consequence of the above Theorem \ref{prop:2.3} is that the effective finiteness theorems \ref{thm:3.4}, \ref{thm:3.9} and \ref{thm:3.2}, \ref{thm:3.8} for $GL_2(\mathbb{Z})$-equivalence classes resp. $\mathbb{Z}$-equivalence classes apply  just as well to Hermite equivalence classes.

   We present here the following,  more precise, explicit variant of Hermite's result in Theorem \ref{thm:2.2}.

\begin{corollary}[of Theorems \ref{thm:3.9} and \ref{thm:3.8}; cf. BEGyRS, 2023]\label{cor:4.7}
\vphantom{}
\begin{enumerate}
    \item[(i)]  Every Hermite equivalence class of polynomials in $\mathbb{Z}[X]$ of degree $n\ge 2$ and of discriminant $D\ne 0$ has a representative with coefficients not exceeding
    $$\exp\{(4^2n^3)^{25n^2}|D|^{5n-3}\}$$
    in absolute value.
    \item[(ii)] Every Hermite equivalence class of monic polynomials in $\mathbb{Z}[X]$ with degree $n\ge 2$ and discriminant $D\ne 0$ has a representative with coefficients not exceeding
    $$\exp\{n^{20}8^{n^2+19}(|D|(\log^\ast |D|)^n)^{n-1}\}$$
    in absolute value.
\end{enumerate}
\end{corollary}

 It is an immediate consequence of Theorem \ref{thm:3.10} that in (i) above  $n \le 3 + 2\log |D|/\log 3$. Further, in (ii), the slightly better inequality $n\le 2 + 2\log |D|/\log 3$  holds.

 The above result implies an effective version of Theorem \ref{thm:2.2}, i.e., for given $n$ and a non-zero integer $D$, one can effectively determine a full system  of representatives for the Hermite equivalence classes of polynomials $f\in\Zz [X]$ of degree $n$ and discriminant $D$. Indeed, one can make a finite list of all polynomials $f\in\Zz [X]$ of height below one of the bounds in Corollary \ref{cor:4.7}. For each polynomial in the list one can check whether it has discriminant $D$. Further,
for each pair of polynomials in the list one can check whether they are Hermite equivalent, by computing the corresponding decomposable forms $[f],[g]$,
and checking whether they are $GL_n(\Zz )$-equivalent, using, e.g., Lemma 18 of Evertse and Gy\H{o}ry (1992a).

   The similarity of Theorems \ref{thm:3.9} , \ref{thm:3.8}  and Corollary \ref{cor:4.7} is only apparent. As was seen in Section \ref{sec:2}, the $GL_2(\mathbb{Z})$-equivalence  and $\mathbb{Z}$-equivalence are in fact much stronger than the Hermite equivalence.

   \begin{rem}
    Every improvement of the bounds in \eqref{eq:3.1} or \eqref{eq:3.2} would yield the same improvement in the bounds of Corollary \ref{cor:4.7}.
\end{rem}

\subsection{Unit equations and $S$-unit equations}
\label{sec:3.3}
\vphantom{a}

\noindent
The unit equations and more general $S$-unit equations play a fundamental role in Diophantine number theory, and in particular in the effective reduction theory of integral polynomials of given discriminant.

First we recall unit equations and $S$-unit equations, and then briefly outline how to apply Baker's theory of logarithmic forms to obtain effective bonds for the solutions of these equations. Then we recall the best known height bounds for the solutions of unit equations and $S$-unit equations over number fields.

For a detailed treatment of unit equations, $S$-unit equations and their further generalizations and applications we refer to our books Evertse and Gy\H ory (2015, 2017, 2022).

Let $K$ be an algebraic number field, $\OO_K$ its ring of integers, $\OO_K^\ast$ the unit group of $\OO_K$, and $M_K$ its set of places, consisting of the finite set of infinite places $S_{\infty}$ of $K$ (corresponding to the real embeddings  and the pairs of conjugate complex embeddings of $K$ in $\Cc$) and the finite places, which we may identify with the prime ideals of $\OO_K$.
To the places in $M_K$ we can associate a set of absolute values $\{ |\cdot |_v:\, v\in M_K\}$, normalized such that if $v$ lies above the place $p\in M_{\Qq}:=\{\infty\}\cup\{\text{primes}\}$, then for $a\in\Qq$
one has $|a|_v=|a|_p^{[K_v:\Qq_p]}$. These absolute values satisfy the product formula $\prod_{v\in M_K} |\alpha |_v=1$ for $\alpha\in K^*$.

Let $a,b$ be given non-zero elements of $K$. Equations of the form
\begin{align}\label{eq:4.3K}
ax+by=1\text{ in unknowns } x,y\in\OO_K^\ast
\end{align}
are called  \textit{unit equations} (\textit{in two unknowns}).
More generally, 
let $S$ be a finite subset of $M_K$ with $S\supseteq S_{\infty}$.
Denote by $\OO_S$ the ring of $S$-integers, i.e., $\{ x\in K:\, |x|_v\leq 1\ \text{for } v\not\in S\}$ and by $\mathcal{O}_S^\ast$ denote the unit group of $\OO_S$, i.e., group of $S$-units. Thus, $\OO_S^\ast=\{ x\in K:\, |x|_v=1\ \text{for } v\not\in S\}$.
For $S=S_{\infty}$ we have $\OO_S^\ast =\OO_K^\ast$. 
Equations of the form
\begin{align}\label{eq:4.3final}
	a x + b y =1\text{ in unknowns } x,y\in\mathcal{O}_S^\ast
\end{align}
are called \textit{$S$-unit equations} (\textit{in two unknowns}). In many cases it is more convenient to consider the unit equations and $S$-unit equations in homogeneous form
\begin{align}\label{eq:4.3a}
	a x+b y+c z=0\text{ in unknowns } x,y,z\in \OO_K^\ast ,\ \text{resp. } \OO_S^\ast,\tag{4.4a}
\end{align}
where $a,b,c$ denote fixed elements of $K\setminus\{0\}$.

For a long time these equations were utilized merely in special cases and in an implicit way. It was implicitly proved by Siegel (1921) for $S=S_\infty$ and by Parry (1950) for any $S$ that equation \eqref{eq:4.3final} has only finitely many solutions. This implies the finiteness of the number of solutions of equation \eqref{eq:4.3a} up to a common proportional factor. Lang (1960) gave a direct proof for a more general version of these finiteness theorems. Their proofs were ineffective.

Generalizing Gelfond's (1935) famous result obtained in the case $m=2$, in the 1960's Baker made a major breakthrough in number theory by giving non-trivial explicit lower bounds for the absolute value of linear forms in logarithms of the form
$$b_1\log\alpha_1+\cdots+b_m\log\alpha_m\ne 0,\,\,m\ge 2$$
where $b_1,\ldots,b_m$ are rational integers, resp. algebraic numbers, $\alpha_1,\ldots,\alpha_m$ are algebraic numbers different from $0$ and $1$, and $\log\alpha_1,\ldots,\log\alpha_m$ denote fixed determination of the logarithms. In case of rational integers $b_1,\ldots,b_m$, this is equivalent to bounding $|\prod\alpha_i^{b_i}-1|$ non-trivially from below. Baker's general effective estimates led to significant applications, and opened a new effective epoch in the theory of Diophantine equations. Baker's quantitative results were later improved, generalized, extended to the $\mathfrak{p}$-adic case and so on by himself and many other authors; for comprehensive overviews we refer to Baker (1990), W\"ustholz, ed. (2002), and Baker and W\"ustholz (2007), and for a shorter overview see Evertse and Gy\H ory (2015), Section 3.2. The last five decades saw the development of an \textit{effective theory} of Diophantine equations.

General effective upper bounds for the solutions of 
\eqref{eq:4.3K} and \eqref{eq:4.3a} in the case $S=S_{\infty}$ were deduced by Gy\H ory (1972a,b, 1973) using an effective result of Baker and Coates (1970), p. 601, on relative Thue equations over number fields. The first \textit{explicit} upper bounds for the solutions of \eqref{eq:4.3K} and \eqref{eq:4.3a} in case $S=S_{\infty}$ were deduced by Gy\H ory (1974) from an explicit inequality of Baker (1968a, Part IV) for linear forms in logarithms of algebraic numbers. For general $S$, Gy\H ory (1979) derived the first explicit bound for the solutions of \eqref{eq:4.3final}, using also the $\mathfrak{p}$-adic version of Baker's theory. Independently, a slightly weaker effective bound was given by Kotov and Trelina (1979).

Let $K$ be an algebraic number field. Given $\alpha_1\kdots \alpha_n\in K$, not all $0$, we define the height of $(\alpha_1\kdots \alpha_n )$ relative to $K$ by
\[
H_K(\alpha_1\kdots \alpha_n):=\prod_{v\in K}\max (|\alpha_1|_v\kdots |\alpha_n|_v).
\]
Recall that the naive height
$H(\alpha )$ of an algebraic number $\alpha$ is given by the maximum of the absolute values of its minimal polynomial, with coefficients having gcd $1$. Then we have
\[
H(\alpha )\leq 2^{\deg\alpha}H_{\Qq (\alpha )}(1,\alpha ). 
\]

Following Section 1.3 from the paper "Solving Diophantine equations by Baker's theory" by Gy\H ory (2002), we briefly \textit{sketch} a proof of the following theorem, by means of Baker's theory. 

\begin{theorem}\label{thm:4.xx}
Let $K$ be a number field, $S$ a finite set of places of $K$ containing $S_{\infty}$, and $a,b $ non-zero elements of $K$.
Let $x,y\in\OO_S^\ast$ 
satisfy \eqref{eq:4.3final}. Then 
\[
\max (H(x) ,H(y))\leq c_4(K,S,a,b),
\]
where $c_4$ is an effectively computable number, depending only on $K,S,a,b$.
\end{theorem}

\begin{proof}[Sketch]
Let $s$ denote the cardinality of $S$.
There is a system of fundamental $S$-units $\{\varrho_1,\ldots,\varrho_{s-1}\}$ in $\OO_S^*$ with heights bounded in terms of $K$ and $S$. Let $x,y$ be a solution of \eqref{eq:4.3final} in $S$-units. Then one can write
$$x=\xi_1\varrho_1^{a_{11}}\cdots \varrho_{s-1}^{a_{1,s-1}},\,\,y=\xi_2\varrho_1^{a_{21}}\cdots\varrho_{s-1}^{a_{2,s-1}},$$
where $\xi_1,\xi_2$ are roots of unity in $K$ and $a_{ij}$ are unknown rational integer exponents. Assume without loss of generality that $A:=\max_j |a_{1j}|\geq \max_j|a_{2j}|$. By elementary means one can show that
\[
A\leq c_5\log \max_{v\in S} |x|_v,
\]
and combining this with $\prod_{v\in S} |x|_v=1$,
one concludes that there is a $v\in S$ such that
$$|x|_v\le c_6\exp\{-c_7A\},$$
where $c_5,c_6,c_7$ can be given explicitly and depend only on $K$ and $S$.  This implies
\begin{align}\label{eq:4.4final}
	0<|\varrho_1^{a_{21}}\cdots\varrho_{s-1}^{a_{2,s-1}}-\alpha|_v\le c_8\exp\{-c_9A\}
\end{align}
with an appropriate $\alpha\in K$ of bounded height. The constants $c_8, c_9$ and $c_{10}$ below depend at most on $K,S$ and $a,b$ and can be given explicitly.

One can now apply the complex or $\mathfrak{p}$-adic version of Baker's theory according as $v\in S_\infty$ or $v\in S\setminus S_\infty$ and this yields
$$\exp\{-c_{10}\log A\}\le |\varrho_1^{a_{21}}\cdots\varrho_{s-1}^{a_{2,s-1}}-\alpha|_v.$$
Comparing this with \eqref{eq:4.4final} we get
\begin{align}\label{eq:4.5final}
	A\le A_0
\end{align}
where $A_0$ can be given explicitly. Finally, we obtain an upper bound for $H(x)$ and $H(y)$ which can also be given explicitly.
\end{proof}

Later, several improvements, effective generalizations, applications and algorithmic results have been obtained for unit and $S$-unit equations by means of Baker's theory; see among others Gy\H ory (1980b, 2002, 2019, 2022), Shorey and Tijdeman (1986), Sprind\v{z}uk (1993), Bugeaud and Gy\H ory (1996), Smart (1998), Ga\'al and Gy\H ory (1999), Hindry and Silverman (2000), W\"ustholz, ed. (2002), Bilu (2002), Bilu, Ga\'al and Gy\H ory (2004), Gy\H ory and Yu (2006), Baker and W\"ustholz (2007), Zannier (2009), Hajdu (2009), B\'erczes, Evertse and Gy\H ory (2009), Evertse and Gy\H ory (2013, 2015 2017, 2022), 
B\'erczes (2015a, 2015b), Bert\'ok and Hajdu (2015, 2018), Bugeaud (2018), Ga\'al (2019), Le Fourn (2020), Alvarado et al. (2021), Gy\H ory and Le Fourn (2024), and the references given there. 

The best known height bound for the solutions of \eqref{eq:4.3K} is due to Gy\H{o}ry and Yu (2006). We formulate it in simplified form. 
As above, let $K$ be a number field of degree $d$ and $r$ the rank of $\OO_K^*$. 
Denote by $h_K,\, R_K$ the class number and regulator of $K$, respectively,
and write again $\log^\ast x:=\max (1,\log x)$.

 Then Gy\H{o}ry and Yu (2006) proved the following.

\begin{theorem}\label{thm:4.wwK}
Let $a,b$ be non-zero elements of $K$. Then for all $x,y\in\OO_K^*$ satisfying \eqref{eq:4.3K} we have
\[
H_K(1,x,y)\leq (3H_K(1,a,b))^A,
\]
where
\[
A=d^5(2r+2)^{4r+40}R_K\log^\ast R_K.
\]
\end{theorem}

\noindent
\textbf{Remark.} The following inequality implies that $A$ can be bounded abouve in terms of $d$ and $D_K$ only:
\[
h_KR_K\leq |D_K|^{1/2}(\log^* |D_K|)^{d-1}.
\]
The first inequality of this type was proved by Landau (1918). For the above version,
see, e.g., Evertse and Gy\H{o}ry (2015, formula (1.5.2)).
\\[0.15cm]

In this section we shall use only Theorem \ref{thm:4.wwK}. Below we formulate a generalization to $S$-unit equations, Theorem \ref{thm:4.ww}, which is not used in this section, but will be needed in Section \ref{sec:9}.  

Consider the general case where $S$ is an arbitrary finite set of places containing $S_{\infty}$.
In terms of $S$, the best known bounds can be found in Gy\H ory (2019), Le Fourn (2020) and Gy\H ory and Le Fourn (2024). We mention here the bound from Gy\H{o}ry (2019) in simplified form. We introduce the necessary notation.
Let as above $K$ be a number field of degree $d$ and $r$ the rank of $\OO_K^*$. Denote by $h_K$, $R_K$ the class number and regulator of $K$, respectively. Further, let $S=S_{\infty}\cup\{ \fp_1\kdots \fp_t\}$, where $\fp_1\kdots \fp_t$ with $t\geq 0$ are the prime ideals in $S$.
Let $s:=\# S$, and denote by $R_S$ the $S$-regulator. It is known that
\[
R_S =R_K\ \text{if } t=0,\ \ R_S=i_SR_K\prod_{i=1}^t\log N_K\fp_i\ \text{otherwise,}
\]
where $i_S$ is a divisor of $h_K$
and $N_K\fa$ denotes the norm of a non-zero ideal $\fa$ of $\OO_K$, i.e., $\#\OO_K/\fa$.
Let $P_S:=1$ if $t=0$ and 
$P_S:=\max_{1\leq i\leq t} N_K\fp_i$ if $t\geq 1$.
Further, put $P_S':=1$ if $t\leq 2$ and $P_S'$ the third largest among the quantities $N_K\fp_i$, $i=1\kdots t$ if $t\geq 3$. Finally, put $\TT_K:=\max (h_K, 160r!\cdot (r+1)^2R_K)$.

\begin{theorem}\label{thm:4.ww} Let $a,b$ be non-zero elements of $K$.
Then for all $x,y\in\OO_S^*$ with \eqref{eq:4.3final} we have 
\[
H_K(1,x,y)\leq (3H_K(1,a,b))^{A_S},
\]
where
\[
A_S:=
 2s^5(16ed)^{4s+3}\TT_K^{t+4}\cdot\frac{P_S'}{\log^* P_S'}\Big(1+\frac{\log^*\log P_S}{\log^*P_S'}\Big)R_S.
\]
\end{theorem}

Observe that for $S=S_{\infty}$, $A$ is much smaller than $A_S$.
Further, $A_S$ can be bounded above in terms of $d$, $|D_K|$, $t$, and $P_S$.

We compare Theorem \ref{thm:4.ww} with the abc-conjecture over number fields.
We first recall the abc-conjecture over $\Qq$, as proposed by Masser (1985), refining an earlier conjecture of Oesterl\'e.
Define the \emph{radical} of a non-zero integer $a$ by $\RR (a):=\prod_{p|a} p$.

\begin{conjecture}[Masser-Oesterl\'e abc-conjecture, 1985]\label{abc-over-Q}
There is a constant $C(\epsilon )>0$ depending on $\epsilon$ such that for all $\epsilon >0$ and all non-zero integers $a,b,c$ with $a+b=c$ and $\gcd (a,b,c)=1$ we have $\max (|a|,|b|,|c|)\leq C(\epsilon )\RR (abc)^{1+\epsilon}$.
\end{conjecture}

There are various proposals to extend this to number fields.
We recall a version of Masser (2002).
Let $K$ be a number field and $D_K$ its discriminant.
Take a non-zero ideal $\fa$ of $\OO_K$. Masser defined the modified
radical of $\fa$ by $\RR_K(\fa ):=\prod_{\fp |\fa} N_K\fp^{e_{\fp}}$, where the product is taken over all prime ideals dividing $\fa$ and  $e_{\fp}$ is the ramification index of $\fp$. Masser considered this modified radical since it has a good behaviour under field extensions, e.g., if $L$ is an extension of $K$ of degree $m$, then $\RR_L(\fa\OO_L)=\RR_K(\fa )^m$.

Recall that the different of $K$ can be expressed as $\fD_K =\prod_{\fp} \fp^{w_\fp}$, where the product is taken over all prime ideals $\fp$ of $\OO_K$ with $e_{\fp}>1$, and where $w_{\fp}\geq e_{\fp}-1$. Further, $|D_K|=N_K\fD_K$. This implies that for any ideal $\fa$ of $\OO_K$, 
\begin{equation}\label{4.yy}
\RR_K'(\fa )\, |\, \RR_K(\fa )\, | \, D_K\cdot\RR_K'(\fa),\ \ \text{where } \RR_K'(\fa ):=\prod_{\fp|\fa } N_K\fp .
\end{equation}

\begin{conjecture}[Masser's uniform abc conjecture over number fields, 2002]\label{conjecture-abc}
There is a constant $C(\epsilon )>0$ depending on $\epsilon$, such that for every  $\epsilon >0$ the following holds. For every number field $K$ of discriminant $D_K$ and every non-zero $\alpha ,\beta ,\gamma\in K$ with $\alpha +\beta =\gamma$, we have
\[
H_K(\alpha ,\beta ,\gamma )\leq C(\epsilon )^{[K:\Qq ]}\big (|D_K|\cdot \RR_K(\fa^{-3}\alpha\beta\gamma)\big)^{1+\epsilon},
\]
where $\fa$ is the fractional ideal generated by $\alpha ,\beta ,\gamma$.
\end{conjecture} 

This implies the following bound for the solutions of the $S$-unit equation \eqref{eq:4.3final} $a x+b y=1$  in $x,y\in\OO_S^\ast$, where again $S$ is a finite set of places of $K$, containing the infinite places and $a,b\in K^*$: let $\RR_S :=1$ if $S=S_{\infty}$ and $\RR_S:=\prod_{i=1}^t N_K\fp_i^{e_{\fp_i}}$,
and put $\RR_K (a,b):=\prod_{\fp} N\fp^{e_{\fp}}$, where the product is taken over all 
$\fp\in M_K\setminus S$ such that $|a|_{\fp}$ and $|b|_{\fp}$ are not both equal to $1$. Then for every solution $x,y\in\OO_S^\ast$ of $ax+by=1$ we have
\[
H_K(1,x,y)\leq C(\epsilon )^d\big( |D_K|\cdot \RR_S\cdot \RR_K (a,b)\big)^{1+\epsilon}H_K(1,a,b)^2.
\] 
See also Gy\H ory (2022), Theorem 3.

Some alternative effective methods were also developed to obtain effective bounds for the solutions of $S$-unit equations. Bombieri (1993, 2002) and Bombieri and Cohen (1997, 2003) worked out such an effective method in Diophantine approximation, based on an extended version of the Thue--Siegel principle, the Dyson Lemma and some geometry of numbers. Bugeaud (1998), following their approach and combining it with estimates for linear forms in logarithms, proved results which are in certain parameters sharper than those of Bombieri and Cohen.

During 1983--95 Frey initiated and developed in several papers the modular degree approach for $S$-unit equations over $\mathbb{Q}$; see e.g. Frey (1997) where he gives height bounds which became unconditional around 2000 when modularity was proved. As is surveyed by von K\"anel (2024), effective bounds over $\mathbb{Q}$ were proved in 2011 independently and simultaneously by von K\"anel (2013, 2014b) and by Murty and Pasten (2013), Pasten (2014).

However, it should be remarked that for most applications of $S$-unit equations, including the reduction theory of integral polynomials treated in our paper, more general results concerning $S$-unit equations of the form \eqref{eq:4.3final} over arbitrary number fields are needed.

\subsection{A brief sketch of the proof of a less precise version of Theorem \ref{thm:3.2}}
\label{sec:3.3a}
\vphantom{a}

\noindent
Consider a monic polynomial $f\in\mathbb{Z}[X]$ of degree $n$ and discriminant $D\ne 0$. In view of Theorem \ref{thm:1.3} we may assume that $n\ge 3$.

First we sketch the proof of assertion (i). Assume that $f$ is irreducible over $\mathbb{Q}$. Let $K=\mathbb{Q}(\alpha)$ for a zero $\alpha$ of $f$, and denote by $D_K$ the discriminant of $K$. Then combining the Minkowski inequality with the fact that $D_K$ divides $D(f)$, i.e. $D$, (i) follows with an appropriate $c_1$. If now $f$ is reducible and $f=f_1\cdots f_t$ with monic irreducible $f_1,\ldots,f_t$, then using $D(f_j)\mid D(f)$ in $\mathbb{Z}$ and applying the just proved (i) for $j=1,\ldots,t$, we obtain (i) in the general case as well.


We now sketch the proof of (ii) in Theorem \ref{thm:3.2}. Its main steps are as follows.
\\[0.15cm]
\textbf{1.}
	Denote by $\alpha_1,\ldots,\alpha_n$ the zeros of $f$, and by $G$ the splitting field of $f$ over $\Qq$. Then $[G:\Qq]\le n!$ and the absolute value $|D_G|$ of the discriminant of $G$ can be estimated from above by a constant $c_{11}(n,|D|)$. Here and below 
 $c_{11},\ldots $ are effectively computable numbers depending only on $n$ and $|D|$. 
\\[0.15cm] 
\textbf{2.} Putting $\Delta_{ij}:=\alpha_i-\alpha_j$ we have
	$$\prod_{1\le i< j\le n}\Delta_{ij}^2=D,$$
	which implies $|N_{G/\Qq}\Delta_{ij}|\le c_{12}(n,|D|)$. It follows that 
	\begin{align}\label{eq:3.3}
		\Delta_{ij} = \delta_{ij}\varepsilon_{ij}\text{, where } H(\delta_{ij})\le c_{13}(n,|D|)
	\end{align}
	and $\varepsilon_{ij} $ is a unit in the ring of integers of $G$.
\\[0.15cm]
\textbf{3.}	
The following identity plays a basic role in the proof:
	\begin{align}\label{eq:3.4}
		\Delta_{ij}+\Delta_{jk}=\Delta_{ik}\text{ for every } i,j,k.
	\end{align}
	Consider the \textit{graph}, whose vertices are $\Delta_{ij}$ $(1\le i\not= j\le n)$ and whose edges are $[\Delta_{ij},\Delta_{ik}]$, $[\Delta_{ij},\Delta_{jk}]$ ($1\leq i\not= j\le n$, 
 $k\not= i,j$). This graph is obviously connected.
\\[0.15cm] 
\textbf{4.} Equations \eqref{eq:3.3} and \eqref{eq:3.4} give rise to a `connected' system of \textit{unit equations}
	\begin{align}\label{eq:3.5}
		\delta_{ijk}\varepsilon_{ijk}+\tau_{ijk}\nu_{ijk}=1,
	\end{align}
	where $\delta_{ijk}:=\delta_{ij}/\delta_{ik}$,
 $\tau_{ijk}:=\delta_{jk}/\delta_{ik}$ are non-zero elements of $G$ with heights effectively bounded above in terms of $n$ and $|D|$ only, and $\varepsilon_{ijk}:=\varepsilon_{ij}/\varepsilon_{ik}$, $\nu_{ijk}:=\varepsilon_{jk}/\varepsilon_{ik}$ are \textit{unknown} units in the ring of integers of $G$.
\\[0.15cm]	
\textbf{5.} 
Applying Theorem \ref{thm:4.wwK}, together with the Remark following it,
we get upper bounds for the heights of the quotients $\Delta_{ij}/\Delta_{ik}=\delta_{ijk}\varepsilon_{ijk}$
  for each triple $\{ i,j,k\}\subset\{ 1\kdots n\}$, depending on $G,n$ and $|D|$, and so eventually only on $n$ and $|D|$, and likewise for $\Delta_{jk}/\Delta_{ik}$.
\\[0.15cm]
\textbf{6.}
  Using the connectedness of the unit equations involved, this yields effective upper bounds for the height of $\Delta_{ij}$ for every $i,j$, depending only on $n$ and $|D|$.
  Indeed, one first obtains an upper bound for the height of any quotient $\Delta_{ij}/\Delta_{kl}$ via 
  $$
  \frac{\Delta_{ij}}{\Delta_{kl}}= \frac{\Delta_{ij}}{\Delta_{ik}}\cdot
   \frac{\Delta_{ik}}{\Delta_{kl}}
   $$
 (using the path $\Delta_{ij}\to\Delta_{ik}\to\Delta_{kl}$ in the graph) and subsequently for the height of each $\Delta_{ij}$ via
$$
\Delta_{ij}^{n(n-1)}=\pm D\cdot\prod_{1\leq k\not= l\leq n}
 \frac{\Delta_{ij}}{\Delta_{kl}}.
$$
\textbf{7.}
Adding the differences $\Delta_{ij}=\alpha_i-\alpha_j$ for fixed $i$ and for $j=1,\ldots,n$, using the fact that $\alpha_1+\cdots +\alpha_n\in\Zz$, putting $\alpha_1+\cdots+\alpha_n=na+a'$ with $a,a'\in\Zz$, $0\le a'<n$, and writing
 \begin{align*}
 	&\beta_i :=\alpha_i-a\text{ for } i=1,\ldots,n,
 \\
 &g(X)=\prod_{i=1}^n(X-\beta_i),
 \end{align*}
 we have that $g(X)=f(X+a)\in\Zz[X]$ and that the height of $g$ has an effective upper bound depending only on $n$ and $|D|$. \qed

\begin{rem}
    We note that for cubic and quartic monic polynomials $f\in\mathbb{Z}[X]$ of given non-zero discriminant, Klaska (2021, 2022) devised another approach for proving Corollary \ref{cor:3.3} via the theory of integral points on elliptic curves.
\end{rem}

\subsection{A brief sketch of the proof of a less precise version of Theorem \ref{thm:3.4}}\label{sec:3.4}
\vphantom{a}

\noindent
Take an integral polynomial $f\in\Zz[X]$ of degree $n$ and discriminant $D\ne 0$. In view of Theorems \ref{thm:1.1}, \ref{thm:1.2} and \ref{thm:1.4} we may assume that $n\ge 4$. The  absolute value of the discriminant of the splitting field of $f$ can be estimated from above in terms of $|D|$, and by the Hermite-Minkowski Theorem, this leaves only a finite, effectively determinable collection of possible splitting fields for $f$. 
So we may restrict ourselves to polynomials $f$ with given splitting field $G$ and ring of integers $\mathcal{O}_G$. 

Take such $f$ and pick a factorization of $f$,
\begin{align}\label{eq:3.6}
	f=\prod_{i=1}^n (\alpha_i X-\beta_i)\,\,\text{over }\overline{\Qq},
\end{align}
such that the number of linear factors with real coefficients is maximal, and the factors with complex coefficients fall apart into complex conjugate pairs.
After multiplying $f$ by a small positive rational integer, which can be effectively bounded in terms of $G$, hence in terms of $n$ and $|D|$ and
which is negligible compared with the other estimates arising from the application of Baker's method, we may assume that $f$ has such a factorization with $\alpha_i,\beta_i\in\mathcal{O}_G$ for $i=1,\ldots,n$. Put 
$$
\Delta_{ij}:=\alpha_i\beta_j-\alpha_j\beta_i\ \ 
\text{for $1\le i,j\le n$.}
$$

We now follow the approach of Evertse and Gy\H{o}ry (2017),
chapters 13 and 14. We outline the main steps of the proof.
\\[0.15cm]
\textbf{1.}
We start with a small variation on the reduction theory of 
Hermite (1848, 1851) and Julia (1917). Let $\tv =(t_1\kdots t_n)$ be a tuple of positive reals such that $t_i=t_j$ for each pair $(i,j)$ such that $\alpha_i,\beta_i$ are the complex conjugates of 
$\alpha_j,\beta_j$. Consider the positive definite quadratic
form
$$
\Phi_{f,\tv}(X,Y):=\sum_{i=1}^n t_i^{-2}(\alpha_iX-\beta_iY)(\overline{\alpha_i}X-\overline{\beta_i}Y).
$$
By Gauss' reduction theory for positive definite binary quadratic forms,
there is $\big(\begin{smallmatrix}a&b\\c&d\end{smallmatrix}\big)\in GL_2(\Zz )$ such that
$\Phi_{f,\tv}(aX+bY,cX+dY)$ is reduced, i.e.,
equal to $AX^2+BXY+CY^2$ with $|B|\leq A\leq C$.
Define the polynomial 
$$
g(X)=(cX+d)^nf\left(\frac{aX+b}{cX+d}\right),
$$
which is $GL_2(\Zz )$-equivalent to $f$.
We denote by $H(g)$ the height of $g$. We recall Theorem 13.1.3 of Evertse and Gy\H{o}ry (2017), and refer for the elementary proof to section 13.1 of that book.

\begin{proposition}\label{prop:5.9}
Let
\begin{align*}
M:=t_1\cdots t_n,\ \ 
R:=\left( \sum_{1\le i<j\le n}\frac{|\Delta_{ij}|^2}{t_i^2t_j^2}\right)^2.
\end{align*}
Then
\begin{align*}
H(g)\le\left(\frac{4}{n\sqrt{3}}\right)^nM^2R^n
\end{align*}
if $f$ has no root in $\Qq$, and
\begin{align*}
H(g)\le\left(\frac{2}{\sqrt{n}}\right)^n\cdot
\left(\frac{2}{\sqrt{3(n-1)}}\right)^{n(n-1)/(n-2)}(M^2R^n)^{(n-1)/(n-2)}
\end{align*}
if $f$ does have a root in $\Qq$.
\end{proposition}
\vskip0.15cm\noindent
\textbf{2.}
For any quadruple $i,j,k,l$ of distinct indices we have the identity
\begin{align}\label{eq:3.7}
	\Delta_{ij}\Delta_{kl}+\Delta_{jk}\Delta_{il}=\Delta_{ik}\Delta_{jl}.
\end{align}
Notice that all terms $\Delta_{ij}$ are in $\mathcal{O}_G$ and divide $D$. Hence $|N_{G/\Qq}(\Delta_{ij})|\le |D|^{[G:\Qq]}$ for all $i,j$ where $[G:\Qq]\le n!$. As above in Section \ref{sec:3.3}, we can express each term $\Delta_{ij}$ as a product of an element of height effectively bounded in terms of $n,D$ and a unit from $\mathcal{O}_G$. By substituting this into the identities \eqref{eq:3.7} we obtain homogeneous unit equations in three terms. Dividing \eqref{eq:3.7} by $\Delta_{ik}\Delta_{jl}$ we get unit equations like in \eqref{eq:3.5} above, and using Theorem \ref{thm:4.wwK}
we obtain effective upper bounds for the heights of the quotients $\Delta_{ij}\Delta_{kl}/\Delta_{ik}\Delta_{jl}$. 
\\[0.15cm]
\textbf{3.}
To obtain an effective upper bound for the height of $g$ in terms of $n$ and $|D|$,
it suffices to effectively estimate the quantities $M$ and $R$ from Proposition \ref{prop:5.9} from above in terms of $n$ and $|D|$, for a suitable choice 
of the $t_i$. For the $t_i$ we choose
$$
t_i:=\left( \prod_{k=1,k\not= i}^n|\Delta_{ik}|\right)^{1/(n-2)}\ \ \text{for } i=1\kdots n.
$$
With this choice,
$$
M=|D|^{1/(n-2)}
$$
and
$$
\frac{|\Delta_{ij}|}{t_it_j}=\left(|D|^{-1}\cdot
\prod_{k,l}\left|\frac{\Delta_{ij}\Delta_{kl}}{\Delta_{ik}\Delta_{jl}}\right|\right)^{1/(n-1)(n-2)},
$$
where the product is taken over all pairs of indices $k,l$
such that $1\le k,l\le n$, $k\not= i,j$, $l \not= i,j$ and $k\not= l$. By inserting the upper bounds for the heights 
of the quantities $\Delta_{ij}\Delta_{kl}/\Delta_{ik}\Delta_{jl}$ obtained in the previous step,
we can estimate from above $M$ and $R$, and subsequently $H(g)$, effectively in terms of $n$ and $|D|$ only.
\qed

\subsection{Conjectural improvements (partly joint work with von K\"anel)}\label{sec:4.4new}
\vphantom{a}

\noindent
This subsection contains important contributions by Rafael von K\"anel.

As was mentioned above, for $n\ge 4$ resp. $n\ge 3$ the proofs of Theorems \ref{thm:3.2}, \ref{thm:3.4}, \ref{thm:3.8} and \ref{thm:3.9} are based on effective results of Gy\H ory on unit equations
 whose proofs depend on Baker's theory of logarithmic forms. The exponential feature of the bounds in \eqref{eq:3.1} and \eqref{eq:3.2} is a consequence of the use of Baker's method. It is likely that the bounds in \eqref{eq:3.1} and \eqref{eq:3.2} can be replaced by some bounds polynomial in terms of $|D|$. This can be achieved if we restrict ourselves to polynomials $f\in\Zz[X]$ having a \textit{fixed splitting field} $G$ over $\Qq$. In this case the bounds in \eqref{eq:3.1} and \eqref{eq:3.2} can be replaced by bounds of the form 
$$c_{14}(n,G)|D|^{c_{15}(n,G)},$$ 
where $c_{14}(n,G), c_{15}(n,G)$ are effectively computable numbers which depend only on $n$ and the discriminant of $G$; see {Gy\H ory} (1984, 1998) resp. {Evertse} and {Gy\H ory} (1991a). The following conjectures seem plausible.

\begin{conjecture}\label{conjecture-y}
	Let $f\in\Zz[X]$ be a monic polynomial of degree $n\ge 3$ and discriminant $D\ne 0$. Then $f$ is 
	$\Zz$-equivalent to a monic polynomial $g$ in $\mathbb{Z}[X]$ such that
	$$H(g)\le c_{16}(n)|D|^{c_{17}(n)}$$
	where $c_{16}(n), c_{17}(n)$ depend only on $n$.
\end{conjecture}

\begin{conjecture}\label{conjecture-x}
	Let $f\in\Zz[X]$ be a polynomial of degree $n\ge 4$ and of discriminant $D\ne 0$. Then $f$ is $GL_2(\Zz)$-equivalent to a polynomial 
	 $g$ in $\mathbb{Z}[X]$ such that
	$$H(g)\le c_{18}(n)|D|^{c_{19}(n)}$$
	where $c_{18}(n), c_{19}(n)$ depend only on $n$.
\end{conjecture}

Conjecture \ref{conjecture-x} has been formulated in Chapter 15 of {Evertse} and {Gy\H ory} (2017).
In fact, Conjecture \ref{conjecture-x} implies Conjecture \ref{conjecture-y}.

\begin{proof}[Conjecture \ref{conjecture-x} $\Longrightarrow$ Conjecture \ref{conjecture-y}]
Let $f\in\Zz [X]$ be a monic polynomial of degree $n\ge 3$ and discriminant $D\not= 0$. 
Consider the polynomial $g(X):= (2X+1)^{n+1}f(\medfrac{X}{2X+1})$. Using that $f$ is monic, one shows by means of a straightforward computation that $g$ has degree $n+1$ and $D(g)=D$.
By Conjecture \ref{conjecture-x} there is 
$\big(\begin{smallmatrix}a&b\\ c&d\end{smallmatrix}\big)\in GL_2(\Zz )$ such that $g^*(X):=(cX+d)^{n+1}g(\medfrac{aX+b}{cX+d})$ has height at most $c_{18}(n+1)|D|^{c_{19}(n+1)}$. A straightforward computation shows that
\begin{align*}
&g^*(X)=(c'X+d')f^*(X),
\\
&\text{with } c'=2a+c,\, d'=2b+c,\ \ f^*(X)=(c'X+d')^nf(\medfrac{aX+b}{c'X+d'}).
\end{align*}
Note that $|c'|, |d'|, H(f^*)\leq c_{20}(n)H(g^*)$. 
 Let $r$ be an integer such that $a':=a+rc'$ satisfies $|a'|\leq\half |c'|$. Then from $ad'-bc'=\pm 1$ it follows that $b':=b'+rd'$ satisfies $|b'|\leq \half |d'|+1$.
 Now define $f^{**}(X):=(-c'X+a')^nf^*(\medfrac{d'X-b'}{-c'X+a'})$. One verifies that
 $f^{**}(X)=f(\pm X\pm r)$ and $H(f^{**})\leq c_{16}(n)|D|^{c_{17}(n)}$.
 \end{proof}

We give some evidence for the conjectures mentioned above.
Evertse proved the following what one may call semi-effective result.

\begin{theorem}[Evertse, 1993]\label{thm:3.8new}
	Let $f\in\Zz[X]$ be a polynomial of degree $n\ge 4$ and of discriminant $D\ne 0$, having splitting field $G$ over $\Qq$. Then $f$ is $GL_2(\Zz)$-equivalent to a polynomial $g$ of height
	$$H(g)\le c_{21}(n,G)|D|^{21/n}.$$
 \end{theorem}
 \noindent
	Here $c_{21}(n,G)$ is a number depending only on $n$ and $G$, which is not effectively computable by the method of proof. For a proof, see also Evertse and Gy\H{o}ry (2017, Chap. 15).
	
The main tool in Evertse's proof is the following theorem.
The constant in this theorem is ineffective. Let $K$ be a number field. Given $\alpha ,\beta ,\gamma\in K$,
we define the height $H_K(\alpha ,\beta ,\gamma ):=\prod_{v\in M_K}\max (|\alpha |_v,|\beta |_v ,|\gamma |_v)$.

\begin{theorem}\label{thm:4.7y1}
Let $\alpha ,\beta ,\gamma$ be non-zero elements of $\OO_K$ with $\alpha +\beta =\gamma$. Then for all $\epsilon >0$ we have
\[
H_K(\alpha ,\beta ,\gamma )\leq c_{22}(K,\epsilon )|N_{K/\Qq}(\alpha\beta\gamma)|^{1+\epsilon},
\]
where $c_{22}(K,\epsilon )$ depends only on $K$ 
and $\epsilon$.
\end{theorem} 

In fact, this is a special case of a general multivariable result of Evertse (1984b, Theorem 1), see also Evertse and Gy\H{o}ry (2015, Theorem 6.1.1).  The proof of this general result is based on Schmidt's Subspace Theorem over number fields. For Theorem \ref{thm:4.7y1} one needs the 
two-dimensional case, which is Roth's Theorem over number fields. Theorem \ref{thm:3.8new} was deduced from Theorem \ref{thm:4.7y1} essentially by following the arguments in Subsection \ref{sec:3.4}, but with various refinements to obtain a bound with an exponent $O(1/n)$ on $|D|$.

In order to prove Conjecture \ref{conjecture-x}, the following variation on Theorem \ref{thm:4.7y1} would suffice:

\begin{conjecture}\label{conjecture-z}
For all number fields $K$ of degree $d\ge 2$ and discriminant $D_K$ and all non-zero
$\alpha ,\beta ,\gamma\in \OO_K$ with  $\alpha +\beta =\gamma$ we have 
\[
H_K(\alpha ,\beta ,\gamma )\leq c_{23}(d) |D_K\cdot 
N_{K/\Qq}(\alpha\beta\gamma)|^{c_{24}(d)},
\]
where $c_{23}(d)$, $c_{24}(d)$ depend
 only on $d$.
\end{conjecture}

This obviously follows from  Masser's uniform abc-conjecture over number fields,
i.e. Conjecture \ref{conjecture-abc}, but is of course much weaker.

\begin{proof}[Conjecture \ref{conjecture-z}
$\Longrightarrow$ Conjecture \ref{conjecture-x} (sketch)]
We follow the argument in Subsection \ref{sec:3.4}, and use the same notation. Let $f\in \Zz [X]$ be a polynomial of degree $n\geq 4$ and discriminant $D\not= 0$. Denote by $G$ the splitting field of $f$. By e.g., Evertse and Gy\H{o}ry (2017, Corollary 13.3.4), there is $a\in\Qq$ with $1\leq |a|\leq c_{25}(n)|D_G|^{c_{26}(n)}$ such that $f_1:=af=\prod_{i=1}^n(\alpha_iX-\beta_i)$ with $\alpha_i,\beta_i\in\OO_G$ for $i=1\kdots n$, and such that the non-real factors among the 
$\alpha_iX-\beta_iY$ can be divided into complex conjugate pairs. Let $D_1:=D(f_1)$. Now define $\Delta_{ij}:= \alpha_i\beta_j-\alpha_j\beta_i$ $(1\leq i<j\leq n)$ and apply Conjecture \ref{conjecture-z} to the identities
\[
\Delta_{ij}\Delta_{kl}+\Delta_{jk}\Delta_{il}=\Delta_{ik}\Delta_{jl}.
\]
Noting that $|N_{G/\Qq}(\Delta_{ij})|\leq |D_1|^{n!}$, it follows that for all quadruples $i,j,k,l$,
\[
H_G(\Delta_{ij}\Delta_{kl},\Delta_{jk}\Delta_{il},\Delta_{ik}\Delta_{jl})\leq c_{27}(n)|D_G \cdot D_1|^{c_{28}(n)}.
\]
This leads to upper bounds for the quantities $|\Delta_{ij}\Delta_{kl}/\Delta_{ik}\Delta_{jl}|$. Following the arguments in part 3 of Subsection \ref{sec:3.4}, applying Proposition \ref{prop:5.9}, one obtains that $f_1$ is $GL_2(\Zz )$-equivalent to a polynomial $g_1$ with 
\[
H(g_1)\leq c_{29}(n) |D_G\cdot D_1|^{c_{30}(n)}.
\]
One can show that $D_G$ divides $D_1^{c_{31}(n)}$. Taking $g:=a^{-1}g_1$ one obtains that $g$ is $GL_2(\Zz )$-equivalent to $f$ and that $H(g)\leq c_{18}(n)|D|^{c_{19}(n)}$.
\end{proof}

We are interested in upper bounds for $H(g)$ that depend as much as possible on $D_G$ and the radical of $D=D(f)$, and as little as possible on $D$ itself. Under assumption of 
 Conjecture \ref{conjecture-abc} (Masser's version of the abc-conjecture over number fields), we deduce the following result for monic polynomials. In fact, it is a modification of some ideas of Rafael von K\"anel, which he kindly shared with us.
Recall that the radical of a non-zero rational integer $a$ is defined by $\RR (a):=\prod_{p|a} p$.

\begin{theorem}\label{thm:4.7y.1}
Under assumption of Conjecture \ref{conjecture-abc}, the following holds.
	Let $f\in\Zz[X]$ be a monic polynomial of degree $n\ge 3$ and of discriminant $D\ne 0$. 
	Let $G$ be the splitting field of $f$ and $D_G$ its discriminant. Then $f$ is 
	$\Zz$-equivalent to a monic polynomial $g\in \Zz [X]$ such that
	\[
	H(g)\leq c_{32}(n)\big( |D_G\cdot\RR (D)|\big)^{c_{33}(n)}\cdot |D|^{1/(n-1)},
	\]
	where $c_{32}(n)$, $c_{33}(n)$ depend only on $n$.
\end{theorem}
\noindent
\textbf{Remark.} With a more elaborate computation, $c_{33}(n)$ can be computed explicitly.

\begin{proof} We use the following notation: we write $A\ll^* B$ if there are positive numbers $c'(n)$, $c''(n)$, depending only on $n$, such that\\ $A\leq c'(n)|D_G\cdot \RR (D)|^{c''(n)}B$.
At each occurrence of $\ll^*$, the constants $c'(n),c''(n)$ may be different.

Write $f(X)=(X-\alpha_1)\cdots (X-\alpha_n)$. Choose a rational integer $a$ such that $|a-(\alpha_1+\cdots +\alpha_n)/n|\leq\half$, and take $g(X):=f(X+a)$. 
This $g$ is clearly $\Zz$-equivalent to $a$. Then
\[
H(g)\leq 2^n\prod_{i=1}^n\max (1,|\alpha_i-a|)\leq 2^n\prod_{i=1}^n\max(1,\half +|\alpha_i-(\alpha_1+\cdots +\alpha_n)/n|),
\]
hence
\begin{equation}\label{eq:4.7a}
H(g)\leq 2^n\prod_{i=1}^n\big( 1+n^{-1}\sum_{j=1}^n |\alpha_i-\alpha_j|\big).
\end{equation}
We prove Theorem \ref{thm:4.7y1} by estimating
the right-hand side from above, and to this end we apply Conjecture \ref{conjecture-abc} to the identities 
\[
(\alpha_i-\alpha_j)+(\alpha_j-\alpha_k)=(\alpha_i-\alpha_k)\ \ \ (i,j,k\in\{ 1\kdots n\}\ \text{pairwise distinct}).
\]
Note that all terms in this sum are algebraic integers in $G$, composed of prime ideals in $\OO_G$ dividing $D$. So by Conjecture \ref{conjecture-abc},
\[
H_G\Big(1,\frac{\alpha_i -\alpha_j}{\alpha_i-\alpha_k}\Big)\leq
H_G(\alpha_i-\alpha_j,\alpha_j-\alpha_k,\alpha_i-\alpha_k)\ll^* 1.
\]
This implies
\[
\Big{|}\frac{\alpha_i-\alpha_j}{\alpha_i-\alpha_k}\Big{|} \ll^* 1\ \ \text{for all pairwise distinct $i,j,k$}
\]
and subsequently, using 
$
\medfrac{\alpha_i-\alpha_j}{\alpha_k-\alpha_l}=-\medfrac{\alpha_i-\alpha_j}{\alpha_i-\alpha_k}
\cdot\medfrac{\alpha_k-\alpha_i}{\alpha_k-\alpha_l}
$,
\[
\Big{|}\frac{\alpha_i-\alpha_j}{\alpha_k-\alpha_l}\Big{|} \ll^* 1\ \ \text{for all pairwise distinct $i,j,k,l$.}
\]
This leads us to
\[
|\alpha_i-\alpha_j|\ll^*\Big(\prod_{1\leq k\not=l\leq n} |\alpha_k-\alpha_l|\Big)^{1/(n(n-1)}=|D|^{1/n(n-1)}\ \ \text{for all } i\not= j.
\]
By inserting this into \eqref{eq:4.7a}, we arrive at $H(g)\ll^* |D|^{1/(n-1)}$.
This completes our proof.	
\end{proof}

Rafael von K\"anel kindly communicated to us a conjecture on monic cubic polynomials of given discriminant that is equivalent to the Masser-Oesterl\'e abc-conjecture over $\Qq$, i.e., Conjecture \ref{abc-over-Q}. To formulate von K\"anel's conjecture, we introduce the weighted height
of $f=X^3+a_1X^2+a_2X+a_3\in\Zz [X]$ by
\[
\Ht (f):=\max (|a_1|,|a_2|^{1/2},|a_3|^{1/3}).
\]
Further, we introduce the quantity
\begin{align}\label{eq:smalldelta}
\delta_f :=&\max \{d\in\Zz :\, d^2|P\ \text{and } d^3|U\},
\\
\notag
&\quad\text{where } P:=a_1^2-3a_2,\ U:=2a_1^3+27a_3-9a_1a_2.
\end{align}
Here $P$ and $U$ are the usual two seminvariants of $f$, which satisfy $4P^3-U^2=27D$,
where $D=D(f)$. Note that $\delta_f^6$ divides $27D$.

\begin{conjecture}[von K\"anel]\label{conjecture-w}
There is a constant $c_{34}(\epsilon )>0$ depending on $\epsilon$,  
such that for every real $\epsilon >0$ the following holds:
\\
For every monic cubic polynomial $f\in\Zz [X]$ of discriminant $D\not= 0$, there is a polynomial $g\in \Zz [X]$ that is $\Zz$-equivalent to $f$ and for which
\[
\Ht (g) \leq c_{34}(\epsilon )\cdot \delta_f\cdot\RR (27D/\delta_f^6)^{1+\epsilon}.
\]
\end{conjecture}

\begin{theorem}[von K\"anel]\label{thm:4.7y.2}
The Masser-Oesterl\'e abc-conjecture over $\Qq$  is equivalent to Conjecture \ref{conjecture-w}. 
\end{theorem}

\noindent
\textbf{Remark.}
It might be possible to extend the proof to prove a version for any number field $K$ without introducing substantial new ideas. However it will be clear that the proof does not work for polynomials of degree $\geq 4$. 

Observing that $\delta_f\cdot \RR (27D/\delta_f^6)$ divides $27D$, this implies at once the following:

\begin{corollary}\label{cor:4.7y.3}
Assume the Masser-Oesterl\'e abc-conjecture over $\Qq$ holds. Then there is a constant $c_{35}(\epsilon )>0$ depending on $\epsilon$, such that for every real$\epsilon >0$ 
the following holds:
\\
For every monic cubic polynomial $f\in\Zz [X]$ of discriminant $D\not= 0$, there is a polynomial $g\in \Zz [X]$ that is $\Zz$-equivalent to $f$ and for which
\[
\Ht (g) \leq c_{35}(\epsilon )\cdot |D|^{1+\epsilon}.
\]
\end{corollary} 
Noting that $H(g)\leq \Ht (g)^3$, Corollary \ref{cor:4.7y.3} immediately implies a version of Conjecture \ref{conjecture-y}.

\begin{proof}[Proof of Theorem \ref{thm:4.7y.2}]
We follow von K\"{a}nel's argument.

It is known (see Bombieri-Gubler (2005, 12.5.12)) that the Masser-Oesterl\'e abc-conjecture over $\Qq$ is equivalent to the following
\begin{conjecture}\label{conjecture-v}
For every real $\epsilon>0$ there is a constant $c_{36}(\epsilon )$ such that all $u,v\in\Zz$ with $w:=u^3-v^2\not= 0$ and
$\gcd(u^3,v^2)$ sixth power-free satisfy 
\[
|u|\leq c_{36}(\epsilon )\cdot \RR(w)^{2+\epsilon}, \quad |v| \leq c_{36}(\epsilon )\cdot \RR(w)^{3+\epsilon}.
\]
\end{conjecture}

Therefore it suffices to show that Conjecture \ref{conjecture-w} is equivalent to 
Conjecture \ref{conjecture-v}. This equivalence is a consequence of Lemmas~\ref{lem1} and \ref{lem2} that we shall prove below.
\end{proof}

In what follows we write $A\ll_{\epsilon} B$ if there is a constant $c(\epsilon )>0$ depending only on $\epsilon$ such that $A\leq c(\epsilon )B$.

\begin{lemma}\label{lem1}
Conjecture \ref{conjecture-v} implies Conjecture \ref{conjecture-w}.
\end{lemma}

\begin{proof}
We assume that Conjecture \ref{conjecture-v} holds and we let $\epsilon>0$ be a real number. 

Let $f\in\Zz [X]$ be a cubic monic polynomial of discriminant $D\neq 0$. 
Write $f=X^3+a_1X^2+a_2X+a_3$ with $a_i\in \Zz$, and let $\delta =\delta_f$, $P$, $U$ be as in \eqref{eq:smalldelta}. We compute
\begin{equation}\label{eq:comp}
f(X-\tfrac{a_1}{3})=X^3+b_2X+b_3, \quad b_2=-\tfrac{P}{3}, \quad b_3=\tfrac{U}{27}, \quad 4P^3-U^2=27D.
\end{equation}
The definition of $\delta$ assures that $P_0=P/\delta^2$ and $U_0=U/\delta^3$ lie in $\Zz$ with $\gcd(P_0^3,U_0^2)$ sixth power-free. Moreover, it follows from \eqref{eq:comp} that $P_0$ and $U_0$ satisfy
$$
(4P_0)^3-(4U_0)^2=16\cdot 27(D/\delta^6).
$$
Next we define $\rho:=\max\{d\in\Zz : \  d^2\mid 4P_0 \textnormal{ and } d^3\mid 4U_0\}$. Then we observe that $u=4P_0/\rho^2$ and $v=4U_0/\rho^3$ lie in $\Zz$ with $\gcd(u^3,v^2)$ sixth power-free, and we obtain
$$u^3-v^2=w, \quad w=\tfrac{16\cdot 27}{\rho^6}(D/\delta^6)\neq 0.$$
Here we used our assumption that $D\neq 0$. It holds that  $\RR (w)\leq 6\cdot \RR(27D/\delta^6)$ since $\rho\in\Zz$ and then an application of Conjecture \ref{conjecture-v}  with $u,v$ leads to 
\begin{equation}\label{eq:uvbound}
\max(|u|^3,|v|^2)\ll_{\epsilon}\RR (w)^{6+\epsilon}\ll_{\epsilon} \RR (27D/\delta^6)^{6+\epsilon}.
\end{equation}
As $\gcd(P_0^3,U_0^2)$ is sixth power-free, the definition of $\rho$ implies  $\rho\mid 2$.  Then, on combining \eqref{eq:uvbound} with the definitions of $b_2,b_3$ and $u,v$, we deduce
\begin{equation}\label{eq:bibound}
\max(|b_2|^{1/2},|b_3|^{1/3})\leq \delta \cdot \max(|u|^{1/2},|v|^{1/3})\ll_{\epsilon} \delta \cdot \RR (27D/\delta^6)^{1+\epsilon}.
\end{equation}
In the case when $-a_1/3\in\Zz$, we can take $g=f(X+\tau)\in\Zz[X]$ for $\tau=-a_1/3\in\Zz$. Indeed $\textnormal{Ht}(g)=\max(|b_2|^{1/2},|b_3|^{1/3})$ by \eqref{eq:comp} and thus \eqref{eq:bibound} gives $\textnormal{Ht}(g)\ll_{\epsilon}\delta \cdot \RR(27D/\delta^6)^{1+\epsilon}$.

Suppose from now on that $-a_1/3\notin\Zz$. Then we may and do choose $\sigma\in\{\tfrac{1}{3},\tfrac{2}{3}\}$ such that $\tau'=-\tfrac{a_1}{3}+\sigma\in\Zz$. Define $g=f(X+\tau')$ and write 
$g=X^3+c_1X^2+c_2X+c_3$ with $c_i\in \Zz$. On using that $g=f((X+\sigma)-\tfrac{a_1}{3})=(X+\sigma)^3+b_2(X+\sigma)+b_3$, we obtain the identities
$$
c_1=3\sigma, \quad c_2=3\sigma^2+b_2, \quad c_3=\sigma^3+b_2\sigma+b_3.
$$
The definition of $\sigma$ gives $|\sigma|\leq 2/3$, and our assumption $D\neq 0$ assures that not both $b_2,b_3$ are zero. Hence we deduce  $\textnormal{Ht}(g)\ll_{\epsilon} \max(|b_2|^{1/2},|b_3|^{1/3})$ which together with \eqref{eq:bibound} implies $\textnormal{Ht}(g)\ll_{\epsilon}\delta \cdot \RR(27D/\delta^6)^{1+\epsilon}$ as desired. This completes the proof of Lemma~\ref{lem1}.
\end{proof}

\begin{lemma}\label{lem2}
Conjecture \ref{conjecture-w} implies Conjecture \ref{conjecture-v}.
\end{lemma}

\begin{proof}
We assume that Conjecture \ref{conjecture-w} holds and we let $\epsilon>0$ be a real number. 

 Let $u,v\in \Zz$ with $\gcd(u^3,v^2)$ sixth power-free and $w=u^3-v^2\neq 0$. We consider
 the monic cubic $f=X^3+a_2X+a_3$ in $\Zz[X]$ where $a_2=-3u$ and $a_3=2v$. A direct computation shows that the discriminant $D$ of $f$ and the seminvariants $P,U$ of $f$ are given by 
 $$
 D=4\cdot 27w, \quad P=9u, \quad U=2\cdot 27v.
 $$
It follows that $D\neq 0$, since $w\neq 0$ by assumption. Moreover our assumption that $\gcd(u^3,v^2)$ is sixth power-free implies that the quantity $\delta$ in \eqref{eq:smalldelta} satisfies $\delta\mid 6$. Then an application of Conjecture \ref{conjecture-w} with $f$ gives that there is $\tau\in \Zz$ such that $g=f(X+\tau)$ satisfies 
\begin{equation}\label{eq:cibound}
\textnormal{Ht}(g)=\max_i |c_i|^{1/i}\ll_{\epsilon}\RR (D)^{1+\epsilon}
\end{equation}
where $g=X^3+c_1X^2+c_2X+c_3$ and $c_i\in \Zz$.  As $f=X^3+0\cdot X^2-3uX+2v$ we obtain that $c_1=3\tau$ and thus $g(X-\tfrac{c_1}{3})=f$. This leads to the following identities
$$
-3u=a_2=-\tfrac{c_1^2}{3}+c_2, \quad 2v=a_3=\tfrac{2}{27}c_1^3-\tfrac{c_1c_2}{3}+c_3.
$$
Thus \eqref{eq:cibound} combined with $D=4\cdot 27w$ implies $|u|\ll_{\epsilon} \RR(w)^{2+\epsilon}$ and $|v|\ll_{\epsilon}\RR(w)^{3+\epsilon}$ as desired. This completes the proof of Lemma~\ref{lem2}.
\end{proof}
 
 We finish this subsection by recalling a function field analogue 
of Conjecture \ref{conjecture-x} that has been proved unconditionally.
Let $\vk$ be an algebraically closed field of characteristic $0$. Define the polynomial ring $A:=\vk [t]$ and its quotient field $L :=\vk (t)$.

Define an absolute value $|\cdot |_{\infty}$ on $L$ as follows: if $a,b\in A$ are two non-zero polynomials, then put $|a/b|_{\infty}:=\exp (\deg a-\deg b)$. Further, define the height of $f(X):=a_0X^n+a_1X^{n-1}+\cdots +a_n\in A[X]$ by $H(f):=\max (|a_0|_{\infty}\kdots |a_n|_{\infty})$. 
Call two polynomials $f,g\in A[X]$ of degree $n$ $GL_2(A)$-equivalent, if $g(X)=u(cX+d)^nf(\medfrac{aX+b}{cX+d})$ for some $u\in\vk^*$ and 
$\big(\begin{smallmatrix}a&b\\ c&d\end{smallmatrix}\big)\in GL_2(A)$.

In his PhD-thesis, Zhuang (2015, Chap. 5, Theorem 5.3.2) proved the following result
 (in fact, Zhuang formulated this in terms of binary forms $F\in A[X,Y]$;
using the correspondence $f(X)=F(X,1)$ one obtains the theorems below).

\begin{theorem}\label{thm:4.7y2} 
Let $f\in A[X]$ be a polynomial of degree $n\geq 3$ and discriminant $D\not= 0$. Assume that $f$ has splitting field $G$ over $L$, and denote by $g_G$ the genus of $G$. Then $f$ is $GL_2(A)$-equivalent to a polynomial $g$ for which
\[
H(g)\leq 
\exp \Big(n^2+6n-7+\frac{(5n-5)(2g_G-1)}{24[G:K]}\Big)\cdot |D|_{\infty}^{21/n}.
\]
\end{theorem}

By estimating $g_G$ from above in terms of $n$ and $|D(f)|_{\infty}$, Zhuang (2015, Chap. 5, Main Theorem) obtained the following, unconditional, function field analogue of Conjecture \ref{conjecture-x}:

\begin{theorem}\label{thm:4.7y4}
Let $f\in A[X]$ be a polynomial of degree $n\geq 3$ and discriminant $D\not= 0$. Then $f$ is $GL_2(A)$-equivalent to a polynomial $g$ for which
\[
H(g)\leq \exp\big( (n-1)(n+6)\big)\cdot |D|_{\infty}^{20+(1/n)}.
\]
\end{theorem} 

The proof of Theorem \ref{thm:4.7y2} is similar to that of Theorem \ref{thm:3.8new}, except that instead of Theorem \ref{thm:4.7y1} Zhuang used the Stothers-Mason abc-Theorem for function fields.

We recall this theorem.
Let $K$ be a function field of transcendence degree $1$ over an algebraically closed field $\vk$ of characteristic $0$. Let $M_K$ be the set of normalized discrete valuations on $K$, i.e., with value group $\Zz$. These valuations satisfy the sum formula $\sum_{v\in M_K} v(x)=0$ for $x\in K^*$.
Denote by $g_K$ the genus of $K$. Define the height of a tuple $(\gamma_1\kdots \gamma_n)\in K^n$ by
$h_K(\gamma_1\kdots \gamma_n):=-\sum_{v\in M_K}\min (v(\gamma_1)\kdots v(\gamma_n))$.

\begin{theorem}\label{thm:4.7y3}
Let $\alpha ,\beta ,\gamma$ be elements of $K\setminus\vk$ such that $\alpha +\beta =\gamma$. Let $s$ denote the number of valuations $v$ of $K$ such that $v(\alpha )$, $v(\beta )$, $v(\gamma )$ are not all equal. Then
\[
h_K(\alpha ,\beta ,\gamma )\leq s+2g_K-2.
\]
\end{theorem}

For a proof, see Mason (1984).

\section{Consequences in algebraic number theory, in particular
for monogenicity and rational monogenicity} 
\label{sec:4}
 
We give some consequences 
of Theorems \ref{thm:3.1}, \ref{thm:3.2} and \ref{thm:3.4} 
in algebraic number theory. 
Of particular interest are applications to monogenicity of number fields and (rational) monogenicity of orders.

Theorem \ref{thm:3.1} due to Birch and Merriman from 1972 has an important ineffective finiteness consequence for algebraic integers of given discriminant; see Theorem \ref{thm:5.1} below.

An effective version of Theorem \ref{thm:5.1} was obtained independently in Gy\H ory (1973), as a consequence of his effective Theorem \ref{thm:3.2} presented above; see Theorem \ref{thm:5.2} below.

Theorem \ref{thm:5.2} as well as its various effective consequences, applications, quantitative variants and generalizations in Gy\H ory (1973, 1974, 1976, 1978a,b, 1980a,b, 1981) led to breakthroughs in the effective theory of number fields. These furnished among others general effective finiteness results for integral elements of given discriminant resp. of given index in number fields and, more generally, in their orders; see Corollaries \ref{cor:5.3} and \ref{cor:5.4}. In particular, as an immediate consequence of his Theorem \ref{thm:5.2}, Gy\H ory provided the \textit{first general effective algorithm} for deciding the \textit{monogenicity} and for determining, at least in principle, \textit{all power integral bases} in number fields and in their orders; see Theorems \ref{cor:5.6} and \ref{thm:5.5} below.

As a consequence of Theorem \ref{thm:3.4} we present from Evertse and Gy\H ory (1991a) a general effective finiteness theorem on algebraic numbers of given discriminant; see Theorem \ref{thm:5.6}.
Finally, we introduce rationally monogenic orders, which are generalizations of monogenic orders, and give an algorithm to determine in principle whether a given order is rationally monogenic, see Theorems \ref{thm:5.7} and \ref{cor:5.8} below.

For convenience, we formulate the above-mentioned effective finiteness results in their simplest form.
For generalizations, further applications and comprehensive treatment of this extensive area, we refer to Gy\H ory (1983, 1984, 1998, 2000, 2006), Evertse and Gy\H ory (1991a, 2017, 2022), BEGyRS (2023), the references given there, and to Sections \ref{sec:5} to \ref{sec:9} of the present paper.

\subsection{Preliminaries}\label{sec:4.1}
\vphantom{a}

\noindent
Throughout this section, $K$ will denote a number field of degree $n\ge 2$ with ring of integers $\mathcal{O}_K$ and discriminant $D_K$. Recall that $K$ has precisely $n$ distinct embeddings in its normal closure over $\Qq$, which we denote by $x\mapsto x^{(i)}$ ($i=1\kdots n$). Here $x^{(1)}=x$.

Let $\MM$ be a free $\Zz$-module in $K$ of rank $n$, and pick a $\Zz$-basis $\{\omega_1\kdots \omega_n\}$ of $\MM$. Then the discriminant of $\MM$ is defined by
\[
D(\MM ):=\Big(\det \big(\omega_i^{(j)}\big)_{i,j=1}^n\Big)^2.
\]
This is a rational number, and it does not depend on the choice of the basis. 

Given two free $\Zz$-modules $\MM_1,\MM_2$ in $K$ of rank $n$ with $\MM_1\supseteq\MM_2$, denote by $[\MM_1:\MM_2]$ the \emph{index} of $\MM_2$ in $\MM_1$, i.e., the cardinality of $\MM_1/\MM_2$. Then 
\begin{equation}\label{eq:4.x1}
D(\MM_2)=[\MM_1:\MM_2]^2D(\MM_1).
\end{equation}
We recall that an \emph{order} of $K$ is a subring of $K$ which as a $\Zz$-module is free of rank $n$.
In particular, $\OO_K$ is an order of $K$, and each other order of $K$ 
is a subring of $\OO_K$. We have $D_K=D(\OO_K)$. Equation \eqref{eq:4.x1} implies that if 
$\MM$ is a submodule of $\OO_K$ of rank $n$, then $D(\MM )\in\Zz$.
 
Let $\alpha$ be a non-zero algebraic integer. Then we denote by $f_\alpha(X)$ the minimal (monic) polynomial of $\alpha$ in $\Zz[X]$. 
Thus, $f_{\alpha}(X)=\prod_{i=1}^ n(X-\alpha^{(i)})$, where $\alpha^{(1)}=\alpha\kdots \alpha^{(n)}$
are the distinct conjugates of $\alpha$ in the splitting field of $f_{\alpha}$. We now define the \emph{discriminant} of $\alpha$
by
\begin{align}\label{eq:4.1}
	D(\alpha) :=D(f_\alpha)=\prod_{1\leq i<j\leq n} (\alpha^{(i)}-\alpha^{(j)})^ 2.
\end{align}
The ring $\Zz [\alpha ]$ is clearly an order of $\Qq (\alpha )$, with $\Zz$-module basis $\{ 1,\alpha \kdots\alpha^{n-1}\}$, so
\begin{align}\label{eq:4.x2}
D(\alpha )=D(\Zz [\alpha ]).
\end{align}

Let now $\mathcal{O}$ be an order of $K$, and $D(\mathcal{O})$ its discriminant. Then $\mathcal{O}$ is a subring of $\mathcal{O}_K$. For a primitive element $\alpha$ of $K$ with $\alpha\in\mathcal{O}_K$ resp. $\alpha\in\mathcal{O}$, we define 
\begin{equation}\label{eq:4.index}
I(\alpha):=[\mathcal{O}_K:\Zz[\alpha]],\ \ I_\mathcal{O}(\alpha):=[\mathcal{O}:\Zz[\alpha]]
\end{equation}
to be the \textit{index} of $\alpha$ in $\mathcal{O}_K$ resp. in $\mathcal{O}$. Then, by \eqref{eq:4.x1}, \eqref{eq:4.x2},
\begin{equation}
\label{eq:4.2}
D(\alpha)=I(\alpha)^2D_K\text{ for }\alpha\in\mathcal{O}_K,\ \ 
D(\alpha)=I_\mathcal{O}(\alpha)^2D(\mathcal{O})\text{ for }\alpha\in\mathcal{O}.
\end{equation}

 Two algebraic integers $\alpha$, $\beta$ are called \emph{$\Zz$-equivalent}
if $\beta =\pm\alpha +a$ for some $a\in\Zz$.
If $\alpha$ and $\beta$ are $\Zz$-equivalent then so are $f_{\alpha}$ and $f_{\beta}$. Conversely,
if $f_{\alpha}$ and $f_{\beta}$ are $\Zz$-equivalent then $\alpha$ is $\Zz$-equivalent to a conjugate of $\beta$.

Clearly, $\Zz$-equivalent elements in $\mathcal{O}_K$ resp. in $\mathcal{O}$ have the same discriminant and hence the same index in $\mathcal{O}_K$ resp. in $\mathcal{O}$.

A number field $K$ is called \textit{monogenic} if $\mathcal{O}_K=\Zz[\alpha]$ for some $\alpha\in\mathcal{O}_K$. This is equivalent to the fact that $I(\alpha)=1$ and that $\{1,\alpha,\ldots,\alpha^{n-1}\}$ is a \textit{power integral basis} in $K$,
i.e., a $\Zz$-module basis of $\OO_K$. Similarly, an order $\mathcal{O}$ of $K$ is said to be \textit{monogenic} if $\mathcal{O}=\Zz[\alpha]$, i.e. if $I_\mathcal{O}(\alpha)=1$ for some $\alpha\in\mathcal{O}$. Clearly, if $\mathcal{O}=\Zz [\alpha ]$ then also $\mathcal{O}=\Zz [\beta ]$ for every $\beta$ that is $\Zz$-equivalent to $\alpha$. 

Further, $K$ resp. $\mathcal{O}$ is called $k\,(\ge 1)$ \textit{times monogenic} if $\mathcal{O}_K$ resp. $\mathcal{O}$ equals $\Zz[\alpha_1]=\cdots=\Zz[\alpha_k]$ for some pairwise $\Zz$-inequivalent $\alpha_1,\ldots,\alpha_k$ in $\mathcal{O}_K$ resp. in $\mathcal{O}$. In case that in the above definition $k$ is maximal, it is called the \textit{multiplicity} of the monogenicity of $K$, resp. of $\OO$.

\subsection{Consequences of Theorems \ref{thm:3.1} and \ref{thm:3.2} for algebraic integers of given discriminant} 
\vphantom{a}

\noindent
From their Theorem \ref{thm:3.1}, Birch and Merriman in 1972 deduced the following ineffective finiteness theorem.
\begin{theorem}[Birch and Merriman, 1972]\label{thm:5.1}
    Up to $\mathbb{Z}$-equivalence, there are only finitely many algebraic integers with given non-zero discriminant.
\end{theorem}

Independently, as a consequence of his Theorem \ref{thm:3.2}, Gy\H ory (1973) proved the following effective version of Theorem \ref{thm:5.1}.

By the \textit{height} $H(\alpha)$ of an algebraic integer $\alpha$ we mean the height $H(f_\alpha)$. 
\begin{theorem}[Gy\H ory, 1973]\label{thm:5.2}
    Let $\alpha$ be an algebraic integer of degree $n\ge 2$ and discriminant $D\ne 0$. Then
    \begin{itemize}
        \item[(i)] $n\le c_1(|D|)$, and
        \item[(ii)] There is an algebraic integer $\beta$, $\mathbb{Z}$-equivalent to $\alpha$ such that
        $$H(\beta)\le c_2(n,|D|),$$
        where $c_1,c_2$ denote the same effectively computable positive numbers as in Theorem \ref{thm:3.2}.
    \end{itemize}
\end{theorem}

This theorem was stated and proved in Gy\H ory (1973) as 'Corollaire 3' of the 'Th\'eor\`eme', cf. Theorem \ref{thm:3.2} above. 

As was mentioned in Section \ref{sec:1}, the cubic case was settled independently by Delone (1930) and Nagell (1930), and the quartic case by Nagell (1967) in an ineffective way. 

Theorems \ref{thm:5.1} resp. \ref{thm:5.2} confirmed in full generality, and in fact Theorem \ref{thm:5.2} in an effective form, a conjecture of Nagell (1967). Further both Theorem \ref{thm:5.1} and Theorem \ref{thm:5.2} imply, Theorem \ref{thm:5.2} in an effective form, that there are only finitely many algebraic units in $\overline{\mathbb{Q}}$ of given discriminant. This gave the effective solution to Problem 19 in the book Narkiewicz (1974).

Finally, we note that Theorem \ref{thm:5.2} easily follows from Theorem \ref{thm:3.2}. Indeed, if $\alpha$ is an algebraic integer with the properties specified in Theorem \ref{thm:5.2}, then by \eqref{eq:4.1}, $D(f_\alpha)=D$ and $\deg f_\alpha=n$. Further, by Theorem \ref{thm:3.2} $f_\alpha$ is $\Zz$-equivalent to some monic $g\in\Zz[X]$ of degree $n$ and discriminant $D$ such that $n\le c_1(|D|)$ and $H(g)\le c_2(n,|D|)$, where $c_1,c_2$ denote the effectively computable numbers occurring in Theorem \ref{thm:3.2}. But then $\alpha$ is $\Zz$-equivalent to a zero of $g$, say $\beta$, whence $\deg\beta\le c_1(|D|)$ and $H(\beta)\le c_2(n,|D|)$ follow.\qed
\\[0.1cm]

The first explicit version of (ii) in Theorem \ref{thm:5.2} was established by Gy\H ory (1974) by means of 
Baker's method. For $c_1(n)$ one can take $2\log |D|/\log 3$. For $c_2(n,|D|)$ we can obtain an explicit bound, using Theorem \ref{thm:3.8} instead of Theorem \ref{thm:3.2}. An even better explicit estimate can be obtained in (ii), observing that in fact we apply Theorem \ref{thm:3.2} (or its explicit version Theorem \ref{thm:3.8}) only to irreducible polynomials $f_\alpha$. The best known bound in (ii) comes from Theorem 6.4.1 of Evertse and Gy\H ory (2017). 

\subsection{Consequences 
for monogenic number fields and orders}\label{sec:5.3new}
\vphantom{a}

\noindent
Let again $K$ be a number field of degree $n\ge 2$ with ring of integers $\mathcal{O}_K$ and discriminant $D_K$.

The following effective corollaries are immediate consequences of Theorem \ref{thm:5.2} (i.e. the `Corollaire 3') of Gy\H ory (1973). Although this was not mentioned by Birch and Merriman in their 1972 paper, it should be remarked that from their Theorem \ref{thm:5.1} one can also deduce in ineffective form the finiteness consequences of the results below.

\begin{corollary}[of Theorem \ref{thm:5.2}]
\label{cor:5.3}
    Let $\mathcal{O}$ be an order of $K$ and $D$ a non-zero integer. Every $\alpha$ in $\mathcal{O}$ of discriminant $D_{K/\mathbb{Q}}(\alpha)=D$ is $\mathbb{Z}$-equivalent to some $\beta\in\mathcal{O}$ such that 
    $$H(\beta)\le c_2(n,|D|),$$
    where $c_2=c_2(n,|D|)$ denotes the same effectively computable positive number as in Theorem \ref{thm:5.2}.
\end{corollary}

This is a special case of Theorem \ref{thm:5.2},
restricted to the 
elements of $\mathcal{O}$. It follows from Corollary \ref{cor:5.3} that up to $\mathbb{Z}$-equivalence, there are only finitely many elements of $\mathcal{O}$ of given non-zero discriminant, and all of them can be, at least in principle, effectively determined. 

As was mentioned above, the first quantitative versions of Theorem \ref{thm:5.2} and Corollary \ref{cor:5.3} were established in Gy\H ory (1974).

\begin{corollary}[of Theorem \ref{thm:5.2}]
 \label{cor:5.4}
    Let $\mathcal{O}$ be an order in $K$ of discriminant $D(\OO )$, and $I_\mathcal{O}$ a positive integer. Every $\alpha$ in $\mathcal{O}$ with index $I_\mathcal{O}(\alpha)=I_\mathcal{O}$ is $\mathbb{Z}$-equivalent to some $\beta\in\mathcal{O}$ such that
    $$H(\beta)\le c_2(n,I^2_\mathcal{O}\cdot |D(\OO )|),$$
    where $c_2$ denotes the same effectively computable positive number as in Theorem \ref{thm:5.2} with $|D|$ replaced by $I^2_\mathcal{O}\cdot |D(\OO)|$.
\end{corollary}

This follows immediately from Corollary \ref{cor:5.3} and the second identity in \eqref{eq:4.2}. Corollary \ref{cor:5.4} implies that up to $\mathbb{Z}$-equivalence there are only finitely many elements in $\mathcal{O}$ with given index and all of them can be, at least in principle, effectively determined.

The next Theorem \ref{cor:5.6} and its more general version Theorem \ref{thm:5.5} are the most influential consequences of Theorem \ref{thm:5.2}. They provided the first general effective algorithm for deciding the monogenicity, the multiplicity of monogenicity, and for determining, at least in principle, all power integral bases in $K$  and in its orders.

Of particular importance are the cases when in Corollaries \ref{cor:5.3}, and \ref{cor:5.4}  $\mathcal{O}$ is just $\mathcal{O}_K$, the ring of integers of $K$. Then Corollary \ref{cor:5.4} implies 

\begin{theorem} [Gy\H{o}ry, 1976]\label{cor:5.6}
Every $\alpha\in \mathcal{O}_K$ with $\mathcal{O}_K=\mathbb{Z}[\alpha]$ is $\mathbb{Z}$-equivalent to some $\beta\in \mathcal{O}_K$ such that
    $$H(\beta)\le c_2(n,|D_K|),$$
    where $c_2$ denotes the same effectively computable positive number as in Corollary \ref{cor:5.4} with $I_{\OO}=1$, $D(\OO )=D_K$.
Consequently, there are only finitely many $\mathbb{Z}$-equivalence classes of $\alpha$ in $\mathcal{O}_K$ such that $\mathcal{O}_K = \mathbb{Z}[\alpha]$, and a full set of representatives of these classes can be, at least in         principle, effectively found.
\end{theorem}

More generally, Corollary \ref{cor:5.4} immediately gives the following:  

\begin{theorem}[Gy\H{o}ry, 1976]
 \label{thm:5.5}
    Let $\mathcal{O}$ be an order of $K$ of discriminant $D(\OO )$. Every $\alpha\in\mathcal{O}$ with $\mathcal{O}=\mathbb{Z}[\alpha]$ is $\mathbb{Z}$-equivalent to some $\beta\in \mathcal{O}$ such that
    $$H(\beta)\le c_2(n,|D(\OO )|),$$
    where $c_2$ denotes the same effectively computable positive number as in Corollary \ref{cor:5.4} with $I_{\OO}=1$.
\end{theorem}


 The first explicit, quantitative versions of Corollary \ref{cor:5.4} and 
 Theorems \ref{cor:5.6} and \ref{thm:5.5} were given in Gy\H ory (1976).

\begin{rem}
	With the above formulation of Corollaries \ref{cor:5.3}, \ref{cor:5.4} and Theorem \ref{thm:5.5} it was easier to point out that these are indeed consequences of Theorems \ref{thm:3.2} and \ref{thm:5.2}. Further, we note that their explicit versions can be easily derived from the explicit variant Theorem \ref{thm:3.8} of Theorem \ref{thm:3.2}. Finally, the corollaries can be deduced with better bounds from less general versions of Theorem \ref{thm:3.2}, where the polynomials $f$ involved are irreducible; for such versions we refer to Gy\H ory (1976, 1998, 2000), Evertse and Gy\H ory (2017) and in fact Corollary \ref{cor:5.3} above.
\end{rem}

\subsection{Reformulation of Corollaries \ref{cor:5.3}, \ref{cor:5.4} and Theorem \ref{thm:5.5} over $\mathcal{O}_K$ in terms of polynomial Diophantine equations over $\mathbb{Z}$}\label{sec:4.3}
\vphantom{a}

\noindent
Let $K$ be an algebraic number field of degree $n\ge 2$ with ring of integers $\mathcal{O}_K$ and discriminant $D_K$. Consider Corollaries \ref{cor:5.3}, \ref{cor:5.4} and Theorem \ref{thm:5.5} for $\mathcal{O}=\mathcal{O}_K.$ Let $\{1,\omega_2,\ldots,\omega_n\}$ be an integral basis of $K$. For $\alpha\in\mathcal{O}_K$ with
$$\alpha=x_1+x_2\omega_2+\cdots+x_n\omega_n,\hspace{3mm}x_1,x_2,\ldots,x_n\in\mathbb{Z},$$
its discriminant
\setcounter{equation}{4}
\begin{align}\label{eq:4.3new}
    D(\alpha) = D(x_2\omega_2+\cdots+x_n\omega_n)
\end{align}
can be regarded as a decomposable form of degree $n(n-1)$ in $x_2,\ldots,x_n$ with coefficients in $\mathbb{Z}$, i.e.,
it is a product of $n(n-1)$ linear forms in $x_2\kdots x_n$ with algebraic coefficients. The form $D(x_2\omega_2+\cdots+x_n\omega_n)$, which was introduced by Kronecker (1882), is called \textit{discriminant form}, while, for $D\ne 0$, the equation
\begin{align}\label{eq:4.4new}
    D(x_2\omega_2+\cdots+x_n\omega_n)=D\text{ in } x_2,\ldots,x_n\in\mathbb{Z}
\end{align}
is called a \textit{discriminant form equation}.

Clearly, Corollary \ref{cor:5.3} implies the following
\begin{corollary}[of Theorem \ref{thm:5.2}]\label{cor:5.4new}
    For given $D\ne 0$, the discriminant form equation \eqref{eq:4.4new} has only finitely many solutions and they can be effectively determined. 
\end{corollary}

The following important fact is due to {Hensel} (1908):
\\
\textit{to the integral basis $\{1,\omega_2,\ldots,\omega_n\}$ of $K$ there corresponds a decomposable form $I(X_2,\ldots,X_n)$ of degree $n(n-1)/2$ in $n-1$ variables with coefficients in $\Zz$ such that for
$\alpha\in\mathcal{O}_K$}
\begin{align}\label{eq:4.3}
	&I(\alpha)=|I(x_2,\ldots,x_n)|\text{ \textit{if} } \alpha=x_1+x_2\omega_2+\cdots+x_n\omega_n\\
	&\text{ \textit{with} } x_1,x_2,\ldots,x_n\in\Zz.\notag
\end{align}
\textit{Here $I(X_2,\ldots,X_n)$ is called an index form, and for given non-zero $I\in\Zz$,}
\begin{align}\label{eq:4.4}
	I(x_2,\ldots,x_n)=\pm I\text{ \textit{in} } x_2,\ldots,x_n\in\Zz
\end{align}
\textit{an} \textit{index form equation}. 

We note that the equations \eqref{eq:4.4new} and \eqref{eq:4.4} are related by the first 
identity in \eqref{eq:4.2}.

In view of \eqref{eq:4.3}, the finiteness assertion of Corollary \ref{cor:5.4} for $\mathcal{O}_K$ implies the following.

\begin{corollary}[of Theorem \ref{thm:5.2}]
\label{cor:4.4}
	For given $I\in\Zz\setminus\{0\}$, the index form equation \eqref{eq:4.4} has only finitely many solutions, and they can be effectively determined.
\end{corollary}

In particular, for $I=1$, we get the following consequence of Theorem \ref{cor:5.6}.

\begin{corollary}[of Theorem \ref{thm:5.2}] 
\label{cor:4.5}
	The index form equation
	\begin{align}\label{eq:4.5}
		I(x_2,\ldots,x_n)=\pm 1\text{ \textit{in} } x_2,\ldots,x_n\in\Zz
	\end{align}
	has only finitely many solutions, and they can be effectively determined.
\end{corollary}

Corollaries \ref{cor:5.4new}, \ref{cor:4.4} and \ref{cor:4.5} were proved in {Gy\H ory} (1976) with explicit upper bounds for the solutions, not only for equations \eqref{eq:4.4new}, \eqref{eq:4.4} and \eqref{eq:4.5} but also for
index form equations related to indices with respect to arbitrary orders $\OO$ of $K$; see also Gy\H ory (2000) and Evertse and Gy\H ory (2017).

The best known upper bound for the solutions of \eqref{eq:4.5} is
\begin{align}\label{eq:4.6}
	\max_{2\le i\le n}|x_i|<\exp\{10^{n^2}(|D_K|(\log |D_K|)^n)^{n-1}\}
\end{align}
which is due to {Evertse} and {Gy\H ory} (2017).

\subsection{A consequence of Theorem \ref{thm:3.4} for algebraic numbers of given discriminant}\label{sec:5.4}
\vphantom{a}

\noindent
Theorem \ref{thm:3.4} can be applied to algebraic numbers that are not necessarily algebraic integers.  Given an algebraic number $\alpha$, we denote by $f_{\alpha}$ its \emph{primitive minimal polynomial}, i.e.,
\begin{align}\label{eq:4.7}
f_{\alpha}=a_0X^n+\cdots +a_n=a_0(X-\alpha^{(1)})\cdots (X-\alpha^{(n)})\in \Zz [X]
\end{align}
where $a_0>0$, $\gcd (a_0\kdots a_n)=1$ and $\alpha^{(1)}=\alpha ,\kdots \alpha^{(n)}$ are the distinct conjugates of $\alpha$. We recall that the height and discriminant of $\alpha$ are defined by those of $f_{\alpha}$,
i.e.,
\[
H(\alpha ):=H(f_{\alpha}),\ \ D(\alpha ):=D(f_{\alpha}).
\]
Two algebraic numbers $\alpha,\beta$ are called $GL_2(\mathbb{Z})\textit{-equivalent}$ if
\begin{align*}
    \beta=\frac{a\alpha+b}{c\alpha+d}\text{\,\, with \,}\begin{pmatrix}
        a&b\\c&d
    \end{pmatrix}\in GL_2(\mathbb{Z}).
\end{align*}
One easily verifies that
if $\alpha ,\beta$ are $GL_2(\Zz )$-equivalent then so are $f_{\alpha}$, $f_{\beta}$ while conversely,
if $f_{\alpha}, f_{\beta}$ are $GL_2(\Zz )$-equivalent, then $\alpha$ is $GL_2(\Zz )$-equivalent to a conjugate of $\beta$.

Consequently, if $\alpha$, $\beta$ are $GL_2(\Zz )$-equivalent, then $D(\alpha)=D(\beta)$. 
Now Theorem \ref{thm:3.4} implies at once:

\begin{theorem}[Evertse and Gy\H ory, 1991a]\label{thm:5.6}
    Every algebraic number $\alpha$ of degree $n\ge 2$ and discriminant $D\ne 0$ is $GL_2(\mathbb{Z})$-equivalent to an algebraic number $\beta$ with
    $$H(\beta)\le c_3(n,|D|),$$
    where $c_3$ denotes the same effectively computable positive number as in Theorem \ref{thm:3.4}.
\end{theorem}

Further, by Th\'eor\`eme 1 of Gy\H ory (1974) we have 
$$n\le 2\log|D|/\log 3.$$

\subsection{Rationally monogenic orders}\label{sec:5.5}
\vphantom{a}

\noindent
Monogenic orders $\Zz [\alpha ]$, where $\alpha$ is an algebraic integer, can be generalized to so-called \emph{rationally monogenic orders} $\Zz_{\alpha}$, where $\alpha$ is not necessarily integral.
We will formulate an analogue of Corollary \ref{cor:5.3} for rationally monogenic orders.
While in the results for monogenic orders, $\Zz$-equivalence of algebraic integers plays an important role, for rationally monogenic orders we have to deal with $GL_2(\Zz )$-equivalence of algebraic numbers. 
Before we define rationally monogenic orders, we briefly go into some history and introduce the necessary terminology.

Let $\alpha$ be a non-zero, not necessarily integral algebraic number of degree $n\geq 3$,
and $f_{\alpha}$ its primitive minimal polynomial, given by \eqref{eq:4.7}.
Define $\Zz_{\alpha}$ to be the $\Zz$-module with basis
$$
1,\ \omega_2:=a_0\alpha,\ \omega_3:=a_0\alpha^2 +a_1\alpha\kdots
\omega_n:=a_0\alpha^{n-1}+a_1\alpha^{n-2}+\cdots +a_{n-2}\alpha .
$$
This $\Zz$-module was introduced by Birch and Merriman (1972),
who observed that it is contained in the ring of integers of $\Qq (\alpha )$, and that for its discriminant we have
\begin{equation}\label{eq:4.8}
D(\Zz_{\alpha})=D(f_{\alpha})=D(\alpha ).
\end{equation}
Nakagawa (1989) showed that $\Zz_{\alpha}$ is in fact an \textit{order}
of the field $\Qq (\alpha )$, i.e.,  closed under multiplication. More precisely, he showed that 
\begin{equation}\label{eq:4.9}
\omega_i\omega_j=-\sum_{\max (i+j-n,1)\leq k\leq i} a_{i+j-k}\omega_k+
\sum_{j<k\leq \min (i+j,n)} a_{i+j-k}\omega_k
\end{equation}
for $i,j=1\kdots n-1$, where $\omega_n:=-a_n$.
This order was further studied by 
Simon (2001, 2003) and Del Corso, Dvornicich and Simon (2005).
They showed that
\begin{equation}\label{eq:4.9a}
\Zz_{\alpha }=\Zz [\alpha ]\cap \Zz [\alpha^{-1}].
\end{equation}
As was very likely known at the time, another description of $\Zz_{\alpha}$ is as follows. Let $\mathcal{M}_{\alpha}$ be the $\Zz$-module generated by $1,\alpha\kdots \alpha^{n-1}$. Then
$\Zz_{\alpha}$ is the \emph{ring of coefficients} of $\MM_{\alpha}$ (see Borevich and Shafarevich (1967), Section 2.2), i.e.,
\begin{equation}\label{eq:4.9b}
\Zz_{\alpha}=\left\{ \xi\in\Qq (\alpha ):\ \xi\mathcal{M}_{\alpha}
\subseteq\mathcal{M}_{\alpha}\right\}.
\end{equation}
We have
\begin{equation}\label{eq:4.9c}
\Zz_{\alpha} =\Zz [\alpha ]\ \ \text{if $\alpha$ is an algebraic integer.}
\end{equation}
Indeed, if $\alpha$ is an algebraic integer of degree $n$, the powers $\alpha^i$ ($i\ge n$) belong to $\mathcal{M}_{\alpha}$, and thus, $\Zz_{\alpha}=\mathcal{M}_{\alpha}=\Zz [\alpha ]$.
Further, for any two non-zero algebraic numbers $\alpha ,\beta$ we have
\begin{equation}\label{eq:4.9d}
\text{$\alpha$, $\beta$\ \,$GL_2(\Zz )$-equivalent}\ \Longrightarrow\ \Zz_{\alpha}=\Zz_{\beta}.
\end{equation}
Indeed, let $\beta =\medfrac{a\alpha +b}{c\alpha +d}$
for some matrix $\big(\begin{smallmatrix}a&b\\ c&d\end{smallmatrix}\big)\in GL_2(\Zz )$.
Then $\mathcal{M}_{\beta}=(c\alpha +d)^{1-n}\mathcal{M}_{\alpha}$ where $n=\deg\alpha$, and thus, $\Zz_{\alpha}=\Zz_{\beta}$.

We call an order $\mathcal{O}$ of a number field $K$
\textit{rationally monogenic} if there is $\alpha$
such that $\mathcal{O}=\Zz_{\alpha}$.
From \eqref{eq:4.9c} it follows that monogenic orders are rationally monogenic.
Below we explain that rationally monogenic orders are in fact special cases of invariant orders
of polynomials.
In particular, $\Zz_{\alpha}$ is the invariant order of $f_{\alpha}$. 

Recall that the index of an algebraic integer was defined in \eqref{eq:4.index}.
Following Simon (2001), we generalize this to not necessarily integral algebraic numbers as follows. Given a non-zero algebraic number $\alpha$, we define the index of $\alpha$ by
\[
I(\alpha ):=[\OO_K : \Zz_{\alpha}],
\]
where $K=\Qq (\alpha )$. In fact, this is the index of $f_{\alpha}$ as it was introduced by Simon. From \eqref{eq:4.x1} and \eqref{eq:4.8} we deduce, analogously to the first identity in \eqref{eq:4.2},
\begin{equation}\label{eq:4.generalindex}
D(\alpha )=I(\alpha)^2D_K.
\end{equation}
For more results and properties of this index, we refer to Simon (2001).

There is a connection between rationally monogenic orders and Hermite  equivalence classes of polynomials, which we explain here without proof. For a non-zero algebraic number $\alpha$,
let $\II_{\alpha}$ be the fractional ideal of $\Zz_{\alpha}$ generated by $1$ and $\alpha$. This is known to be invertible, see Simon (2003). It is called also the \emph{invariant ideal} of $f_{\alpha}$. 

\begin{theorem}[BEGyRS, 2023]\label{thm:B5}
Let $f,g\in \Zz [X ]$ be two primitive, irreducible polynomials. 
Then the following three assertions are equivalent:
\begin{itemize}
\item[(i)] $f$ and $g$ are Hermite equivalent;
\item[(ii)] $f$ has a root $\alpha$ and $g$ a root $\beta$ such that $\MM_{\beta}=\lambda\MM_{\alpha}$
for some non-zero $\lambda\in\Qq (\alpha )$;
\item[(iii)] $f$ has a root $\alpha$ and $g$ a root $\beta$ such that
$\Zz_{\alpha}=\Zz_{\beta}$ and $\II_{\alpha}$ and $\II_{\beta}$
lie in the same ideal class of $\Zz_{\alpha}$.
\end{itemize}
\end{theorem}

In the particular case that $f$ and $g$ are monic, we have $\alpha\in \Zz [\alpha ]=\Zz_{\alpha}$
and $\II_{\alpha}=\Zz_{\alpha}$ and likewise for $g$ and $\beta$. This leads to the
following corollary.

\begin{corollary}[BEGyRS, 2023]\label{cor:B5}
Let $f,g\in \Zz [X ]$ be two monic, irreducible polynomials.
Then $f$ and $g$ are Hermite equivalent if and only if $f$ has a root $\alpha$ and $g$ a root $\beta$ such that $\Zz [\alpha ]=\Zz [\beta ]$.
\end{corollary}

In BEGyRS (2023) an example of a quartic algebraic number field $K$ was given such that $\OO_K =\Zz_{\alpha}=\Zz_{\beta}$ for certain $\alpha ,\beta\in K$, but $f_{\alpha}$, $f_{\beta}$ lie
in different Hermite equivalence classes. So far, we haven't been able to find similar examples for algebraic number fields of degree $\geq 5$.

An order $\OO$ of a number field $K$ is called \emph{primitive} if there are no integer $a>1$ and order $\OO'$ such that $\OO =\Zz +a\OO'$. It is not difficult to show that a rationally monogenic order is primitive. It follows from work of Delone and Faddeev (1940) that every primitive order of a cubic number field is rationally monogenic.
Simon (2001) gave various examples of number fields of degree $\geq 4$ that are not rationally monogenic, i.e., whose rings of integers
are not rationally monogenic.

In Evertse (2023) the following was shown:

\begin{theorem}\label{thm:5.6a}
Every number field $K$
of degree $\ge 3$ has infinitely many orders that are rationally monogenic but not monogenic.
\end{theorem}

We finally arrive at the main result of this subsection, which follows directly from
Theorem \ref{thm:5.6} and \eqref{eq:4.8}:

\begin{theorem}\label{thm:5.7}
Let $\OO$ be an order of a number field $K$, and 
denote by $D(\OO )$ its discriminant.
Then every $\alpha$ such that $\Zz_{\alpha}=\OO$ is $GL_2(\Zz )$-equivalent
to some $\beta\in K$ of height $H(\beta )\le c_3(n,|D(\OO )|)$, where $c_3$ denotes the same effectively computable positive number as in Theorem \ref{thm:5.6}. 
\end{theorem}

This implies

\begin{theorem}\label{cor:5.8}
Let $\OO$ be an order of a number field $K$.
Then it can be effectively decided whether there is $\alpha$ such that $\OO =\Zz_{\alpha}$.
\\
Moreover, there are only finitely many $GL_2(\Zz )$-equivalence classes of
$\alpha\in K$ such that $\Zz_{\alpha}=\mathcal{O}$,
and a full system of representatives of those can be effectively determined. 
\end{theorem}

\begin{proof}[Idea of proof]
Suppose $K$ is effectively given in the form
$\Qq (\gamma )$, with an algebraic number $\gamma$ of degree $n$.
Thus, each element of $K$ has a represention as a $\Qq$-linear combinations of $1,\gamma\kdots \gamma^{n-1}$,
and we can express all computations on $K$
in terms of such representations.

Let the order $\OO$ be given by a $\Zz$-module basis $1,\theta_2\kdots \theta_n$ (with representions as described above). Using Theorem \ref{thm:5.6},
one can effectively determine a full system of representatives for the $GL_2(\Zz )$-equivalence classes of those $\alpha\in K$ with $D(\alpha )=D(\OO )$. To check whether such a representative $\alpha$ satisfies $\Zz_{\alpha}=\OO$, one can proceed as follows.
Verify that $\OO \subseteq \Zz_{\alpha}$ by checking $\theta_i\MM_{\alpha}\subseteq \MM_{\alpha}$ for $i=2\kdots n$.
If so, we have in fact $\OO =\Zz_{\alpha}$ by \eqref{eq:4.8}.
\end{proof}

The rationally monogenic orders introduced above are in fact special cases of \emph{invariant orders} or \emph{invariant rings} of binary forms, for which there is now a vast general theory. 
Although outside the scope of this paper, we give some background on these rings.

Let $A$ be a commutative ring (with $1$),
and $a_0\kdots a_n\in A$. Then the \emph{invariant ring} (order if $A=\Zz$) associated with $(a_0\kdots a_n)$, or rather with the binary form $F(X,Y)=a_0X^n+\cdots +a_nY^n$ (but we allow here that $a_0=0$ or even $a_0=\cdots =a_n=0$)
is given by the $A$-algebra $A_F$ with $A$-module basis $1,\omega_2\kdots \omega_n$ satisfying the multiplication table \eqref{eq:4.9}. This is in fact a commutative, associative $A$-algebra.
The name `invariant ring' (invariant order if $A=\Zz$) comes from the following invariance property: if
$F,G$ are two $GL_2(A)$-equivalent binary forms, i.e., $G(X,Y)= uF(aX+bY,cX+dY)$ for some 
$u\in A^*$, 
$\big(\begin{smallmatrix}a&b\\ c&d\end{smallmatrix}\big)\in GL_2(A)$, then $A_G\cong A_F$ as $A$-algebras. Any ring that is the invariant ring of a binary form is called a binary ring (or binary order if $A=\Zz$).

Thus, $\Zz_{\alpha}$ is the invariant order of $F_{\alpha}(X,Y):=Y^nf_{\alpha}(X/Y)$. In other words, a rationally monogenic order is the invariant order of a primitive, irreducible binary form.

From work of Delone and Faddeev (1940), later extended by Gan, Gross, and Savin (2002) and Deligne (unpublished) (see also section 16.3 of Evertse and Gy\H{o}ry (2017)) it follows that for every commutative ring $A$, the map $F\mapsto A_F$ gives a one-to-one correspondence between $GL_2(A)$-equivalence classes of binary cubic forms in $A[X,Y]$ and isomorphism classes of free cubic $A$-algebras, i.e., commutative, associative, unital $A$-algebras that as an $A$-module are free of rank $3$. Wood (2011) gave a geometric interpretation of invariant rings of binary forms.

\section{Algorithmic resolution of index form equations, application to (multiply) monogenic number fields}\label{sec:5}

As above, $K$ will denote a number field of degree $n\ge 3$ with ring of integers $\mathcal{O}_K$ and discriminant $D_K$. For an index form $I(X_2,\ldots,X_n)$ associated with an integral basis $\{1,\omega_2,\ldots,\omega_n\}$ of $K$, consider again the above \textit{index form equation \eqref{eq:4.5}}.

The exponential bound \eqref{eq:4.6} for the solutions of \eqref{eq:4.5} is too large for practical use. In the 1990's, there were new breakthroughs, leading to the complete resolution of certain index form equations.
In fact,  practical methods were elaborated for solving equation \eqref{eq:4.5} when $|D_K|$ is not too large, and the degree $n$ of $K$ is $\le 6$.
Further, \eqref{eq:4.5} was solved
for some special {higher degree number fields} $K$ up to about degree $15$ and for some relative extensions of degree $\le 4$.

\subsection{The case $\mathbf{n=3}$ and $\mathbf{4}$. Approach via Thue equations of degree $3$ and $4$.}\label{sec:5.1}
\vphantom{a}

\noindent
As will be seen, in the case $n=3$ equation \eqref{eq:4.5} can be reduced to a cubic Thue equation while, in the case $n=4$, to a cubic and some quartic Thue equations, that is to equations of the form
\begin{align}\label{eq:6.1}
    F(x,y)=m\text{ in } x,y\in\mathbb{Z},
\end{align}
where $m$ is a non-zero integer and $F\in \mathbb{Z}[X,Y]$ is a binary form of degree $3$ or $4$ with pairwise non-proportional linear factors over $\overline{\mathbb{Q}}$. By a general theorem of Thue (1909), every equation of the type \eqref{eq:6.1} with degree $\ge 3$ has only finitely many solutions, and Baker (1968b) gave an explicit upper bound for their solutions in terms of $|m|$ and the height and degree of $F$. The best known bound is due to Bugeaud and Gy\H ory (1996b). However, in concrete cases this bound is too large for practical use. For solving concrete Thue equations, general practical methods were developed in Peth\H o and Schulenberg (1987) for $m=1$, and in Tzanakis and de Weger (1989) for arbitrary $m$. Later, these methods were made even more efficient in Bilu and Hanrot (1996, 1999) and Hanrot (1997). Their algorithms are based on Baker's method and certain reduction and enumeration techniques. Hence we possess efficient algorithms for solving equation \eqref{eq:4.5} for $n=3$ and $4$. However, this approach cannot be applied in general to index form equations in number fields $K$ of degree $n>4$, except for $n=6,8,9$ when $K$ has a quadratic or cubic subfield; then equation \eqref{eq:4.5} leads to relative cubic or quartic Thue equations.

For $\mathbf{n=3}$, {Ga\'al} and {Schulte} (1989) reduced equation \eqref{eq:4.5} to a cubic Thue equation with $m=1$. Then, using the algorithm elaborated for solving cubic Thue equations, they determined all power integral bases of cubic fields $K$ with discriminant $-300\le D_K\le 3137$. Their computations were later extended in Schulte (1989, 1991). 

For $\mathbf{n=4}$, {Gaál}, {Peth\H o} and {Pohst} (1993, 1996) first reduced the  equation \eqref{eq:4.5} to a cubic Thue equation and a pair of ternary quadratic equations. Then the quadratic equations were themselves reduced to quartic Thue equations. Finally, by means of efficient algorithms for solving such Thue equations, they computed the solutions of equation \eqref{eq:4.5} for quartic number fields with not too large discriminant.  They obtained several interesting tables on the distribution of minimal indices and about the average behaviour of minimal indices.

\subsection{The cases $\mathbf{n=5}$ and $\mathbf{6}$. Refined version of the general approach via unit equations, combined with reduction and enumeration algorithms}\label{sec:5.2}
\vphantom{a}

\noindent
For $\mathbf{n\ge 5}$, the approach via {Thue equations} does not work, in general. 
For $\mathbf{n=5}$ and $\mathbf{6}$ a refined version of the general approach involving unit equations is needed. Since by \eqref{eq:4.3}, \eqref{eq:4.2} and \eqref{eq:4.1} we have
for $\alpha\in\OO_K$ with $K=\Qq (\alpha )$,
\[
 \eqref{eq:4.5}\Leftrightarrow D(\alpha)=D_K\Leftrightarrow D(f_\alpha)=D_K\ \text{in\ }\alpha\in\mathcal{O}_K
\]
 where $f_\alpha\in\Zz[X]$ is the minimal polynomial of $\alpha$, in case of concrete equations \eqref{eq:4.5} a refinement of the proof of Theorem \ref{thm:3.2} for irreducible $f_\alpha$'s must be combined with some reduction and enumeration algorithms.

The refined version of the general method for solving index form equations \eqref{eq:4.5} consists of the following steps:
\begin{itemize}
	\item[{\bf 1.}] Reduction to unit equations but in considerably smaller subfields of the normal closure $G$ of $K$, of which the unit rank is much smaller than that of $G$, i.e., at most $n(n-1)/2-1$ (note that the unit rank of $G$ may be as large as $n!-1$); cf. {Gy\H ory} (1998, 2000). Then in the unit equation corresponding to \eqref{eq:3.5}, one can write 
	$\varepsilon_{ijk}=\zeta_{ijk}\rho_1^{a_{ijk,1}}\cdots \rho_r^{a_{ijk,r}}$, 
	with a root of unity $\zeta_{ijk}$ and a fundamental system of units $\rho_1\kdots\rho_r$ of bounded height, and
in concrete cases one can bound the exponents $|a_{ijk,l}|$ by {Baker's method}. Here the estimate of Baker and W\"ustholz (1993) for linear forms in logarithms of algebraic numbers is very practical to apply in calculations.
	\item[{\bf 2.}] The bounds in concrete cases are still too large. Hence a \textit{reduction algorithm} is needed, {reducing} the {Baker's bound} for $|a_{ijk,l}|$ in several steps if necessary by a refined version of the \textit{$L^3$-algorithm}; cf. {de Weger} (1989), Tzanakis and de Weger (1989), {Wildanger} (1997) and {Ga\'al} and {Pohst} (1996).
	\item[{\bf 3.}] The last step is to apply an \textit{enumeration algorithm}, determining the small solutions {under the reduced bound}; cf. {Wildanger} (1997, 2000), {Ga\'al} and {Gy\H ory} (1999) and {Bilu}, {Ga\'al} and {Gy\H ory} (2004).
\end{itemize} 

Combining the refined version of the general approach with reduction and enumeration algorithms, for $n=5,6$ and for not too large $|D_K|$, {Ga\'al} and {Gy\H ory} (1999), resp. {Bilu}, {Ga\'al} and {Gy\H ory} (2004) gave algorithms for determining all power integral bases and hence checking the monogenicity and determining the multiplicity of the monogenicity of $K$.

We note that the use of the refined version of the general approach is particularly important in the application of the enumeration algorithm.

To perform computations, {algebraic number theory packages}, a {computer algebra system} and in some cases a {supercomputer} were needed.

\subsection{Examples: resolutions of index form equations of the form \eqref{eq:4.5} for $\mathbf{n=3,4,5,6}$ in the most difficult case}\label{sec:5.3}
\vphantom{a}

\noindent
In the examples below, the authors resolved concrete index form equations of the form \eqref{eq:4.5} for $n=3,4,5,6$. The number fields $K$ of degree $n$ are given by irreducible monic polynomials $f(X)\in\Zz[X]$, a zero of which generates the corresponding $K$ over $\Qq$. In each case all power integral bases in $K$, and therefore the \textit{multiplicity of the monogenicity of $K$}, denoted by $mm(K)$, are computed by the method outlined above. For the lists of the power integral bases, we refer to the original papers and to {Evertse} and {Gy\H ory} (2017) and {Ga\'al} (2019).
\begin{itemize}
	\item[] $\mathbf{n=3}$, $f(X)=X^3-X^2-2X+1$,  $mm(K)=9$ ({Ga\'al} and {Schulte}, 1989);
 \\
	\item[] $\mathbf{n=4}$, $f(X)=X^4-4X^2-X+1$, $mm(K)=17$ ({Ga\'al}, {Peth\H o} and {Pohst}, 1990's);
 \\
	\item[] $\mathbf{n=5}$, $f(X)=X^5-5X^3+X^2+3X-1$, $mm(K)=39$ ({Ga\'al} and {Gy\H ory}, 1999); 
 \\
	\item[] $\mathbf{n=6}$, $f(X)=X^6-5X^5+2X^4+18X^3-11X^2-19X+1$, $mm(K)=45$ ({Bilu}, {Ga\'al}, and {Gy\H ory}, 2004). 
 \end{itemize}
 
We note that from the point of view of computation, the above examples belong to the most difficult cases for $n=3,4,5,$ and $6$, $K$ being in each case totally real with Galois group $S_n$. In these cases the number of exponents in the unit equations involved is the largest possible.
\\[0.2cm]
\textbf{Remark.} The general procedure outlined above to solve any concrete equation \eqref{eq:4.5} for $n=6$ requires considerable CPU-time. In certain special cases (e.g., if $n=6$ and $K$ has a quadratic subfield), there are faster algorithms, see Ga\'{a}l (2024, 2025).
However, some of these algorithms determine only the ``small'' solutions, and do not exclude the existence
of ``large" solutions.

For $n\ge 7$, the above mentioned algorithms do not work in general. Then the number of fundamental units, $\rho_1,\ldots,\rho_r$ involved can be $\ge \frac{7\cdot 6}{2}-1=20$ which is too large to use the enumeration algorithm.

\begin{prob}
	For $\mathbf{n= 7,8,\ldots}$, give a practical algorithm for solving equation \eqref{eq:4.5} in case of \textbf{any} number field $K$ of degree $\mathbf{n}$ with not too large discriminant.
\end{prob}

    \section{Power integral bases and canonical number systems in number fields}\label{sec:6}

Number systems and their generalizations have been intensively studied for a long time. Here we present an important generalization for the number field case, point out its close connection with power integral bases and formulate an application of the above Theorem \ref{cor:5.6} to this generalization.

Let $K$ be an algebraic number field with ring of integers $\mathcal{O}_K$, and let $\alpha\in\mathcal{O}_K$ with $|N_{K/\mathbb{Q}}(\alpha)|\ge 2$. Then $\{\alpha,\mathcal{N}(\alpha)\}$ with
$$\mathcal{N}(\alpha)=\{0,1,\ldots,|N_{K/\mathbb{Q}}(\alpha)|-1\}$$
is called a \textit{canonical number system}, in short CNS, in $\mathcal{O}_K$, if every non-zero element of $\mathcal{O}_K$ has a unique representation of the form
$$a_0+a_1\alpha+\cdots+a_k\alpha^k\text{ with } a_i\in\mathcal{N}(\alpha)\text{ for } i=0,\ldots,k,\, a_k\ne 0.$$
Then $\alpha$ is called the \textit{base} and $\mathcal{N}(\alpha)$ the \textit{set of digits} of the number system. This concept is a generalization of the radix representation considered in $\mathbb{Z}$.

B. Kov\'acs (1981) proved the following fundamental theorem.

\begin{theorem}[B. Kov\'acs, 1981]\label{thm:6.1}
In $\mathcal{O}_K$ there exists a canonical number system if and only if $\mathcal{O}_K$ has a power integral basis.   
\end{theorem}

Together with the above Theorem \ref{cor:5.6} of Gy\H{o}ry (1976) this implies that it is effectively decidable whether there exists a CNS in $\mathcal{O}_K$. Theorem \ref{cor:5.6} provides even a general algorithm to determine all power integral bases in $\mathcal{O}_K$. Using this, B. Kov\'acs and Peth\H o (1991) proved as follows.

\begin{theorem}[B. Kov\'acs and Peth\H o, 1991]\label{thm:6.2}
Up to $\mathbb{Z}$-equivalence, there are only finitely many CNS's in $\mathcal{O}_K$, and all of them can be effectively determined.    
\end{theorem}

In fact, using Theorem \ref{thm:5.5},    they extended their result to any order $\mathcal{O}$ of $K$ as well. In an order $\mathcal{O}$, a canonical number system $\{\alpha,\mathcal{N}(\alpha)\}$ is defined in a similar way as in $\mathcal{O}_K$.

We note that Brunotte (2001) considerably improved the procedure of B. Kov\' acs and Peth\H o (1991) and gave an efficient algorithm for finding all such CNS's, provided that one has an efficient algorithm for determining all power integral bases in $\mathcal{O}_K$, resp. in $\mathcal{O}$. As was seen in Section \ref{sec:5}, such an algorithm is known for number fields $K$ of degree at most $4$ if their discriminants are not too large in absolute value.

B. Kov\'acs and Peth\H o (1991) gave also a complete, effective characterization of CNS's in number fields and in their orders.

\begin{theorem}[B. Kov\'acs and Peth\H o, 1991]\label{thm:7.3new}
Let $\mathcal{O}$ be an order in a number field $K$. There exist $\alpha_1,\ldots,\alpha_t\in\mathcal{O}$, $n_1,\ldots,n_t\in\mathbb{Z}$, $N_1,\ldots,N_t$ finite subsets of $\mathbb{Z}$, which are all effectively computable, such that $\{\alpha,\mathcal{N}(\alpha)\}$ is a CNS in $\mathcal{O}$, if and only if $\alpha=\alpha_i-h$ for some integers $i,h$ with $1\le i\le t$ and $h\ge n_i$ or $h\in N_i$.
\end{theorem}

Several generalizations and applications have been obtained. Peth\H o and Varga (2017) generalized the result of B. Kov\'acs to CNS's over imaginary quadratic Euclidean domains. Peth\H o and Thuswaldner (2018) study CNS's in relative extensions. Most of the results of B. Kov\'acs and Peth\H o (1991) are generalized to this situation. Further generalizations are in Evertse, Gy\H ory, Peth\H o and Thuswaldner (2019) over general orders.

Peth\H o (1991) introduced the notion of CNS polynomials. The monic polynomial $P(X)\in\mathbb{Z}[X]$ is called \textit{CNS polynomial} if $|P(0)|\ge 2$ and for every $0\ne Q(X)\in\mathbb{Z}[X]$ there exist unique integers $\ell\ge 0$, $q_1,\ldots,q_\ell\in\{0,1,\ldots,|P(0)|-1\}$ such that
$$Q(X)\equiv \sum_{j=0}^\ell q_jx^j\pmod{P(X)}.$$
He proved that \textit{if $P(X)$ is irreducible and monic and $\alpha$ is one of the zeros of $P(X)$, then $P(X)$ is a CNS polynomial if and only if $\{\alpha,0,1,\ldots,|P(0)|-1\}$ is a \textnormal{CNS} in $\mathbb{Z}[\alpha]$}.

A. Kov\'acs (2001) computed all CNS polynomials with $P(0)=2$ up to degree $8$. This computation was extended up to degree $14$ in Burcsi and A. Kov\'acs (2008).

Akiyama, Borb\'ely, Brunotte, Peth\H o and Thuswaldner (2005) defined the \textit{shift} radix system (SRS). It is a discrete dynamical system, which is a common generalization of CNS polynomials and some kind of $\beta$ representations of real numbers. Many properties of SRS were also described.

For surveys, we refer to Brunotte (2001), Peth\H o (2004), Brunotte, Huszti and Peth\H o (2006), Komornik (2011), Evertse, Gy\H ory, Peth\H o and Thus\-waldner (2019) and the references given there.

\section{Further consequences and applications of the reduction theory
}\label{sec:8}

\noindent
The main results from the effective reduction theory for polynomials discussed before, i.e., Theorems \ref{thm:3.2} and \ref{thm:3.4}, as well as their various versions led to many applications. Some of them were treated in Sections \ref{sec:3} to \ref{sec:6}. Below we briefly present some 
others in their simplest form. For further applications, we refer to the survey paper Gy\H ory (2006), the books Gy\H ory (1980b), Smart (1998), Evertse and Gy\H ory (2017) and the references given there.

\subsection{Applications to classical Diophantine equations.}\label{sec:8.1new}
\vphantom{a}

\noindent
Theorem \ref{thm:3.2} can be applied to superelliptic equations and the Schinzel--Tijdeman equation.
\begin{itemize}

    \item  Let $f\in\mathbb{Z}[X]$ be a monic polynomial of degree $n\ge 3$ with discriminant $D(f)\ne 0$, and $m\ge 2$ an integer. Consider the solutions $x,y\in\mathbb{Z}$ of the equation
    \begin{align}\label{eq:9.5}
        f(x)=y^m.
    \end{align}
                Applying various variants of Theorem \ref{thm:3.2} to the polynomial $f$ and then using Baker's method for the reduced equation, Trelina (1985) and, for $n=3, m=2$, Pint\'er (1995) gave effective upper bounds for $|y|$ that depend on $m,n$ and $|D(f)|$, but \textit{not on the height} of $f$. We recall that the height of $f$ can be arbitrarily large with respect to $|D(f)|$. Furthermore, Gy\H ory and Pint\'er (2008) showed that for each solution $x,y$ of \eqref{eq:9.5} with $\gcd(y,D(f))=1$, $|y|^m$ can be effectively bounded in terms of the radical of $D(f)$, i.e. the product of the distinct prime factors of $D(f)$. It should be noted that $|D(f)|$ can be arbitrarily large with respect to its radical.  Brindza, Evertse and Gy\H ory (1991), Haristoy (2003) and Gy\H ory and Pint\'er (2008) gave upper bounds even for $m$ that depend only on $n$ and $|D(f)|$.
   \item Consider now an application of Theorem \ref{thm:3.2} to \textit{equations of discriminant type}
   \begin{align}\label{eq:8.2new}
   	D(x_1,\ldots,x_n)=D\text{ in } x_1,\ldots,x_n\in\mathbb{Z},
   \end{align}
   where $D(x_1,\ldots,x_n):=D(f(X))$ is the discriminant of the polynomial $f(X)=X^n+x_1X^{n-1}+\cdots + x_n$ in $X$, and $D\ne 0$ is a given rational integer. If $(x_1,\ldots,x_n)$ is a solution of \eqref{eq:8.2new} then so is
   \begin{align*}
   	(x_1^\ast,\ldots,x_n^\ast)=\left( \frac{f^{(n-1)}(a)}{(n-1)!},\ldots,f(a) \right)\text{ for any } a\in\mathbb{Z},
   \end{align*}
   where $X^n+x_1^\ast X^{n-1}+\cdots + x_n^\ast=:f^\ast (X)=f(X+a)$. Such a set of solutions of \eqref{eq:8.2new} is called a \textit{family of solutions}. Using a quantitative version of his Theorem \ref{thm:3.2}, Gy\H ory (1976) proved that \eqref{eq:8.2new} has only finitely many families of solutions and a representative of every family can be effectively determined. Theorem \ref{thm:3.8} gives a considerable improvement of this result of Gy\H ory (1976).
\end{itemize}

The binary form variant of Theorem \ref{thm:3.4} can be applied to Thue equations, Thue inequalities and Thue--Mahler equations.
\begin{itemize}
    \item Let $F\in\mathbb{Z}[X,Y]$ be an irreducible binary form of degree $n\ge 3$ and discriminant $D$, let $p_1,\ldots,p_s\,\,(s\ge 0)$ be distinct primes not exceeding $P$, and let $m$ be a positive integer coprime with $p_1,\ldots,p_s$. There are several upper bounds for the \textit{number} of solutions $x,y$ of the Thue equation
    \begin{align}\label{eq:9.6}
        F(x,y)=m,
    \end{align}
    the Thue inequality
    \begin{align}\label{eq:9.7}
        0<|F(x,y)|\le m
    \end{align}
    and the Thue--Mahler equation
    \begin{align}\label{eq:9.8}
        F(x,y)=mp_1^{z_1}\cdots p_s^{z_s}\text{, with } (x,y)=1,
    \end{align}
    where $z_1,\ldots,z_s$ are also unknown non-negative integers.

Using a quantitative binary form version of Theorem \ref{thm:3.4}, e.g. the general effective Theorem 1 of Evertse and Gy\H ory (1991a) on binary forms of given degree and given discriminant over $\mathbb{Z}$, the previously obtained upper bounds for the \emph{number} of solutions of these equations were substantially improved under the assumptions that $n,D,m,s$ and $P$ satisfy some additional conditions. Such improved upper bounds were derived in Stewart (1991) for \eqref{eq:9.8} with $\gcd (x,y)=1$ when $m>C_1$, in Brindza (1996) for \eqref{eq:9.6} with $\gcd (x,y)=1$ when $m>C_2$, and in Thunder (1995) for \eqref{eq:9.7} when $m>C_3$, where $C_1,C_2,C_3$ are effectively computable numbers such that $C_1$ depends on $n,|D|,P,s$ and $C_2,C_3$ on $n$ and $|D|$. Further, Evertse and Gy\H ory (1991b) showed that if $|D|>C_4$, then the number of coprime solutions of \eqref{eq:9.7} is at most $6n$ if $n>400$, and by Gy\H ory (2001) it is at most $28n+6$ if $|D|>C_5$ and $3\le n\le 400$. For $m=1$ and $|D|>C_6$, this was later improved by Akhtari (2012) to $11n-2$. Here $C_4,C_5,C_6$ are effectively computable numbers such that $C_4,C_5$ depend on $m$ and $n$, and $C_6$ on $n$. Together with the above mentioned quantitative version of Theorem \ref{thm:3.4}, these imply that for given $n\ge 3$ and $m\ge 1$, there are only finitely many $GL_2(\mathbb{Z})$-equivalence classes of irreducible binary forms $F\in\mathbb{Z}[X,Y]$ of degree $n$ for which the number of coprime solutions of \eqref{eq:9.7} exceeds $28n+6$ or $11n-2$ if $m=1$.

\item The quantitative version of Theorem \ref{thm:3.4}, proved in Evertse and Gy\H ory (1991a) was also applied in Evertse (1993) to bound the number of solutions of some resultant inequalities, and in Ribenboim (2006) to binary forms with given discriminant, having additional conditions on the coefficients. 

We remark that using the improved and completely explicit version Theorem \ref{thm:3.9} of Evertse and Gy\H ory (2017), the above quoted applications can be made more precise.
\end{itemize}

\subsection{Some other applications of Theorems \ref{thm:3.2} and \ref{thm:3.4}}\label{sec:8.2new}

\begin{itemize}
    \item In Evertse and Gy\H ory (2017), as a consequence of Theorem \ref{thm:3.9}, we derived for any separable polynomial $f\in\mathbb{Z}[X]$ of degree $n\ge 4$ an improvement of the previous bounds for the minimal root distance of $f$.   
    \item Some applications of Theorem \ref{thm:3.2} were given to the reducibility of a general class of polynomials of the form $g(f(X))$ where $f,g$ are monic polynomials, $g(X)$ is irreducible with CM splitting field. For given prime $p$ and $g\in \Zz [X]$, there are  up to $\mathbb{Z}$-equivalence only finitely many $f\in \Zz [X]$ of degree $p$ with distinct real zeros for which $g(f(X))$ is reducible; see Gy\H ory (1976, 1982).
    \item In Evertse and Gy\H ory (1991a), a quantitative binary form variant of Theorem \ref{thm:3.4} was utilized to give effective upper bounds for the minimal non-zero absolute value of binary forms at integral points.
    \item For an application of an earlier version of Theorem \ref{thm:5.2} (ii) to integral valued polynomials, see Peruginelli (2014).
    \item For an application of Theorem \ref{thm:5.2} to so-called binomially equivalent numbers, see Yingst (2006).
\end{itemize}

As will be seen in the next section, the various generalizations presented there of Theorems \ref{thm:3.2} and \ref{thm:3.4} have also several applications.

\section{Generalizations and their consequences, applications}\label{sec:9}

In Sections \ref{sec:3} to \ref{sec:8} we presented the most significant results and consequences/applications of the effective reduction theory of integral polynomials over $\mathbb{Z}$. In the last decades this effective theory has been generalized by the authors among others for the number field case, more precisely for the case of integral, resp. $S$-integral polynomials over number fields. In the monic case, they have obtained even more general effective results for polynomials over finitely generated domains of characteristic $0$ which may contain transcendental elements, too. These provided many important consequences and applications, and yielded a further advancement in the theory.

In this section we formulate some typical general effective theorems on integral polynomials over number fields and finitely generated domains, including various generalizations of Theorems \ref{thm:3.2}, \ref{thm:3.4} and their consequences.

For simplicity, we present them in qualitative forms. For explicit versions and further results and applications, we refer to our original works or our books Evertse and Gy\H ory (2017, 2022). The proofs depend explicitly or implicitly on an effective finiteness theorem of Gy\H ory (1979) or its improvements by Bugeaud and Gy\H ory (1996a), Gy\H ory and Yu (2006), Evertse and Gy\H ory (2015) or Gy\H{o}ry (2019), see Theorem \ref{thm:4.ww} above on $S$-unit equations, resp. of Evertse and Gy\H ory (2013) on unit equations over finitely generated domains.

For convenience, the monic and non-monic cases are treated separately in the Subsections \ref{sec:9.1} and \ref{sec:9.2} below.



\subsection{Generalizations: the monic case}\label{sec:9.1}
\subsubsection{Results over number fields}\label{sec:9.1.1}
\vphantom{a}

\noindent
Let $L$ be a number field with ring of integers $\mathcal{O}_L$, and $S$ a finite set of places on $L$ containing all infinite places $S_\infty$. The ring of $S$-integers of $L$, denoted by $\mathcal{O}_S$, consists of those elements of $L$ which are integral at every finite placee outside $S$.
A fractional ideal of $\OO_S$ is a subset $\fa$ of $L$ such that there is non-zero $\delta\in L$
such that $\delta\fa$ is an ideal of $\OO_S$. 
Given a subset $\VV\not=\{ 0\}$ of $L$ such that
$\delta\VV\subset \OO_S$ for some non-zero $\delta\in \OO_S$, we denote by $(\VV )_S$ the fractional ideal of $\OO_S$ generated by $\VV$. 
Lastly, we denote by $\mathcal{O}^\ast_S$ the unit group of $\mathcal{O}_S$.

Two \textit{monic} polynomials $f,g\in\mathcal{O}_S[X]$ of degree $n$ are called $\mathcal{O}_S$-\textit{equivalent} if
$$g(X)=\varepsilon^nf(\varepsilon^{-1}X+a)\text{ for some }\varepsilon\in\mathcal{O}^\ast_S\text{ and } a\in\mathcal{O}_S,$$
and \textit{strongly $\mathcal{O}_S$-equivalent} if
$$g(X)=f(X+a)\text{ for some } a\in\mathcal{O}_S.$$
In this case $D(g)=\varepsilon^{n(n-1)}D(f)$, resp. $D(g)=D(f)$.

For a polynomial $g\in\overline{\mathbb{Q}}[X]$, we denote by $H(g)$ the absolute height of the vector whose coordinates are the coefficients of $g$.

\begin{theorem}[Gy\H ory, 1978b, 1984]\label{thm:9.1new}
    Let $\delta\in\mathcal{O}_S\setminus\{0\}$, and let $f\in\mathcal{O}_S[X]$ be a monic polynomial of degree $n\ge 2$ with discriminant $D(f)\in\delta\mathcal{O}^\ast_S$. Then $f$ is $\mathcal{O}_S$-equivalent to a monic polynomial $g\in\mathcal{O}_S[X]$ for which
    $$H(g)<C_1(L,S,(\delta )_S,n)$$
    where $C_1$ is an effectively computable number depending only on $L,S,\delta$ and $n$.
\end{theorem}

For $L=\mathbb{Q}$, $S=S_\infty$, this is just Theorem \ref{thm:3.2}, (ii), where the bound given for $H(g)$ is in fact independent of $n$. We note that this is not the case in general, see Evertse and Gy\H ory (2017), p. 155.

For the best known, completely explicit bound $C_1$ see also Evertse and Gy\H ory (2017), Theorem 8.2.3.

Theorem \ref{thm:9.1new} implies the following effective finiteness results.

\begin{corollary}\label{cor:9.2}
    For given integer $n\ge 2$ and $\delta\in\mathcal{O}_S\setminus\{0\}$, there are only finitely many $\mathcal{O}_S$-equivalence classes of monic polynomials $f$ in $\mathcal{O}_S[X]$ of degree $n$ and with $D(f)\in\delta\mathcal{O}_S^\ast$. Further, there exists an algorithm that for any $n\ge 2$ and any effectively given $L,S$ and $\delta$ computes a full set of representatives of these classes.
\end{corollary}

Theorem \ref{thm:9.1new} gives also in an effective form that there are only finitely many strong $\mathcal{O}_S$-equivalence classes of monic polynomials $f\in\mathcal{O}_S[X]$ of given degree $n\ge 2$ and with given discriminant $D(f)=\delta\ne 0$. For a quantitative and explicit version, see Corollary 8.2.6 in Evertse and Gy\H ory (2017).

We recall that for definitions of \textit{effectively given} concepts, structures, etc. we referred in Subsection \ref{sec:0.1} to the corresponding sections of our books Evertse and Gy\H ory (2015, 2017, 2022).

We now present another version of Theorem \ref{thm:9.1new} which is more convenient to apply.

With the above notation, let $\mathcal{L}=\mathcal{O}_S^\ast\cap\mathcal{O}_L$. If $\mathfrak{p}_1,\ldots,\mathfrak{p}_t$ denote the prime ideals of $\mathcal{O}_L$ corresponding to the finite places of $S$, then $\mathcal{L}$ is just the multiplicative semigroup of non-zero elements of $\mathcal{O}_L$ which are not divisible by any prime ideal different from $\mathfrak{p}_1,\ldots,\mathfrak{p}_t$. The set $\mathcal{L}$ contains obviously the unit group $\mathcal{O}_L^\ast$ of $\mathcal{O}_L$, and, for $t=0$, $\mathcal{L}=\mathcal{O}_L^\ast$.

We say that the monic polynomials $f,g\in\mathcal{O}_L[X]$ are \textit{strongly\, $\mathcal{O}_L$-equivalent} if
$$g(X)=f(X+a)\text{ for some } a\in\mathcal{O}_L.$$

The next theorem was proved in Gy\H ory (1978b) in a quantitative form.
\begin{theorem}[Gy\H ory, 1978b]\label{thm:9.3new}
    Let $L,\mathcal{L}$ be as above and let $\delta$ be a non-zero element of $\mathcal{O}_L$. If $f\in\mathcal{O}_L[X]$ is a monic polynomial of degree $n\ge 2$ with discriminant $D(f)\in\delta\mathcal{L}$, then it is strongly $\mathcal{O}_L$-equivalent to a polynomial of the form $\eta^ng(\eta^{-1}X)$, where $\eta\in\mathcal{L}$, $g\in\mathcal{O}_L[X]$ and
    $$H(g)\le C_2(L,\mathcal{L},(\delta )_S,n),$$
    where $C_2$ is an effectively computable number depending only on $L, \mathcal{L},(\delta )_S$ and $n$.
\end{theorem}

If $L=\mathbb{Q}$ and $t=0$, then $\mathcal{L}=\{\pm1\}$, and Theorem \ref{thm:9.3new} gives again Theorem \ref{thm:3.2} (ii). For a more general version of Theorem \ref{thm:9.3new} with not necessarily non-zero $\delta$, see also Theorem 2 in Gy\H ory (1981).

We present now some applications of Theorem \ref{thm:9.3new} to algebraic integers whose \textit{discriminants} resp. \textit{indices} over $L$ belong to $\delta\mathcal{L}$.

For an algebraic integer $\alpha$ of degree $n\ge 2$ over $L$, $f_{\alpha ,L}$ will denote the monic minimal polynomial of $\alpha$ over $L$, i.e., the monic polynomial in $\OO_L[X]$ of minimal degree of which $\alpha$ is a zero. We define the discriminant of $\alpha$ relative to $L$ by 
\[
D_L(\alpha ):=D(f_{\alpha ,L})=\prod_{1\leq i<j\leq n}(\alpha^{(i)}-\alpha^{(j)})^2,
\]
where $\alpha^{(1)}\kdots\alpha^{(n)}$ are the conjugates of $\alpha$ over $L$. 
The algebraic integers $\alpha$ and $\beta$ are said to be \textit{strongly $\mathcal{O}_L$-equivalent} over $L$ when $\alpha-\beta\in\mathcal{O}_L$. In this case their minimal polynomials over $L$ are also strongly $\mathcal{O}_L$-equivalent.

We denote by $H(\beta)$ the absolute height of an algebraic number $\beta$.

Corollary \ref{cor:9.4} is an immediate consequence of Theorem \ref{thm:9.3new}. It was proved in Gy\H ory (1978b) in quantitative form.

\begin{corollary}[Gy\H ory, 1978b]\label{cor:9.4}
    Let $L,\mathcal{L}$ and $\delta$ be as in Theorem \ref{thm:9.3new}, and let $\alpha$ be an algebraic integer with degree $n\ge 2$ and discriminant $D_L(\alpha)\in\delta\mathcal{L}$ over $L$. Then $\alpha$ is strongly $\mathcal{O}_L$-equivalent to an algebraic integer of the form $\eta\beta$, where $\eta\in\mathcal{L}$ and $\beta$ is an algebraic integer satisfying
    $$H(\beta)<C_3(L,\mathcal{L},\delta,n)$$
    with an effectively computable number $C_3$ depending only on $L,\mathcal{L},\delta$ and $n$.
\end{corollary}

This is a considerable effective generalization of Theorem \ref{thm:5.2} in two different directions, for the number field case and for the $\mathfrak{p}$-adic case. We note that in the special case $L=\mathbb{Q}$, Corollary \ref{cor:9.4} was proved independently by Trelina (1977a).

A simple consequence of Corollary \ref{cor:9.4} is that up to the obvious multiplications by elements of $\mathcal{L}$ and translations by integers of $L$, there are only finitely many algebraic integers $\alpha$ with given degree $n$ and discriminant $D_L(\alpha)\in\delta\mathcal{L}$ over $L$ and they can be effectively determined. As is remarked in Gy\H ory (1978b), the first, finiteness part can be deduced, in an \textit{ineffective} form, from the \textit{ineffective}  Theorem \ref{thm:3.1} of Birch and Merriman (1972) over number fields and from the finiteness of the number of solutions of the generalized Thue--Mahler equation; cf. Parry (1950). 

As is pointed out in Gy\H ory (1978b), p. 177, if in Corollary \ref{cor:9.4} we restrict ourselves to integers $\alpha$ of a fixed algebraic number field $K$ of degree $n\ge 3$ over $L$, then the proof of Corollary \ref{cor:9.4} in Gy\H ory (1978b) gives the following in quantitative form.
\begin{corollary}[Gy\H ory, 1978b, 1981]\label{cor:9.5}
    Let $L,\mathcal{L},\delta$ and $K$ be as above, and let $\alpha$ be a primitive integral element of $K$ with discriminant $D_{K/L}(\alpha)\in\delta\mathcal{L}$ over $L$. Then $\alpha$ is strongly $\mathcal{O}_L$-equivalent to an algebraic integer of the form $\eta\beta$, where $\eta\in\mathcal{L}$, and $\beta$ is an algebraic integer in $K$ such that
    $$H(\beta)<C_4(L,K,\mathcal{L},\delta,n)$$
    with an effectively computable number $C_4$ which depend only on $L,K,\mathcal{L}$ and $\delta$.  
\end{corollary}

Keeping the above notations, we present some consequences of Corollary \ref{cor:9.5}. Consider an order $\mathcal{O}$ of the field extension $K/L$ (i.e. let $\mathcal{O}$ be a subring of $\mathcal{O}_K$, the ring of integers of $K$, that has the full dimension $n$ as an $\mathcal{O}_L$-module). Denote by
 $\fD_{K/L}(\mathcal{O})$ the \textit{discriminant ideal} of $\mathcal{O}$. Then we have (cf. Fr\"ohlich, 1967)
$$(D_{L}(\alpha))=\fI^2_\mathcal{O}(\alpha)\cdot \fD_{K/L}(\mathcal{O})$$
for any $\alpha\in\mathcal{O}$ such that $L(\alpha)=K$. Here $\fI_\mathcal{O}(\alpha)$ is an integral ideal which is called the \textit{index} of $\alpha$ in $\mathcal{O}$. It is clear that if $\alpha,\beta\in\mathcal{O}$ are strongly $\mathcal{O}_L$-equivalent then $\fI_\mathcal{O}(\alpha)=\fI_\mathcal{O}(\beta)$.
\begin{corollary}[Gy\H ory, 1981]\label{cor:9.6}
    If $\alpha\in\mathcal{O}$ has index $\fI_\mathcal{O}(\alpha)$ not divisible by any prime ideal different from $\mathfrak{p}_1,\ldots,\mathfrak{p}_t$, then $\alpha$ is strongly $\mathcal{O}_L$-equivalent to an algebraic integer of the form $\eta\beta$, where $\eta\in\mathcal{L},\beta\in\mathcal{O}$, and
    $$H(\beta)<C_5(L,K,\mathcal{O},\mathfrak{p}_1,\ldots,\mathfrak{p}_t),$$
    where $C_5$ is an effectively computable number depending only on\linebreak $L,K,\mathcal{O},\mathfrak{p}_1,\ldots,\mathfrak{p}_t$.
\end{corollary}

In the case $\mathcal{O}=\mathcal{O}_K$, a prime ideal $\mathfrak{p}$ in $L$ is called a \textit{common index divisor} of $K/L$ if $\mathfrak{p}$ divides $\fI_{\mathcal{O}_K}(\alpha)$ for every primitive integral element $\alpha$ of $K/L$. The number of common index divisors is finite and a well-known theorem of Hasse (1980) gives an elegant characterization of these divisors. It is interesting to apply Corollary \ref{cor:9.6} to the case when $\mathfrak{p}_1,\ldots,\mathfrak{p}_t$ are just the common index divisors of $K/L$. There are relative extensions of arbitrary high degree in which there exists no element $\alpha$ with index not divisible by prime ideals different from the common index divisors; cf. Pleasants (1974). Corollary \ref{cor:9.6} provides and effective algorithm for deciding whether such an element $\alpha$ exists and for determining all $\alpha$ having this property.

Corollaries \ref{cor:9.5} and \ref{cor:9.6} allowed Gy\H ory (1981) to get some information about the arithmetical structure of those non-zero algebraic integers resp. non-zero integral ideals in $L$ which are discriminants resp. indices of elements of $\mathcal{O}_K$ over $L$.

We now present generalizations of Theorems \ref{cor:5.6} and \ref{thm:5.5} for the relative case.

Let again $L$ be a number field, $K$ an extension of degree $n\ge 2$ of $L$, and $\mathcal{O}$ an order of $K$ over $L$. Then $\mathcal{O}=\mathcal{O}_L[\alpha]$ for some $\alpha\in\mathcal{O}$ if and only if $\fI_\mathcal{O}(\alpha)=\mathcal{O}_L$. In this case $\{1,\alpha,\ldots,\alpha^{n-1}\}$ is a $\OO_L$-module basis for $\mathcal{O}$. There exists an extensive literature of such power bases of orders of number fields and related topics; we refer the reader to the works of Hensel (1908), Hasse (1980), Narkiewicz (1974), Gy\H ory (1978a, 1978/79), and Evertse and Gy\H ory (2017), and thence to the literature mentioned there.

We say that $\alpha,\beta\in\mathcal{O}$ are $\mathcal{O}_L$\textit{-equivalent} if $\beta=a+\varepsilon\alpha$ for some $a\in\mathcal{O}_L$ and unit $\varepsilon$ in $\mathcal{O}_L$. If $\alpha$ is a generator of $\mathcal{O}$ over $\mathcal{O}_L$, i.e. $\mathcal{O}=\mathcal{O}_L[\alpha]$ then so is every $\beta$ which is $\mathcal{O}_L$-equivalent to $\alpha$.

The following fundamental theorem is a consequence of Corollary \ref{cor:9.6}.

\begin{theorem}[Gy\H ory, 1981]\label{thm:9.7}
    Let $\mathcal{O}$ be an order of $K/L$, and suppose that $\mathcal{O}=\mathcal{O}_L[\alpha]$ for some $\alpha\in\mathcal{O}$. Then there is $\beta\in\OO$ that is $\OO_L$-equivalent to $\alpha$   
and for which
    $$H(\beta)<C_6(L,K,\mathcal{O}),$$
    where $C_6$ is an effectively computable number depending only on $L,K$ and $\mathcal{O}$.
\end{theorem}

For $L=\mathbb{Q}$, this gives Theorems \ref{cor:5.6} and \ref{thm:5.5} above. 
In the case $\mathcal{O}=\mathcal{O}_K$, Theorem \ref{thm:9.7} was proved in Gy\H ory (1978a) with a completely explicit bound corresponding to $C_6$. For the best known explicit bound in Theorem \ref{thm:9.7}, see Corollary 8.4.13 in Evertse and Gy\H ory (2017).

Theorem \ref{thm:9.7} provides a general effective algorithm for deciding whether a relative extension $K/L$ resp. an order $\mathcal{O}$ of $K$ over $L$ is monogenic or not, and for determining all $\alpha\in\mathcal{O}_K$ resp. all $\alpha\in\mathcal{O}$ for which $\mathcal{O}_K=\mathcal{O}_L[\alpha]$ resp. $\mathcal{O}=\mathcal{O}_L[\alpha]$.

We now present a very important consequence of Theorem \ref{thm:9.7}.
Let again $L$ be a number field, $K$ an extension of degree $n\ge 2$, and $\mathcal{O}_K, \mathcal{O}_L$ the rings of integers of $K$ resp. $L$. Pleasants (1974) gave an explicit formula which enables one to compute a positive integer $m(\mathcal{O}_K,\mathcal{O}_L)$ such that if $r(\mathcal{O}_K,\mathcal{O}_L)$ denotes the minimal number of generators of $\mathcal{O}_K$ as $\mathcal{O}_L$-algebra then
    \begin{align*}
        m(\mathcal{O}_K,\mathcal{O}_L)\le r(\mathcal{O}_K,\mathcal{O}_L)\le\max\{m(\mathcal{O}_K,\mathcal{O}_L),2\}.
    \end{align*}
    Pleasants proved that if $L=\mathbb{Q}$, there are number fields $K$ of arbitrarily large degree over $\mathbb{Q}$ such that $m(\mathcal{O}_K,\mathbb{Z})=1$ and $\mathcal{O}_K$ is not monogenic. Consequently, his theorem does not make it possible to decide whether the ring of integers of a number field is monogenic. Together with Pleasants' result, our Theorem \ref{thm:9.7} above gives the following
    
    \begin{customcor}{of Theorem \ref{thm:9.7}}[and of Pleasants (1974)]
        There is an algorithm for determining the least number of elements of $\mathcal{O}_K$ that generate $\mathcal{O}_K$ as an $\mathcal{O}_L$-algebra.
    \end{customcor}

            Chapter 11 of Evertse and Gy\H ory (2017)
        considers more generally $\mathcal{O}_S$-orders of finite \'etale $L$-algebras, and gives a method to determine a system of $\mathcal{O}_S$-algebra generators of minimal cardinality of such an order.
    This was basically work of Kravchenko, Mazur and Petrenko (2012), worked out in more detail in a special case.

We give an overview of generalizations of some of the results from the previous sections.
\begin{itemize}
    \item In Section \ref{sec:4}, several results have reformulations in terms of polynomial Diophantine equations; see equations \eqref{eq:4.4new}  and \eqref{eq:4.4}. In the present section the above extensions of the results from Section \ref{sec:4} have also reformulation in terms of discriminant form equations and index form equations over number fields and in the $\mathfrak{p}$-adic case.
        \item Corollaries  \ref{cor:5.4new}, \ref{cor:4.4} were first extended to the case when $D$ resp. $I$ is replaced by $p_1^{u_1}\cdots p_s^{u_s}$, where $p_1,\ldots,p_s$ are fixed primes and $u_1,\ldots,u_s$ are unknown non-negative integers; see Gy\H ory (1978b, 1981), Trelina (1977a, 1977b), Gy\H ory and Papp (1977). These results yielded e.g. explicit lower bounds for the greatest prime factor of discriminant and index of an integer of a number field. For generalizations for the number field case, see Gy\H ory (1980a, 1981).
    \item Corollary \ref{cor:5.4new} on discriminant form equations was generalized for more general decomposable form equations of the form
    \begin{align}\label{eq:9.1}
        F(x_1,\ldots,x_m)=F\text{ in } x_1,\ldots,x_m\in\mathbb{Z},
    \end{align}
    where $F\in\mathbb{Z}\setminus\{0\}$ and $F(X_1,\ldots,X_m)$ is a decomposable form with coefficients in $\mathbb{Z}$ which factorizes into linear factors over $\overline{\mathbb{Q}}$ such that these factors form a so-called triangularly connected system (i.e. \eqref{eq:9.1} can be reduced to a connected system of three terms unit equations); see Gy\H ory and Papp (1978) and, more generally, Gy\H ory (1998).

    For discriminant form equations and more general decomposable form equations, see also Evertse and Gy\H ory (2017), Chapters 6, 8 and 10, and Evertse and Gy\H ory (2022), Chapters 2 and 4.    
    \item Corollary \ref{cor:5.4new} was generalized for the `inhomogeneous' case by Ga\'al (1986).
    \item Analogous results were established over function fields by Gy\H ory (1984, 2000); Ga\'al (1988), Mason (1988), Shlapentokh (1996).
\end{itemize}

\subsubsection{Results over finitely generated domains}\label{sec:9.1.2}
\vphantom{a}

\noindent
We now present two general finiteness theorems where the ground ring is an integrally closed integral domain $A$ of characteristic $0$ that is finitely generated over $\mathbb{Z}$ as a $\Zz$-algebra, i.e., $A=\Zz [z_1\kdots z_r]$,
where we allow some of the $z_i$ to be transcendental.

We say that the monic polynomials $f,g\in A[X]$ are \textit{strongly $A$-equivalent} if $g(X)=f(X+a)$ with some $a\in A$. Then $f$ and $g$ have the same discriminant.

\begin{theorem}[Gy\H ory, 1982]\label{thm:9.2}
Let $G$ be a finite extension of the quotient field of $A$.
    Up to strong $A$-equivalence, there are only finitely many monic $f(X)$ in $A[X]$ with given non-zero discriminant $\delta$ having all their zeros in $G$.
\end{theorem}

This was made effective by Gy\H ory (1984) in a special case, and in full generality by Evertse and Gy\H ory (2017), provided that $A, G$ and $\delta$ are given effectively in the sense defined in Evertse and Gy\H ory (2017, 2022).

\begin{theorem}[Evertse and Gy\H ory, 2017, 2022]\label{thm:9.3}
    Let $A,G,\delta$ be as above. Up to strong $A$-equivalence, there are only finitely many monic $f(X)$ in $A[X]$ with $D(f)=\delta$, and if $A,G,\delta$ are effectively given, all these $f$ can be effectively determined.
\end{theorem}

\begin{prob}\label{prob:2}
    Are Theorems \ref{thm:9.2} and \ref{thm:9.3} true without fixing the splitting field $G$?
\end{prob}

Several results of the theory have been extended to the case of \'etale algebras in Evertse and Gy\H ory (2017, 2022).
\begin{itemize}
    \item Let $K$ be a number field with ring of integers $\mathcal{O}_K$, and $D\ne 0$ an integer. As was seen above, up to strong $\Zz$-equivalence,  the equation
    \begin{align}\label{eq:9.2}
        D(\alpha)=D\text{ in } \alpha\in\mathcal{O}_K
    \end{align}
    has only finitely many solutions, and all of them can be effectively determined.

    Let $A=\mathbb{Z}[z_1,\ldots,z_r]$ be an integral domain of characteristic $0$ with algebraic or transcendental generators
    $z_1\kdots z_r$, $L$ its quotient field, and
    $\Omega$ a \textit{finite \'etale $L$-algebra} (i.e. a direct product of finite extensions $K_1,\ldots,K_t$ of $L$).
    Denote by $A_{\Omega}$ the integral closure of $A$ in $\Omega$.
    The \textit{discriminant} of $\alpha\in A_{\Omega}$ over $L$ with
    $\Omega =L[\alpha ]$ is given by $D_{L}(\alpha ):=D(f_{\alpha ,L})$,
    where $f_{\alpha ,L}$ is the monic minimal polynomial of $\alpha$
    over $L$.

    Let $\mathcal{O}$ be an $A$\textit{-order} of $\Omega$, i.e. an $A$-subalgebra of $A_{\Omega}$ which spans $\Omega$ as an $L$-vector space.
    We say that $\alpha,\beta\in \mathcal{O}$ are \textit{strongly    $A$-equivalent} if $\beta-\alpha\in A$.
    One verifies that if $\alpha ,\beta\in\mathcal{O}$
    are strongly $A$-equivalent then $f_{\alpha ,L}$, $f_{\beta ,L}$
    are also strongly $A$-equivalent, and thus, $D_L(\beta )=D_L(\alpha )$.
    
Let $\delta$ be a non-zero element of $L$.     
    Consider the following generalization of equation \eqref{eq:9.2}:
    \begin{align}\label{eq:9.3}
        D_{L}(\alpha)=\delta\text{ in }\alpha\in\mathcal{O}.
        \end{align}
        For an integral domain $B$, denote by $B^+$ the additive group of $B$.

\begin{theorem}[Evertse and Gy\H ory, 2022]\label{thm:9.4}
    If
    \begin{align}\label{eq:9.4}
        (\mathcal{O}\cap L)^+/A^+ \text{\textit{ is finite,}}
    \end{align}
    then the set of $\alpha\in\mathcal{O}$ with \eqref{eq:9.3} is a union of finitely many strong $A$-equivalence classes. Moreover, if $A, \Omega,\mathcal{O}$ and $\delta$ are given effectively in a well-defined way, one can determine a set consisting of precisely one element from each of these classes.
\end{theorem}

The condition \eqref{eq:9.4} is necessary and decidable.

For $A=\mathbb{Z}$, $L=\mathbb{Q}$, $\Omega=$ number field $K$, $\mathcal{O}=\mathcal{O}_K$, Theorem \ref{thm:9.4} gives the above theorem concerning equation \eqref{eq:9.2}.
\end{itemize}

\subsection{Generalizations: the non-monic case}\label{sec:9.2}
\vphantom{a}

\noindent
As was seen above, Theorem \ref{thm:3.2} (ii) and its consequences in Sections \ref{sec:3} and \ref{sec:4} were later extended to the number field case and $\mathfrak{p}$-adic case. Theorem \ref{thm:3.4} was already generalized for the same generality by the authors in the first paper on the subject Evertse and Gy\H ory (1991a).

We present now this general theorem from the non-monic case which corresponds to Theorem \ref{thm:9.1new} above.

Keeping the above notations, let again $L$ be a number field, and $S$ a finite set of places on $L$ containing all infinite places. We denote by $\mathcal{O}_S$ the ring of $S$-integers and by $\mathcal{O}^\ast_S$ the group of $S$-units. Two polynomials $f,g\in\mathcal{O}_S[X]$ of degree $n$ are said to be \emph{$GL_2(\mathcal{O}_S)$-equivalent} if
$$g(X)=\varepsilon(cX+d)^n f\left(\frac{aX+b}{cX+d}\right)\text{ with some }\begin{pmatrix}
    a&b\\c&d\end{pmatrix}\in GL_2(\mathcal{O}_S)\text{ and }\varepsilon\in\mathcal{O}^\ast_S.$$

As above, for a polynomial $g(X)\in\overline{\mathbb{Q}}[X]$ we denote by $H(g)$ the absolute height of the vector whose coordinates are the coefficients of $g$.

\begin{theorem}[Evertse and Gy\H ory, 1991a]\label{thm:9.11}
    Let $\delta\in\mathcal{O}_S\setminus\{0\}$, and let $f\in\mathcal{O}_S[X]$ be a polynomial of degree $n\ge 2$ and of discriminant $D(f)\in\delta\mathcal{O}^\ast_S$. Then $f$ is $GL_2(\mathcal{O}_S)$-equivalent to a polynomial $g\in\mathcal{O}_S[X]$ such that
    $$H(g)<C_7(L,S,(\delta )_S,n),$$
    where $C_7$ is an effectively computable number, given explicitly in terms of $L,S,(\delta )_S$ and $n$.
\end{theorem}

For $L=\mathbb{Q}, \mathcal{O}_S=\mathbb{Z}$, when $\mathcal{O}^\ast_S=\{\pm 1\}$, Theorem \ref{thm:9.11} gives Theorem \ref{thm:3.4}. For the best known, completely explicit bound $C_7$, see Theorem 14.2.2 in Evertse and Gy\H ory (2017).   

The binary form variant of Theorem \ref{thm:3.4} was later generalized for decomposable forms in more than two variables in Evertse and Gy\H ory (1992) and Gy\H ory (1994).

Let $K$ be an extension of $L$ of degree $n\geq 3$. 
Let $\alpha$ be a primitive element of $K/L$, i.e., $K= L(\alpha )$.
We would have liked to define the discriminant of $\alpha$ over $\OO_S$ to be the discriminant of $f$, where $f$ is a minimal polynomial of $\alpha$ in $\OO_S[X]$ whose coefficients generate the unit ideal. But in case that $\OO_S$ is not a principal ideal domain, such a minimal polynomial need not exist. Instead, we give a more subtle definition.
Denote by $\PP_S(\alpha )$ the set of polynomials $f\in\OO_S[X]$ such that $f$ is irreducible in $L[X]$
and $f(\alpha )=0$, and  
define the discriminant ideal of $\alpha$ with respect to $\OO_S$ by 
\[
\fd_S(\alpha ):=( D(f):\, f\in \PP_S(\alpha ) )_S.
\]
Given $f(X)=a_0X^n+\cdots +a_n\in\PP_S(\alpha )$, let $\fc_S(f):=(a_0\kdots a_n)_S$ denote its content. Then
\begin{equation}\label{eq;9.x} 
\fd_S(\alpha )= D(f)\cdot \fc_S^{2-2n}.
\end{equation}

Two elements $\alpha ,\beta$ of $K$ are called $GL_2(\OO_S)$-equivalent if $\beta =\medfrac{a\alpha +b}{c\alpha +d}$ for some $\big(\begin{smallmatrix}a&b\\c&d\end{smallmatrix}\big)\in GL_2(\OO_S)$;
such elements 
satisfy $\fd_S(\alpha )=\fd_S(\beta )$.

Theorem \ref{thm:9.11} has the following consequence.

\begin{theorem}\label{thm:9.12}
Let $\alpha$ with $K= L(\alpha)$. Then $\alpha$ is $GL_2(\OO_S)$-equivalent to an element $\beta\in K$ with
\[
H(\beta )\leq C_8(L,S,\fd_S(\alpha ),n),
\]
where $C_8$ is an effectively computable number, given explicitly in terms of $L,S,\fd_S(\alpha )$ and $n$.
\end{theorem}

\begin{proof}[Idea of proof] 
Choose a finite set of ideals of $\OO_S$ that form a full system of representatives for the ideal classes of $\OO_S$. This depends only on $L$ and $S$. There is $f\in\PP_S(\alpha )$ such that $\fc_S(f)=\fa$, where $\fa$ belongs to this finite set of ideals. By Theorem \ref{thm:9.11}, there is $g\in\OO_S[X]$, $GL_2(\OO_S)$-equivalent to $f$,
such that 
\[
H(g) <C_7(L,S,(D(f) )_S,n)=C_7(L,S,\fa^{2n-2}\fd_S(\alpha ),n).
\]
Now $g$ has a zero $\beta$ that is $GL_2(\OO_S)$-equivalent to $\alpha$, and for this $\beta$ we have
$H(\beta )<C_8(L,S,\fd_S(\alpha ),n)$.
\end{proof}

\section{Multiply monogenic and rationally monogenic orders}\label{sec:10}

In this section we consider `Diophantine equations'
\begin{align}\label{eq:10.1}
&\Zz [\alpha ]=\OO\ \ \text{in algebraic integers $\alpha$,}
\\
\label{eq:10.2}
&\Zz_{\alpha }=\OO\ \ \text{in algebraic numbers $\alpha$,}
\end{align}
where $\OO$ is a given order of a number field.
As observed before, from the effective reduction theory for polynomials one can deduce effective finiteness results for the collection of $\Zz$-equivalence classes of algebraic integers $\alpha$ with
\eqref{eq:10.1}, 
respectively the collection of $\GL_2(\Zz )$-equivalence classes of algebraic numbers $\alpha$ with \eqref{eq:10.2}.
Although this does not strictly belong to the effective reduction theory for polynomials,
in this section, we give an overview of results with upper bounds for the \emph{number} of these classes, i.e., for the multiplicity of (rational) monogenicity for the order $\OO$ under consideration. An important feature of these bounds is their uniformity, i.e., they depend at most on the rank of $\OO$. We have included outlines of the proofs of the main results. The main tools are upper bounds for the number of solutions of  equations $ax+by=1$ in algebraic units $x,y$. 

\subsection{Monogenic orders}\label{sec:appA}
\vphantom{a}

\noindent
In this subsection, we consider \eqref{eq:10.1}.
Let $K$ be a number field with ring of integers $\mathcal{O}_K$, and $\mathcal{O}$ an arbitrary order of $K$, i.e., a subring of $\mathcal{O}_K$ with quotient field $K$. It follows from Theorem \ref{thm:5.5} above (in an effective form) that up to $\mathbb{Z}$-equivalence, there are only finitely many $\alpha\in \mathcal{O}$ with $\mathcal{O}=\Zz [\alpha ]$. The order $\mathcal{O}$ is said to be \textit{$k$-times monogenic/precisely $k$ times monogenic/at most $k$ times monogenic} if there are at least/precisely/at most $k$ pairwise $\mathbb{Z}$-inequivalent such generators $\alpha$ of $\mathcal{O}$ over $\mathbb{Z}$.

It is easy to see that every order of a quadratic number field is precisely one time monogenic.

For fixed $n\ge 3$, we denote by $M(n)$ the smallest integer $k$ such that for every number field $K$ of degree $n$
and every order $\mathcal{O}$ of $K$, the order $\mathcal{O}$ is at most $k$ times monogenic.
We start with recalling an old result of ours.

\begin{theorem}[Evertse and Gy\H{o}ry, 1985]\label{thm:10.1}
Let $K$ be a number field of degree $n\geq 3$, and suppose that its normal closure has degree $g$.
Then every order of $K$ is at most $(3\times 7^{2g})^{n-2}$ times monogenic.

In particular, $M(n)$ is finite, and $M(n)\leq (3\times 7^{2n!})^{n-2}$.
\end{theorem}

This was deduced from an upper bound for the number of solutions of $S$-unit equations, obtained shortly before by the first author,
see Evertse (1984a).

There are now much better upper bounds for $M(n)$.
The problem of estimating $M(3)$ can be reduced via index form equations to estimating the number of integer solutions of a Thue equation $|F(x,y)|=1$
with $F$ an integral cubic binary form. Bennett (2001) proved that such an equation has up to sign at most $10$ solutions. This gives the following.

\begin{theorem}[Bennett, 2001]\label{thm:A1}
We have $M(3)\le 10$.
\end{theorem}

For $n\ge 4$, the first author improved the bound of Theorem \ref{thm:10.1} as follows.

\begin{theorem}[Evertse, 2011]\label{thm:A2}
    For $n\ge 4$, $M(n)\le 2^{4(n+5)(n-2)}$ holds.
\end{theorem}

The main tool in the proof is an important improvement and generalization  of the first author's result from 1984, due to Beukers and Schlickewei (1996), see Theorem \ref{thm:10.2}
in Section \ref{sec:10.2}.

In the case of quartic number fields, Bhargava (2022) substantially improved Evertse's bound by proving the following theorem.

\begin{theorem}[Bhargava, 2022]\label{thm:A3}
    We have $M(4)\le 2760$ (and $M(4)\le 182$ if $|D(\mathcal{O})|$ is sufficiently large).
\end{theorem}

Bhargava proved his theorem via a parametrization of quartic rings and their cubic resolvent rings, and utilized Akhtari's recent upper bound (see the Appendix of Bhargava (2022)) for the number of solutions of quartic Thue equations.

Akhtari (2022) gave another, more direct proof for Theorem \ref{thm:A3}, following the approach of Ga\'{a}l, Peth\H{o} and Pohst (1996)  (which in fact is going into the same direction as Bhargava's approach but is less general), and combining this with her own upper bound
for the number of solutions of quartic Thue equations.

Theorem \ref{thm:A2} is probably far from best possible in terms of $n$. We pose the following problem:

\begin{prob}[Gy\H ory, 2000]
Do there exist absolute constants $c_1,c_2$ such that $M(n)< c_1n^{c_2}$ for all $n\ge 4$ ?    
\end{prob}

The best lower bound we could find
is due to Miller-Sims and Robertson (2005).
Let $p$ be a prime number, $\zeta_p$ a primitive $p$-th root of unity, and $K_p$ the associated real cyclotomic field, i.e., $\Qq (\zeta_p +\zeta_p^{-1})$.
Then $K_p$ has degree $(p-1)/2$, and its ring of integers is
$\OO_p:=
\Zz [\zeta_p+\zeta_p^{-1}]$.
They proved that if $p\geq 7$ then $\Zz [\alpha ]=\OO_p$ is satisfied
by $\alpha =\zeta_p^k+\zeta_p^{-k}$, $(\zeta_p^k+\zeta_p^{-k}+b)^{-1}$ ($b=-1,0,1,2$,
$k=1\kdots (p-1)/2$).
If $p=7$ then among these numbers there are precisely nine pairwise
$\Zz$-inequivalent ones and by a result of Ga\'{a}l and Schulte (1989)
these are up to $\Zz$-equivalence the only numbers $\alpha$ 
with $\Zz [\alpha ]=\OO_7$. If $p\geq 11$
then all these numbers are pairwise $\Zz$-inequivalent
and thus, the order $\OO_p$ is $5(p-1)/2=5[K_p:\Qq ]$ times
monogenic.

We now fix a number field $K$ of degree $\geq 3$, and consider only orders of $K$. As it turned out, most orders of $K$ have only small multiplicity of monogenicity, bounded above independently even of the degree of $K$.
In 2013, we proved the following result with B\'erczes:

\begin{theorem}[B\'erczes, Evertse and Gy\H ory, 2013]\label{thm:A4}
    Let $K$ be an algebraic number field of degree $\ge 3$. Then $K$ has only finitely many orders that are three times monogenic.
\end{theorem}

To see that this is optimal, let $K$ be a non-CM number field of degree $\geq 3$. Then the ring of integers of $K$ has infinitely many units $\varepsilon$ with $K=\Qq (\varepsilon )$. For every of these $\varepsilon$ we obtain  a two times monogenic order $\Zz [\varepsilon ] =\Zz [\varepsilon^{-1}]$ of $K$.

Theorem \ref{thm:A4} is proved
by means of a reduction to unit equations in more than two unknowns, and a use of ineffective finiteness theorems for these equations. So Theorem \ref{thm:A4} 
is ineffective, in the sense that its proof does not allow to determine the exceptional orders.

\begin{prob}
    Make Theorem \ref{thm:A4} effective.
\end{prob}

This seems to be completely out of reach. At present, it is not known how to make the results on unit equations in more than two unknowns effective.

\subsection{Outlines of the proofs of 
Theorems \ref{thm:A2} and \ref{thm:A4}}\label{sec:10.2}
\vphantom{a}

\noindent
We start with recalling some auxiliary results
from the literature.

\begin{theorem}[Beukers and Schlickewei, 1996]\label{thm:10.2}
Let $F$ be a field of characteristic $0$, and $\Gamma$
a multiplicative subgroup of $F^*\times F^*$ of rank $r$. Then the equation $x+y=1$ has at most $2^{8r+8}$ solutions $(x,y)\in\Gamma$.
\end{theorem}

\begin{corollary}\label{thm:10.2.1}
Let $F$ be a field of characteristic $0$, let
$m\geq 1$, and let $\Gamma$ be a multiplicative subgroup of $(F^*)^{2m}$ of rank $r$. Then there are at most $2^{8(r+2m-1)}$ tuples $(x_1,y_1\kdots x_m,y_m)\in\Gamma$ satisfying
\begin{equation}\label{eq:10.2.1}
x_i+y_i=1\ \ \text{for } i=1\kdots m.
\end{equation}
\end{corollary}

This result is easily deduced from Theorem \ref{thm:10.2}
using induction on $m$, see Evertse (2011), or Evertse and Gy\H{o}ry (2017), Corollary 4.3.5.

\begin{theorem}\label{thm:10.2.2}
Let $F$ be a field of characteristic $0$, let $m\geq 1$, and let $\Gamma$ be a multiplicative subgroup of $(F^*)^m$. Then there are at most finitely many tuples $(x_1\kdots x_m)\in\Gamma$ satisfying
\begin{equation}\label{eq:10.2.2}
\left\{ 
\begin{array}{l}
x_1+\cdots +x_m=1,
\\
x_{i_1}+\cdots +x_{i_t}\not= 0\ \text{for each non-empty subset $\{ i_1\kdots i_t\}$ of $\{ 1\kdots m\}$.}
\end{array}\right.
\end{equation}
\end{theorem}

This was proved by Evertse (1984b) and van der Poorten and Schlickewei (1982, 1991), combining Schmidt's and Schlickewei's Subspace Theorem from Diophantine approximation with a specialization argument. We note that Theorem \ref{thm:10.2.2} is ineffective, hence so are its consequences. Although we will not need these here, we mention that there are explicit upper bounds for the number of solutions of \eqref{eq:10.2.2} depending only on $m$ and on $r:=\rank\Gamma$, see Evertse, Schlickewei and Schmidt (2002) or Amoroso and Viada (2009), who obtained the up to now best upper bound $(8m)^{4m^4(m+r+1)}$. 

\begin{theorem}\label{thm:10.2.3}
Let $F$ be a field of characteristic $0$, and $\Gamma$ a multiplicative subgroup of $F^*\times F^*$. Then there are only finitely many pairs $(a,b)\in F^*\times F^*$ such that $a+b=1$, and such that $ax+by=1$ has three solutions $(x,y)\in\Gamma$,
the pair $(1,1)$ included.
\end{theorem}

\begin{proof}[Idea of proof]
This is basically a result of Evertse, Gy\H{o}ry, Stewart, and Tij\-de\-man (1988), see also Evertse and Gy\H{o}ry (2015), Theorem  6.1.6. The idea is as follows. Suppose that there are $(x_1,y_1), (x_2,y_2)\in\Gamma$, distinct from each other and distinct from $(1,1)$, such that $ax_i+by_i=1$ for $i=1,2$. Then
\[
\left|
\begin{array}{lll} 1&1&1\\ 1&x_1&y_1 \\ 1&x_2&y_2\end{array}
\right| =0.
\]
Expand the determinant, divide by a term to obtain a five term sum equal to $1$, 
consider all possible partitions into minimal vanishing subsums, and apply Theorem \ref{thm:10.2.2}
to each of them.
\end{proof}

Let $K$ be a number field of degree $n\geq 3$. Denote by $x\mapsto x^{(i)}$ ($i=1\kdots n$) the embeddings of $K$ in $G$, where $G$ is the normal closure of $K$. For $\alpha$ with $\Qq (\alpha )=K$ and $i=3\kdots n$, define
\[
x_i(\alpha )=
\frac{\alpha^{(i)}-\alpha^{(1)}}{\alpha^{(2)}-\alpha^{(1)}},\ \
y_i(\alpha )=
 \frac{\alpha^{(2)}-\alpha^{(i)}}{\alpha^{(2)}-\alpha^{(1)}}  
\]
and the tuple
\[
\kappa (\alpha ):=(x_3(\alpha ), y_3(\alpha )\kdots
x_n(\alpha ),y_n(\alpha )).
\]

In addition, we need a few simple lemmas.
We call $\alpha ,\beta$ $\Qq$-equivalent if $\beta =\lambda \alpha +a$ for some $\lambda\in\Qq^*$, $a\in\Qq$.

\begin{lemma}\label{lem:10.2.4}
Let $\alpha ,\beta$ with 
$\Qq (\alpha )=\Qq (\beta )=K$.
\\[0.1cm]
(i) $\kappa (\alpha )=\kappa (\beta )$ $\Longleftrightarrow$ $\alpha ,\beta$ are $\Qq$-equivalent.
\\[0.1cm]
(ii) Assume in addition that $\Zz [\alpha ]=\Zz [\beta ]$ and that $\alpha ,\beta$ are $\Qq$-equivalent. Then $\alpha ,\beta$ are $\Zz$-equivalent.
\end{lemma}

\begin{proof} (i) Clearly, $\kappa (\alpha )=\kappa (\beta )$
if and only if $(\alpha^{(i)},\beta^{(i)})$ ($i=1\kdots n$) are collinear, i.e., $\beta^{(i)}=\lambda \alpha^{(i)} +a$ ($i=1\kdots n$) for some $\lambda\in G^*$, $a\in G$. One easily shows that this is posible only if $\lambda , a$ are invariant under Galois action, i.e., lie in $\Qq$.

(ii) Our assumption $\Qq (\alpha )=\Qq (\beta )=K$
implies that $\beta =f(\alpha )$ for some unique polynomial $f\in\Qq [X]$ of degree $<n$, and then $\Zz [\alpha ]=\Zz [\beta ]$ implies $f\in \Zz [X]$.
So if $\alpha ,\beta$ are $\Qq$-equivalent, then $\beta =\lambda \alpha +a$ with $\lambda ,a\in\Zz$.
By interchanging the role of $\alpha ,\beta$ we see that $\lambda^{-1}\in\Zz$, hence $\lambda =\pm 1$.
\end{proof}

\begin{lemma}\label{lem:10.2.5}
Let $\alpha ,\beta$ be such that $\Qq (\alpha )=\Qq (\beta )=K$ and $\Zz [\alpha ]=\Zz [\beta ]$.
Then
\[
\frac{\beta^{(i)}-\beta^{(j)}}{\alpha^{(i)}-\alpha^{(j)}}\in\OO_G^*\ \text{for } i,j=1\kdots n,\, i\not =j.
\]
\end{lemma}

\begin{proof}
Use $\beta =f(\alpha )$, $\alpha =g(\beta )$ for some $f,g\in\Zz [X]$.
\end{proof} 

\begin{proof}[Sketch of the proof of Theorem \ref{thm:A2}]
Let $\OO$ be an order of $K$. Note that for $\alpha\in K$ with $K=\Qq (\alpha )$
we have relations
\[
x_i(\alpha)+y_i(\alpha )=1\ \ (i=3\kdots n)
\]
where
\[
x_i(\alpha )=
\frac{\alpha^{(i)}-\alpha^{(1)}}{\alpha^{(2)}-\alpha^{(1)}},\ \
y_i(\alpha )=
 \frac{\alpha^{(2)}-\alpha^{(i)}}{\alpha^{(2)}-\alpha^{(1)}}  
.
\]
It was proved in Evertse (2011),
see also Evertse and Gy\H{o}ry (2017), pages 206--208, that if one restricts to $\alpha$ with $\Zz [\alpha ]=\OO$,
then the set of tuples
\[
\{ \kappa (\alpha ):\, \Zz [\alpha ]=\OO\}
\]
generates a multiplicative subgroup of $(G^*)^{2n-4}$ of rank at most $n(n-1)/2$.
In the deduction of this we used a refinement of Lemma 
\ref{lem:10.2.5}.
Now an application of Corollary \ref{thm:10.2.1} and Lemma \ref{lem:10.2.4} implies Theorem \ref{thm:A2}.
\end{proof}

In the proof of Theorem \ref{thm:A4} we need the 
following lemma.
Call $\alpha_1$ $k$-special if $\alpha_1\in O_K$, $K=\Qq (\alpha_1)$ and there are $\alpha_2\kdots \alpha_k$ such that $\alpha_1\kdots\alpha_k$ are pairwise $\Zz$-inequivalent and $\Zz [\alpha_1]=\Zz [\alpha_2]=\cdots =\Zz [\alpha_k]$.

\begin{lemma}\label{lem:10.2.6}
Let $\CC$ be a $\Qq$-equivalence class of $2$-special numbers. Then $\CC$ is the union of finitely many $\Zz$-equivalence classes.
\end{lemma}

\begin{proof}
For the somewhat involved argument we refer to 
B\'{e}rczes, Evertse and Gy\H{o}ry (2013) or Evertse and Gy\H{o}ry (2017), Lemma 9.5.6.
\end{proof}

\begin{proof}[Sketch of the proof of Theorem \ref{thm:A4}]

We have to prove that there are only finitely many orders $\Zz [\alpha ]$ such that $\alpha$ is $3$-special. It suffices to show that the $3$-special $\alpha$ lie in finitely many $\Zz$-equivalence classes. We sketch the argument.

Let $\alpha\in O_K$ be $3$-special. Pick $\beta ,\gamma$ such that $\alpha ,\beta , \gamma$ are pairwise $\Zz$-inequivalent, and $\Zz [\alpha ]=\Zz [\beta ]=\Zz [\gamma ]$. For any three distinct indices $i,j,k$ from $\{ 1\kdots n\}$, define
\[
\varepsilon_{ijk}=
\frac{(\beta^{(i)}-\beta^{(j)})/(\beta^{(i)}-\beta^{(k)})}
{(\alpha^{(i)}-\alpha^{(j)})/(\alpha^{(i)}-\alpha^{(k)})},\ \ 
\eta_{ijk}=\frac{(\gamma^{(i)}-\gamma^{(j)})/(\gamma^{(i)}-\gamma^{(k)})}
{(\alpha^{(i)}-\alpha^{(j)})/(\alpha^{(i)}-\alpha^{(k)})}.
\]

Then by Lemma \ref{lem:10.2.5}, the equation
\[
\frac{\alpha^{(i)}-\alpha^{(j)}}{\alpha^{(i)}-\alpha^{(k)}}x+\frac{\alpha^{(j)}-\alpha^{(k)}}{\alpha^{(i)}-\alpha^{(k)}}y=1\ \ \text{in } x,y\in \OO_G^*
\]
has three solutions
\[
(1,1),\ (\varepsilon_{ijk},\varepsilon_{kji}),\ 
(\eta_{ijk},\eta_{kji}).
\]
If for all $i,j,k$ and all $\alpha ,\beta ,\gamma$ as above these three solutions were distinct, we could conclude  
from Theorem \ref{thm:10.2.3} that there is a finite set $\SS$ such that $\alpha^{(i,j,k)}\in\SS$ for all $i,j,k$ and all $3$-special $\alpha$. It need not be true, however, that in all cases these three solutions are distinct. However, by means of a combinatorial argument,
worked out in B\'{e}rczes, Evertse and Gy\H{o}ry (2013) or Evertse and Gy\H{o}ry (2017), pp. 211--216
we deduce that the existence of a finite set $\SS$ as above still holds. Now Lemma \ref{lem:10.2.4} (i) implies that the $3$-special numbers $\alpha$ lie in only finitely many $\Qq$-equivalence classes.                                                                    
Finally, Lemma \ref{lem:10.2.6} implies that the 
$3$-special numbers lie in only finitely many
$\Zz$-equivalence classes. 
\end{proof}

\subsection{Generalizations for rationally monogenic orders}
\label{sec:appB}
\vphantom{a}

\noindent
The theorems stated in Subsection \ref{sec:appA} have analogues for rationally monogenic orders. 
For the necessary terminology and properties we refer to Subsection \ref{sec:5.5}.

For a not necessarily integral algebraic number $\alpha$ of degree $n\geq 3$ we define
\begin{align*}
\MM_{\alpha} &:=\{ x_0+x_1\alpha+\cdots +x_{n-1}\alpha^{n-1}:\, x_0\kdots x_{n-1}\in\Zz\},
\\
\Zz_{\alpha} &:=\{ \xi\in\Qq (\alpha ):\, \xi\MM_{\alpha}\subseteq\MM_{\alpha}\}.
\end{align*}
Recall that $\Zz_{\alpha}=\Zz_{\beta}$ if $\alpha$ and $\beta$ are $GL_2(\Zz )$-equivalent.

An order $\OO$ of a number field $K$ is called \emph{rationally monogenic} if $\OO =\Zz_{\alpha}$ for some algebraic number $\alpha$. As observed in Subsection \ref{sec:5.5}, if $\alpha$ is an algebraic integer, then $\Zz_{\alpha}=\Zz [\alpha ]$.
Thus, monogenic orders are rationally monogenic. We further recall that rationally monogenic orders are primitive, i.e., they cannot be expressed as $\Zz +a\OO'$ for some integer $a>1$ and order $\OO'$.

We say that an order $\mathcal{O}$ of a number field $K$ is \textit{$k$ times/precisely $k$ times/at most $k$-times rationally monogenic} if up to $GL_2(\Zz )$-equivalence there are at least/precisely/at most $k$ numbers $\alpha$ such that $\mathcal{O}=\Zz_{\alpha}$.
Denote by $RM(n)$ the least number $k$ such that for every number field $K$ of degree $n$ and every order $\mathcal{O}$ of $K$, the order $\mathcal{O}$ is at most $k$ times rationally monogenic.

From work of Delone and Faddeev (1940) it follows that $RM(3)\leq 1$, that is, every order of a cubic number field 
is at most one time rationally monogenic (and in fact precisely one time if the order is primitive).
From a result of B\'{e}rczes, Evertse and Gy\H{o}ry (2004)
the following analogue of Theorem \ref{thm:A2} can be deduced:

\begin{theorem}\label{thm:B2}
For every $n\ge 4$,
$RM(n)$ is finite and in fact,
$RM(n)\leq n\times 2^{24n^3}$.
\end{theorem}

Similarly to Theorem \ref{thm:A2} the proof uses Theorem \ref{thm:10.2} of Beukers and Schlickewei (1996) mentioned above.

This bound has been improved. The best bounds to date are as follows:

\begin{theorem}\label{thm:B3}
We have
\begin{itemize}
\item[(i)] $RM(4)\leq 40$ (Bhargava (2022));
\item[(ii)] $RM(n)\leq 2^{5n^2}$ for $n\ge 5$ (Evertse and Gy\H{o}ry (2017)).
\end{itemize}
\end{theorem}

The proof of part (ii) is similar to that of Theorem \ref{thm:B2} but with a combinatorial improvement in the argument. The proof of part (i) also uses 
a parametrization of quartic rings and their cubic resolvent rings.

Recently, the following analogue of Theorem \ref{thm:A4} for rationally monogenic orders
was proved:

\begin{theorem}[Evertse, 2023]\label{thm:B4}\phantom{ }
\begin{itemize}
\item[(i)] Let $K$ be a number field of degree $4$.
Then $K$ has only finitely many three times rationally monogenic orders.
\item[(ii)] Let $K$ be a number field of degree $\geq 5$
such that the normal closure of $K$ is $5$-transitive. 
Then $K$ has only finitely many two times rationally monogenic orders.
\end{itemize}
\end{theorem}

Part (i) is best possible in the sense that there are quartic number fields having infinitely many two times rationally monogenic orders. 
In fact, B\'{e}rczes, Evertse and Gy\H{o}ry (2013, end of Section 1) give the following
construction:
\\[0.15cm]
Let $r,s$ be integers such that $f(X)=(X^2-r)^2-X-s$ is irreducible, and let $K=\Qq (\alpha )$,
where $\alpha$ is a root of $f$. Then $K$ has infinitely many orders $\OO_m$ ($m=1,2,\ldots$)
with the following property: $\OO_m=\Zz [\alpha_m]=\Zz [\beta_m]$, where $\beta_m=\alpha_m^2-r_m$,
$\alpha_m =\beta_m^2-s_m$ for some integers $r_m,s_m$.
\\[0.15cm]
It is clear that $\alpha_m,\beta_m$ in the above theorem are not $\GL_2(\Zz )$-equivalent.
We would like to pose the following problem:

\begin{prob}\label{prob:x1} Does every quartic number field have infinitely many orders that are two times rationally monogenic? If not, can we characterize those quartic number fields that do? Do the two times rationally monogenic orders have a particular structure?
\end{prob}

Similary to Theorem \ref{thm:A4}, Theorem \ref{thm:B4} has been proved
by means of a reduction to unit equations in more than two unknowns, and a use of ineffective finiteness theorems for such equations. So likewise, Theorem \ref{thm:B4} is ineffective.

It is not clear whether the $5$-transitivity condition on the Galois closure of $K$ in part (ii) is necessary; this was just a technical condition needed for the proof.
We are interested in the following problem:

\begin{prob}\label{prob:x2}
Is it true that every number field of degree $n\geq 5$ has only finitely many orders that are two times rationally monogenic? If not, can we characterize those number fields that do?
\end{prob}

Combining Theorems \ref{thm:B5} and  \ref{thm:B4} one can deduce the following counterpart of
Theorem \ref{thm:2.4}. For a number field $K$, let $\PP\II (K)$
denote the set of primitive, irreducible polynomials in $\Zz [X]$ having a root
$\alpha$ such that $K=\Qq (\alpha )$.

\begin{corollary}[Evertse, 2023]\label{cor:B6}\phantom{ }
\begin{itemize}
\item[(i)]
Let $K$ be a quartic number field. Then $\PP\II (K)$ has only finitely many Hermite equivalence classes that split into more than two $GL_2(\Zz )$-equivalence classes.
\item[(ii)]
Let $K$ be a number field of degree $\ge 5$ whose normal closure is $5$-transitive.
Then $\PP\II (K)$ has only finitely many Hermite equivalence classes that split into more than one $GL_2(\Zz )$-equivalence class.
\end{itemize}
\end{corollary}

Part (ii) was conjectured in BEGyRS (2023), without the $5$-transitivity condition.

\subsection{Outlines of the proofs of Theorems \ref{thm:B3} (ii) and \ref{thm:B4}}\label{sec:10.4}
\vphantom{a}

\noindent
The main new tool is the following result.

\begin{theorem}\label{thm:10.17}
Let $F$ be a field of characteristic $0$, and $\Gamma$ a finitely generated subgroup of $F^*$. Then there is a finite subset $\SS$ of $F^*$ with $1\in\SS$, such that for the set of solutions $(x_1,x_2,x_3,y_1,y_2,y_3)\in\Gamma^6$ of
\begin{equation}\label{eq:10.40}
(x_1-1)(x_2-1)(x_3-1)=(y_1-1)(y_2-1)(y_3-1)
\end{equation}
at least one of the following holds:
\\[0.2cm]
{\bf (i)} at least one of $x_1\kdots y_3$ belongs to $\mathcal{S}$;
\\
{\bf (ii)} there are $\eta_1,\eta_2,\eta_3\in\{ \pm 1\}$ such that
$(y_1,y_2,y_3)$ is a permutation of $(x_1^{\eta_1},x_2^{\eta_2},x_3^{\eta_3})$;
\\
{\bf (iii)} one of the numbers in
$\{x_ix_j,\, x_i/x_j,\, y_iy_j,\, y_i/y_j:\, 1\leq i<j\leq 3\}$
is equal to either $-1$, or to a primitive cube root of unity.
\end{theorem}

\begin{proof}
This is Proposition 8.1 of B\'{e}rczes, Evertse and Gy\H{o}ry (2013). The proof is basically to expand \eqref{eq:10.40}, divide by one term to get an equation of type \eqref{eq:10.2.2} in $16$ terms equal to $1$, consider all possible partitions into minimal vanishing subsums, and apply Theorem \ref{thm:10.2.2} to each of them (by using symmetric properties we can substantially reduce the number of cases).
\end{proof}

We need some other lemmas. Let $K$ be a number field of degree $n\geq 4$. Denote by $G$ the normal closure of $K$ and by $x\mapsto x^{(i)}$ ($i=1\kdots n$) the embeddings of $K$ in $G$. For $\alpha$ with $\Qq (\alpha )=K$ we define the \emph{cross ratios}
\[
\cross_{ijkl}(\alpha ):=\frac{(\alpha^{(i)}-\alpha^{(j)})(\alpha^{(k)}-\alpha^{(l)})}
 {(\alpha^{(i)}-\alpha^{(k)})(\alpha^{(j)}-\alpha^{(l)})}
\]
for any four distinct indices $i,j,k,l\in\{ 1\kdots n\}$ and we define the tuple
\[
\lambda (\alpha ):=(\cross_{123i}(\alpha ),\cross_{1i32}(\alpha ):\, i=4\kdots n).
\]
We call $\alpha ,\beta\in K$ $GL_2(\Qq )$-equivalent
if $\beta =\medfrac{a\alpha +b}{c\alpha +d}$ for some
$\big(\begin{smallmatrix} a&b\\c&d\end{smallmatrix}\big)\in GL_2(\Qq )$.

\begin{lemma}\label{lem:10.18}
Let $\alpha ,\beta$ with $\Qq (\alpha )=\Qq (\beta )=K$.
\\
(i) $\lambda (\alpha )=\lambda (\beta )$ 
$\Longleftrightarrow$ $\alpha,\beta$ are $GL_2(\Qq )$-equivalent.
\\
(ii) If $\Zz_{\alpha }=\Zz_{\beta}$ and $\alpha ,\beta$ are $GL_2(\Qq )$-equivalent, then $\alpha ,\beta$ are $GL_2(\Zz )$-equivalent.
\end{lemma}

\begin{proof} (i). $\Longleftarrow$ is straighforward. As for $\Longrightarrow$, from elementary projective geometry it follows that if $\lambda (\alpha )=\lambda (\beta )$ then there is a projective transformation of $\Pp^1$ defined over $G$ that maps $\alpha^{(i)}$ to $\beta^{(i)}$, for $i=1\kdots n$. It is easy to verify that this projective transformation is invariant under Galois action, hence defined over $\Qq$.

(ii). See for instance Lemma 2.6 of Evertse (2023).
\end{proof}

\begin{lemma}\label{lem:10.19}
Let $\alpha ,\beta$ with $\Qq(\alpha )=\Qq (\beta )=K$ and $\Zz_{\alpha }=\Zz_{\beta }$. Then for all distinct $i,j,k,l\in\{ 1\kdots n\}$ we have $\cross_{ijkl}(\beta )/\cross_{ijkl}(\alpha )\in \OO_G^*$.
\end{lemma}

\begin{proof} This is Lemma 2.4 of Evertse (2023).
\end{proof}

\begin{proof}[Sketch of the proof of Theorem
\ref{thm:B3} (ii)]
Let $\OO$ be an order of $K$. Note that for $\alpha\in K$ with $K=\Qq (\alpha )$
we have relations
\[
\cross_{123i}(\alpha)+\cross_{1i32}(\alpha )=1\ \ (i=4\kdots n).
\]
By Lemma 17.7.3 of Evertse and Gy\H{o}ry (2017), the set of tuples
\[
\{ \lambda (\alpha ):\, \Zz_{\alpha}=\OO\}
\]
generates a multiplicative subgroup of $(G^*)^{2n-6}$ of rank at most $n(n-1)/2$.
In the deduction of this we used a refinement of Lemma \ref{lem:10.19}.
Now an application of Corollary \ref{thm:10.2.1} and Lemma \ref{lem:10.18} implies Theorem \ref{thm:B3} (ii).
\end{proof}

Call $\alpha_1$ with $\Qq (\alpha_1)=K$ \emph{$k$-special} if there are $\alpha_2\kdots \alpha_k\in K$ such that $\alpha_1\kdots \alpha_k$ are pairwise $GL_2(\Zz )$-inequivalent and $\Zz_{\alpha_1}=\cdots =\Zz_{\alpha_k}$.
We should mention here that if $K$ has degree $3$ then there are no $2$-special numbers in $K$. 

In the proof of Theorem \ref{thm:B4}, we need the following lemma.

\begin{lemma}\label{lem:10.21}
Assume $n\geq 4$.
Let $\CC$ be a $GL_2(\Qq )$-equivalence class of $2$-special numbers. Then $\CC$ is the union of finitely many $GL_2(\Zz )$-equivalence classes.
\end{lemma}

\begin{proof}[Idea of proof] This is Proposition 5.1 of Evertse (2023). Its proof is fairly complicated. We give a brief outline.

We define $\cross_{ijkl}(\CC ):=\cross_{ijkl}(\alpha )$ for any $\alpha\in\CC$. This is well-defined since $GL_2(\Qq )$-equivalent algebraic numbers have the same cross ratios. Let $\alpha\in\CC$, let $\beta\in K$ be such that $\Zz_{\beta}=\Zz_{\alpha}$ and $\beta$ is not $GL_2(\Zz )$-equivalent to $\alpha$, and let $\DD$ be the $GL_2(\Qq )$-equivalence class of $\beta$. Then $\DD\not=\CC$ by Lemma \ref{lem:10.18} (ii). Clearly, $\cross_{ijkl}(\beta )=:\cross_{ijkl}(\DD )$
depends only on $\DD$. By Lemma \ref{lem:10.19} we have $\cross_{ijkl}(\DD )/\cross_{ijkl}(\CC )\in\OO_G^*$ for all $i,j,k,l$. Further,
\[
1=\cross_{ijkl}(\beta )+\cross_{ilkj}(\beta )
=\cross_{ijkl}(\CC )\cdot\frac{\cross_{ijkl}(\DD )}{\cross_{ijkl}(\CC )}+
\cross_{ilkj}(\CC )\cdot\frac{\cross_{ilkj}(\DD )}{\cross_{ilkj}(\CC )}
\]
for all $i,j,k,l$. Now by Theorem \ref{thm:10.2},
for given $\CC$ there are only finitely many possible values for each $\cross_{ijkl}(\DD )$ and thus,
by Lemma \ref{lem:10.18}, at most finitely many possibilities for $\DD$. It follows that $\CC$ is the union of finitely many sets
\[
\CC (\DD ):=\{ \alpha\in\CC :\, \text{there is $\beta\in\DD$ with $\Zz_{\alpha}=\Zz_{\beta}$}\}
\]
where $\DD\not=\CC$ is a $GL_2(\Qq )$-equivalence class of $2$-special numbers. So it suffices to prove that each set $\CC (\DD )$ is the union of finitely many $GL_2(\Zz )$-equivalence classes.

Now fix $\DD$, $\alpha\in\CC (\DD )$, and 
$\beta\in\DD$ such that $\Zz_{\alpha}=\Zz_{\beta}$.
Let $\alpha '$ be any other element of $\CC (\DD )$.
Then we can write
\[
\alpha '=\frac{a\alpha +b}{c\alpha +d}\ \text{with } a,b,c,d\in\Zz,\, \gcd (a,b,c,d)=1,\ ad-bc=:\Delta\not= 0.
\]
We have to prove that the numbers $\alpha'\in\CC (\DD )$ lie in only finitely many $GL_2(\Zz )$-equivalence classes.
Recall that there are $U\in GL_2(\Zz )$ and 
$a' ,b' ,d'\in\Zz$ with $a'd'=\Delta$, $|b'|\leq |d'|/2$ such that $U\big(\begin{smallmatrix}a&b\\ c&d\end{smallmatrix}\big)=\big(\begin{smallmatrix}a'&b'\\ 0&d'\end{smallmatrix}\big)$. Hence $\alpha '$ is $GL_2(\Zz )$-equivalent to $\alpha^*:= (a'\alpha +b')/d'$. It suffices to prove that there are only finitely many possibilities for $\alpha^*$.
It is in fact sufficient to prove that $\Delta$ is bounded, since for given $\Delta$ there are only finitely many possibilities for $(a',b',d')$. 
The boundedness of $\Delta$ is provided by the following elementary lemma, which is Proposition 4.1 of Evertse (2023). We refer to that paper for the rather lengthy proof.

\begin{lemma}\label{lem:10.22}
Let $D$ be the discriminant of $\Zz_{\alpha}$,
and let $\fa (\alpha ,\beta)$ be the ideal of $\OO_G$
generated by the numbers $\cross_{ijkl}(\beta )/\cross_{ijkl}(\alpha )-1$ ($1\leq i<j<k<l\leq n$).
Then $\Delta$ divides $D^5\cdot \fa (\alpha ,\beta )^2$.
\end{lemma}

\end{proof}

\begin{proof}[Sketch of the proof of Theorem \ref{thm:B4}]
For $\alpha ,\beta$ with $\Qq (\alpha )=\Qq (\beta )=K$ and distinct $i,j,k,l\in\{ 1\kdots n\}$ we put
\[
\ve_{ijkl}(\alpha ,\beta ):=\frac{\cross_{ijkl}(\beta )}{\cross_{ijkl}(\alpha )}.
\]
Lemma \ref{lem:10.19} implies that if $\Zz_{\alpha}=\Zz_{\beta}$, then $\ve_{ijkl}(\alpha ,\beta )\in \OO_G^*$.
\\[0.15cm]

\emph{The case $n=4$.} We have to show that there are only finitely many orders $\OO$ of $K$ such that $\OO=\Zz_{\alpha}$ for some $3$-special $\alpha$. It clearly suffices to show that the $3$-special numbers $\alpha$ lie in only finitely many $GL_2(\Zz )$-equivalence classes. 

Let $\alpha\in K$ be $3$-special, and choose $\beta ,\gamma\in K$ such that $\alpha ,\beta ,\gamma$ are pairwise $GL_2(\Zz )$-inequivalent, and $\Zz_{\alpha}=\Zz_{\beta}=\Zz_{\gamma}$. Let $(i,j,k,l)$ be a permutation of $(1,2,3,4)$. Then the equation
\[
\cross_{ijkl}(\alpha )x+\cross_{ilkj}(\alpha )y=1
\]
has three distinct solutions $(x,y)\in\OO_G^*\times\OO_G^*$,
i.e., $(1,1)$, $(\ve_{ijkl}(\alpha ,\beta ),\ve_{ilkj}(\alpha ,\beta ))$,
$(\ve_{ijkl}(\alpha ,\gamma ),\ve_{ilkj}(\alpha ,\gamma ))$. Now Theorem \ref{thm:10.2.3} implies that $\cross_{ijkl}(\alpha )$ can assume only finitely many values. From Lemma \ref{lem:10.19} it now follows that the $3$-special $\alpha\in K$ lie in only finitely $GL_2(\Qq )$-equivalence classes. Finally, from Lemma \ref{lem:10.21} it follows that they lie in finitely many $GL_2(\Zz )$-equivalence classes.
\\[0.15cm]

\emph{The case $n\geq 5$.} 
We have to show that there are only finitely many orders $\OO$ of $K$ such that $\OO=\Zz_{\alpha}$ for some $2$-special $\alpha$. It clearly suffices to show that the $2$-special numbers $\alpha$ lie in only finitely many $GL_2(\Zz )$-equivalence classes.

Let $\alpha\in K$ be $2$-special, and choose $\beta$
such that $\alpha ,\beta$ are $GL_2(\Zz )$-inequivalent and $\Zz_{\alpha}=\Zz_{\beta}$.
Henceforth, we write $\ve_{ijkl}$ for $\ve_{ijkl}(\alpha ,\beta )$.
Let $i,j,k,l$ be distinct indices from $\{ 1\kdots n\}$. Then
\[
\cross_{ijkl}(\alpha )+\cross_{ilkj}(\alpha )=1,\ \ 
\cross_{ijkl}(\alpha )\ve_{ijkl}+\cross_{ilkj}(\alpha )\ve_{ilkj}=1,
\]
$\ve_{ilkj}/\ve_{ijkl}=\ve_{iljk}$, which imply
\begin{equation}\label{eq:10.22}
\cross_{ijkl}(\alpha )=\frac{\ve_{ilkj}-1}{\ve_{ilkj}-\ve_{ijkl}},\ \ 
\cross_{ijkl}(\beta )=\ve_{ijkl}\cross_{ijkl}(\alpha )=\frac{\ve_{ilkj}-1}{\ve_{iljk}-1}.
\end{equation}
Now picking a fifth index $m$, 
and using 
$\medfrac{\cross_{jmlk}(\beta )\cross_{ijkm}(\beta )}{\cross_{ijkl}(\beta)}=1$,
we obtain
\begin{equation}\label{eq:10.23}
\frac{\ve_{jklm}-1}{\ve_{jkml}-1}\cdot
\frac{\ve_{imkj}-1}{\ve_{imjk}-1}\cdot\frac{\ve_{iljk}-1}{\ve_{ilkj}-1}=1.
\end{equation}
We apply Theorem \ref{thm:10.17} to \eqref{eq:10.23}
for all $i,j,k,l,m$. Our assumption that the Galois group of $G$ is $5$-transitive implies various conjugacy relations between the $\ve_{ijkl}$. 
Using all of these, we infer that for each quadruple $i,j,k,l$ there are only finitely many possible values for $\ve_{ijkl}$ (we should mention here that without the $5$-transitivity assumption, we do not know how to prove this). Now \eqref{eq:10.22} implies that there are only finitely many possible values for $\cross_{ijkl}(\alpha )$, if $\alpha$ runs through the $2$-special numbers of $\alpha$, and thus, by Lemma \ref{lem:10.18}, that the $2$-special $\alpha\in K$ lie in only finitely many $GL_2(\Qq )$-equivalence classes. Finally, from Lemma \ref{lem:10.21} it follows that they
lie in only finitely many $GL_2(\Zz )$-equivalence classes. 
\end{proof}

\addtocontents{toc}{\SkipTocEntry}
\section*{Appendix: Related topics}
\addcontentsline{toc}{section}{\numberline{}\vskip0.2cm {\sc Appendix: Related topics}}
\addtocontents{toc}{\vspace{0.1cm}}

\setcounter{section}{0}
\renewcommand{\thesection}{\Alph{section}}


We briefly discuss some further topics related to monogenic number fields and monogenic orders and generalizations thereof that do not strictly belong to the reduction theory of integral polynomials.

\section{Monogenicity, class group and Galois group}\label{sec:appNewC}
Recently, surprising results have been obtained in precise and quantitative form that imply that on average, the monogenicity of a number field has an altering effect on the structure of its $2$-class group, 
see 
Bhargava, Hanke and Shankar (2020), Siad (2021), Swaminathan (2023), Shankar, Siad and Swaminathan (2025), and
Bhargava, Shankar and Swaminathan (2025).
The $2$-class group $Cl_2(K)$ of a number field $K$ is the group of ideal classes of $K$ whose order divides $2$.

To illustrate this, we recall some results from the literature. A \emph{monogenized number field} is a pair $(K,\alpha )$ consisting of a number field $K$ and $\alpha\in\OO_K$ such that $\OO_K=\Zz [\alpha ]$.
Two monogenized number fields $(K_1,\alpha_1)$, $(K_2,\alpha_2)$ are called isomorphic if there are a field isomorphism $\varphi :\, K_1\to K_2$ and a rational integer $a$ such that $\alpha_2 =\pm\varphi (\alpha_1)+a$. 

We now restrict ourselves to monogenized cubic fields.
We define the height of a monogenized cubic field $(K,\alpha )$ as follows. Let $f=X^3+aX^2+bX+c\in\Zz [X]$ be the minimal polynomial of $\alpha$. Then the height of $(K,\alpha )$ is
\begin{align*}
H(K,\alpha ) := &\max (|I(f)|^3,J(f)^2/4),
\\
&\quad\text{where }
I(f):=a^2-3b,\ \ J(f):=-2a^3+9ab-27c .
\end{align*}
One can show that isomorphic monogenized cubic fields have the same height. Further, the pair $(I(f),J(f)^2)$ uniquely determines an isomorphism class. Lastly, the discriminant of $f$ is $D(f)=\medfrac{1}{27}\big( 4I(f)^3-J(f)^2\big)$.

\begin{theorem}[Bhargava, 2005]\label{thm:Bx1}
Let $K$ run through the cubic number fields, ordered by discriminant.
\\
(i) The average size of $Cl_2(K)$ over the totally real cubic  fields is $5/4$.
\\
(ii) The average size of $Cl_2(K)$ over the complex cubic fields is $3/2$.
\end{theorem}

\begin{theorem}[Bhargava, Hanke and Shankar, 2020]\label{thm:Bx2}
Let $(K,\alpha )$ run through the monogenized cubic number fields whose Galois closure has Galois group isomorphic to $S_3$, ordered by height.
\\
(i) The average size of $Cl_2(K)$ over the totally real monogenized cubic fields is $3/2$.
\\
(ii) The average size of $Cl_2(K)$ over the complex monogenized cubic fields is $2$.
\end{theorem}

Siad (2021) proved a generalization of the last theorem for number fields of odd degree $\geq 5$.

We briefly discuss some other topics.
Recently, Arpin, Bozlee, Herr and Smith (2023a,b) introduced and studied twisted monogenic relative extensions $K/L$.  They proved  that $L$ has trivial class group (this is the case if e.g. $L=\mathbb{Q}$) if and only if every twisted monogenic extension of $L$ is monogenic.

Another topic worth of study is the connection between (multiplicity of) monogenicity of the ring of integers of a number field $K$ and the size of the Galois group of its Galois closure.
The examples of number fields $K$ of degree $n = 3, 4, 5, 6$ in Section \ref{sec:5} show that the multiplicity of monogenicity of $\OO_K$ can be relatively large if the Galois group of the Galois closure of $K$ is $S_n$, i.e. if its size  is large relative to $n$.

\section{Distribution of monogenic and non-monogenic number fields}\label{sec:appC}
As is well-known, all quadratic number fields and cyclotomic fields are monogenic. For degree $n=3$, the first example of a non-monogenic number field was given by Dedekind (1878).
For every $n\ge 3$, there are infinitely many isomorphism classes of \textit{monogenic}, cf. Kedlaya (2012), and infinitely many isomorphism classes of \textit{non-monogenic}
number fields of degree  $n$. 

Let $K$ be a number field, and $\{ 1,\omega_2\kdots \omega_n\}$ a $\Zz$-module basis of $\OO_K$. Denote by
$I(X_2\kdots X_n)$ the associated index form, as introduced in Subsection \ref{sec:4.3}. Thus, if $\alpha =x_1+x_2\omega_2+\cdots +x_n\omega_n$ with $x_1\kdots x_n\in\Zz$, then $[\OO_K:\Zz [\alpha ]]=|I(x_2\kdots x_n)|$.
Consequently, $K$ is monogenic if and only if $I(x_2\kdots x_n)=\pm 1$ is solvable in $x_2\kdots x_n\in\Zz$.
We say that $K$ has \emph{no local obstruction to being monogenic} if for every prime number $p$, the equation $I(x_2\kdots x_n)=\pm 1$ has a solution $x_2\kdots x_n$ in the $p$-adic integers. This notion does not depend on the choice of $\omega_2\kdots\omega_n$. We recall some recent results.

\begin{theorem}[Alp\"{o}ge , Bhargava and Shnidman, 2025]\label{thm:B.x1}
Let $K$ run through the isomorphism classes of cubic fields, ordered by their absolute discriminant. Then a positive proportion of them are not monogenic, and yet have no local obstruction to being monogenic.
\end{theorem}

Subsequently, but published earlier, the same authors proved the following result for quartic fields. Recall that a number field is called rationally monogenic if its ring of integers is rationally monogenic.

\begin{theorem}[Alp\"{o}ge, Bhargava and Shnidman, 2024]\label{thm:B.x2}
Let $K$ run through the isomorphism classes of quartic fields, ordered by their absolute discriminant. Then a positive proportion of them are not \emph{rationally} monogenic, and  yet have no local obstruction to being monogenic.
\end{theorem}

    For $n=3,4,6$, tables of Ga\'al (2019) suggest that the density of monogenic number fields $K$ of degree $n$ decreases with the absolute value of the discriminant $|D_K|$.
 
Bhargava, Shankar and Wang established the following pioneering result.

\begin{theorem}[Bhargava, Shankar and Wang, 2022]
\label{thm:D1}
Denote by $M_n(X)$ the number of isomorphism classes of monogenic number fields $K$ of degree $n$ with $|D_K|\leq X$ and with associated Galois group $S_n$. Then
for every $n\ge 2$ we have
$$M_n(X)\gg X^{1/2+1/n}\ \ \text{as } X\to\infty.$$
\end{theorem}

The authors conjecture that the exponent on $X$ is optimal.

In Part II of their paper, the authors proved a corresponding result for rationally monogenic number fields:
 
 \begin{theorem}[Bhargava, Shankar and Wang, 2025]
\label{thm:D2}
Denote by $RM_n(X)$ the number of isomorphism classes of rationally monogenic number fields $K$ of degree $n$ with $|D_K|\leq X$ and with associated Galois group $S_n$. Then
for every $n\ge 3$ we have
$$RM_n(X)\gg X^{1/2+1/(n-1)}\ \ \text{as } X\to\infty.$$
\end{theorem}

Let $N_n(X)$ denote the number of isomorphism classes of number fields $X$ of degree $n$ with $|D_K|\leq X$. It is conjectured that $N_n(X)\asymp X$ as $X\to\infty$.
This is easy for $n=2$. Davenport and Heilbronn (1971) proved this conjecture for $n=3$
and Bhargava (2005, 2010) for $n=4,5$. 



\section{Arithmetic characterization of monogenic and multiply monogenic number fields}\label{sec:appD}
The following problem continues to attract considerable attention:
\\[0.15cm]
\textbf{Hasse's problem} (1960's): \textit{give an arithmetic characterization of monogenic number fields.}
\\

In this direction there are many important results for deciding the \textit{monogenicity} or \textit{non-monogenicity} of number fields from certain special infinite classes, including quadratic, cyclotomic, abelian, cyclic, pure, composite number fields, certain quartic, sextic, multiquadratic number fields and relative extensions, and parametric families of number fields defined by binomial, trinomial,\ldots\ irreducible polynomials.

In their proofs various types of tools are used, among others De\-de\-kind's criterion; Newton polygons; Montes' algorithm; Ore's theorem; Engstr\"om's theorem; Gr\"obner basis approach; reduction to binomial Thue equations; elliptic curve approaches, irreducible monic polynomials with square-free discriminant; non-squarefree discriminant approach; infinite parametric families of number fields; use of the index form equation approach with ``small" solutions.

For details, we refer to Dedekind (1878) and to the books Hensel (1908), Hasse (1963), Narkiewicz (1974), Evertse and Gy\H ory (2017), Ga\'al (2019) and the references given there. For some recent developments, see also the survey article Ga\'al (2024) 
with many interesting special results, and the recent interesting papers Kaur, Kumar and Remete (2025), Sharma and Sarma (2025), 
Gu\`{a}rdia and Perdet (2025), Ga\'{a}l (2025), Harrington and Jones (2025), Yakkou, Aghzer and Boua (202?),
and K\"{o}nig (2025).
We note that Hasse's problem has not yet been solved in full generality.

A more precise version of Hasse's problem is as follows.

\begin{prob}\label{prob:4}
For $m\ge 1$, give an arithmetic characterization of those number fields whose ring of integers is $m$ times monogenic.
\end{prob}

Clearly, Hasse's problem and Problem \ref{prob:4} do not properly belong to the reduction theory of integral polynomials.

Dedekind's necessary condition for monogenicity of a number field was generalized by Del Corso, Dvornicich and Simon (2005) to a condition for rational monogenicity. Perhaps this provides a tool to construct more examples of number fields that are not rationally monogenic.

\addtocontents{toc}{\SkipTocEntry}

\end{document}